\documentclass[12pt]{amsart}

\usepackage{amssymb}
\usepackage{epsf}

\newcommand{\DOT}{\setlength{\unitlength}{1pt}\begin{picture}(2.5,2)(1,1)
\put(1,2){\circle*{2}}\end{picture}}

\newcommand{\Edot}{{E_{\,\DOT}}}
\newcommand{\Epdot}{{E_{\,\DOT}'}}

\newcommand{\Etdot}{{\wt{E}_{\DOT}}}
\newcommand{\Etpdot}{{\wt{E}_{\DOT}'}}

\newcommand{\Fdot}{{F\!_{\DOT}}}

\newcommand{\Gdot}{{G_{\DOT}}}
\newcommand{\Gpdot}{{G_{\DOT}'}}

\newcommand{\ol}[1]{{\overline{#1}}}
\newcommand{\wh}[1]{{\widehat{#1}}}
\newcommand{\wt}[1]{{\widetilde{#1}}}

\newcommand{\precdot}{{\prec\!\!\!\cdot\,}}

\newcommand{\Len}{{\mathcal L}}

\newcommand{\Span}[1]{{\langle #1 \rangle}}
\newcommand{\Perm}[1]{{\langle #1 ]}}
\newcommand{\supp}{{\mbox{\rm supp}}}

\newcommand{\Sln}{{\mbox{\it SL}_n}}
\newcommand{\Spn}{{\mbox{\it Sp}_{2n}{\mathbb C}}}
\newcommand{\Son}{{\mbox{\it So}_{2n+1}{\mathbb C}}}

\newcommand{\tridot}%
{{\,\vartriangleleft\hspace{-8.5pt}\cdot\hspace{3.5pt}}}
\newcommand{\tridoto}%
{{\,\vartriangleleft\hspace{-8.5pt}\cdot\hspace{2pt}_0\,\,}}
\newcommand{\tridotk}%
{{\,\vartriangleleft\hspace{-8.5pt}\cdot\hspace{2pt}_k\,\,}}

\font\parenfont=cmr7
\font\subfont=cmr5

\newcommand{\pcdot}%
{{\,\prec\hspace{-6.5pt}\cdot\hspace{-2.pt}\raisebox{1.26pt}{\parenfont )}\,}}
\newcommand{\prcc}{{\,\prec\hspace{-5.2pt}\raisebox{1.26pt}{\parenfont )}\,}} 
\newcommand{\prcs}{{\prec\hspace{-1.5pt}\raisebox{.6pt}{\subfont )}\,}}

\newcommand{\QED}{
\setlength{\unitlength}{1.0pt}%
\begin{picture}(10,7.5)
\put(5,-2.5){\rule{2.5pt}{2.5pt}}
\put(2.5,0){\rule{5pt}{2.5pt}}
\put(0,2.5){\rule{10pt}{2.5pt}}
\end{picture}\vspace{10pt}}

\newtheorem{thm}{Theorem}[section]
\newtheorem{prop}[thm]{Proposition}
\newtheorem{lem}[thm]{Lemma}
\newtheorem{cor}[thm]{Corollary}
\newtheorem{example}[thm]{Example}
\newtheorem{alg}[thm]{Algorithm}
\newtheorem{defth}[thm]{Definition-Theorem}

\newenvironment{ex}{\begin{example}\sl}{\end{example}}
\newtheorem{definition}[thm]{Definition}
\newenvironment{defn}{\begin{definition}\sl}{\end{definition}}
\newtheorem{remark}[thm]{Remark}
\newenvironment{rem}{\begin{remark}\rm}{\end{remark}}

\begin{document}

\title[Pieri-type formula]{A Pieri-type formula for isotropic flag manifolds} 

\author{Nantel Bergeron \and Frank Sottile}

\address{Department of Mathematics and Statistics\\
        York University\\
        Toronto, Ontario M3J 1P3\\
	CANADA}
\email[Nantel Bergeron]{bergeron@mathstat.yorku.ca}
\urladdr[Nantel Bergeron]{http://www.math.yorku.ca/bergeron}
\address{Department of Mathematics\\
        University of Wisconsin\\
        Van Vleck Hall\\
        480 Lincoln Drive\\
        Madison, Wisconsin 53706-1388\\
        USA}
\email[Frank Sottile]{sottile@math.wisc.edu}
\urladdr[Frank Sottile]{http://www.math.wisc.edu/\~{}sottile}
\date{2 July 1998}
\thanks{First author supported in part by NSERC and CRM grants}
\thanks{Second author supported in part by NSERC grant  OGP0170279}
\thanks{\ research at MSRI supported by NSF grant DMS-9701755}

\begin{abstract}
We give the formula for  multiplying a Schubert class
on an odd orthogonal or symplectic flag manifold by a special Schubert class
pulled back from a Grassmannian of maximal isotropic subspaces.
This is also the formula for
multiplying a type $B$ (respectively, type $C$) Schubert
polynomial by the Schur 
$P$-polynomial $p_m$ (respectively, the Schur $Q$-polynomial $q_m$).
Geometric constructions and intermediate results allow us to ultimately
deduce this from formulas for the classical flag
manifold. 
These intermediate results are concerned with the Bruhat order of the Coxeter
group ${\mathcal B}_\infty$,  
identities of the structure constants for the Schubert basis of cohomology,
and intersections of Schubert varieties.
We show these identities follow from the Pieri-type
formula, except some `hidden symmetries' of the
structure constants. 
Our analysis leads to a new partial order on the Coxeter
group ${\mathcal B}_\infty$ and formulas for many of these structure constants.
\end{abstract}


\maketitle

\mbox{ \ }\vspace{-20pt}

\tableofcontents

\section*{Introduction}
The cohomology of a flag manifold $G/B$ has an integral basis of Schubert
classes ${\mathfrak S}_w$ indexed by elements $w$ of the Weyl group of $G$.
The algebraic structure of these rings is known~\cite{Borel_cohomology} with
respect to a monomial basis, and there are methods (Schubert polynomials)
for expressing the  ${\mathfrak S}_w$ in terms of this
basis~\cite{BGG,BH,De74,FK_Bn,Fu96,LS82a,PrRa97}.
Moreover, their multiplicative structure with respect to the Schubert basis
is determined by Chevalley's formula~\cite{Chevalley91}.
Despite this, it remains an open problem to give a closed or bijective
formula for the integral structure constants 
$c^w_{u\,v}$ defined by the
identity 
$$
{\mathfrak S}_u\cdot{\mathfrak S}_v \ =\ 
\sum_w c^w_{u\,v}\,{\mathfrak S}_w.
$$
These $c^w_{u\,v}$ are non-negative as they count the flags in a suitable
triple intersection of Schubert varieties.
They are expected to be related to the enumeration of chains in
the Bruhat order of the Weyl group (see~\cite{BS98a} and
the references therein).

Of particular interest are Pieri-type formulas which describe the constants
$c^w_{u\,v}$ when ${\mathfrak S}_v$ is a special Schubert class pulled back
from a Grassmannian projection ($G/P$, $P$ maximal parabolic), as these
determine the ring structure for the cohomology of $G/P$ when $P$ is any
parabolic subgroup.
When $G$ is $\Sln{\mathbb C}$, a Pieri-type formula for multiplication by a
special Schubert class was described~\cite{LS82a} in terms of the Weyl group
element $wu^{-1}$. 
A formula in terms of chains in the Bruhat order was conjectured~\cite{BB}
and given a geometric proof~\cite{Sottile96}. 
Our main results are the analogous formulas when  
$G$ is $\Spn$ or $\Son$ and ${\mathfrak S}_v$
is a special Schubert class pulled back from a Grassmannian of maximal
isotropic subspaces.
These are common generalizations of the Pieri-type formulas for 
$\Sln{\mathbb C}$, Chevalley's formula, and Pieri-type formulas for
Grassmannians of maximal isotropic
subspaces~\cite{Hiller_Boe}. 
\smallskip

Our proof uses results on the Bruhat order, identities
of these structure constants, a decomposition of intersections of Schubert
varieties, and formulas in the cohomology of the
$\Sln{\mathbb C}$-flag manifold to explicitly determine a
triple intersection of Schubert varieties.
This shows the coefficients in the Pieri-type formula are the intersection
number of a linear space with a collection of quadrics.
Some intermediate results, including a fundamental identity and some
additional `hidden symmetries' of 
the structure constants, are deduced from
constructions on $\Sln{\mathbb C}$-flag manifolds~\cite{BS98a}.
This analysis leads to other results,
including a new partial order on the infinite Coxeter group 
${\mathcal B}_\infty$
and a monoid for chains in this order as in~\cite{BS_monoid}. 
We show how the Pieri-type formula implies our fundamental identity, use the
identities to express many structure constants in terms of the
Littlewood-Richardson coefficients for the multiplication of Schur $P$-- (or
$Q$--) functions~\cite{Stembridge_shifted}, and apply the hidden symmetries 
to the enumeration of chains in the Bruhat order.

\section{Statement of results}

Schubert classes in the cohomology of the flag manifolds 
$\Son/B$ and $\Spn/B$ form integral bases indexed by
elements of the Weyl group ${\mathcal B}_n$.
We represent ${\mathcal B}_n$ as the group
of permutations $w$ of $\{-n,\ldots,-2,-1,1,\ldots,n\}$
satisfying $w(-a)=-w(a)$ for $1\leq a\leq n$.
Let ${\mathfrak B}_w$ denote the Schubert class indexed by 
$w\in{\mathcal B}_n$  in 
$H^*\Son/B$  and ${\mathfrak C}_w$ that
in $H^*\Spn/B$.
The degree of these classes is $2\cdot\ell(w)$, where the {\em length}
$\ell(w)$ of $w$ is
$$
\#\{0< i<j\leq n\mid w(i)>w(j)\}\quad + \quad
\sum_{i>0>w(i)}|w(i)|.
$$
For an integer $i$, let $\ol{\imath}$ denote $-i$.
For each $1\leq m\leq n$, 
define $v_m\in {\mathcal B}_n$ by 
$$
\ol{m}=v_m(1)<0<v_m(2)<\cdots<v_m(n).
$$
This indexes a (maximal isotropic) special Schubert class in either cohomology
ring, written as $p_m:={\mathfrak B}_{v_m}$ and  
$q_m:={\mathfrak C}_{v_m}$.
We state the Pieri-type formula for the products
${\mathfrak B}_w \cdot p_m$ and ${\mathfrak C}_w \cdot q_m$ in terms of
chains in the Bruhat order on ${\mathcal B}_n$.
For this, we need a definition.

\begin{defn}\label{def:0-Bruhat}
The {\em $0$-Bruhat order} $\leq_0$ on ${\mathcal B}_n$ is
defined recursively as follows:
$u\lessdot_0 w$ is a cover in the $0$-Bruhat order if and only if
\begin{enumerate}
\item[(1)] $\ell(u)+1 = \ell(w)$, and
\item[(2)] $u^{-1}w$ is a reflection of the form $(\ol{\imath},\,i)$ or 
      $(\ol{\imath},\,j)(\ol{\jmath},\,i)$ for some $0<i<j\leq n$.
\end{enumerate}
\end{defn}

Chevalley's formula~\cite{Chevalley91} may be stated as follows:
\begin{equation}\label{eq:Chevalley}
\begin{array}{rcl}
{\mathfrak B}_u\cdot p_1 &=& 
	{\displaystyle \sum_{u\lessdot_0 w} {\mathfrak B}_w}\\
{\mathfrak C}_u\cdot q_1 &=&   \rule{0pt}{19pt}
{\displaystyle
\sum_{u\lessdot_0 w} \chi(u^{-1}w){\mathfrak C}_w},
\end{array}
\end{equation}
where $\chi(u^{-1}w)$ is the number of transpositions in the reflection
$u^{-1}w$.

We enrich the 0-Bruhat order in two complementary ways.
Write the two types of covers in the 0-Bruhat order as
$u\lessdot_0(\ol{\beta},\,\beta)u$ and 
$u\lessdot_0(\ol{\beta},\,\ol{\alpha})(\alpha,\,\beta)u$
where $0<\alpha<\beta\leq n$.
The {\em labeled $0$-Bruhat r\'eseau} is a labeled directed multigraph with
vertex set ${\mathcal B}_n$ and labeled edges between covers in the 
0-Bruhat order given by the following rule:
If $u\lessdot_0(\ol{\beta},\,\beta)u$, then a single edge is drawn
with label $\beta$.
If $u\lessdot_0(\ol{\beta},\,\ol{\alpha})(\alpha,\,\beta)u$, then
two edges are drawn with respective labels $\ol{\alpha}$ and $\beta$.
Thus if $u\lessdot_0w$, then $\chi(u^{-1}w)$ counts the edges from $u$ to
$w$ in this  0-Bruhat r\'eseau.
The {\it labeled $0$-Bruhat order} is obtained from this r\'eseau by
removing edges with negative integer labels.

Given a (saturated) chain $\gamma$ in either of these structures, let 
${\rm end}(\gamma)$ denote the endpoint of $\gamma$.
A {\it peak} in a chain $\gamma$ is an index $i\in\{2,\ldots,m{-}1\}$
with $a_{i-1}<a_i>a_{i+1}$, where $a_1,a_2,\ldots, a_m$ is the sequence of
edge labels in $\gamma$.  A {\it descent} is an index $i<m$ with
$a_i>a_{i+1}$ and an {\it ascent} is an index $i<m$ with $a_i<a_{i+1}$.
\medskip

\noindent{\bf Theorem A. }(Pieri-type formula)
{\em 
Let $u\in {\mathcal B}_n$ and $m>0$.
Then
\begin{enumerate}
\item[I.] {\rm (Odd-orthogonal Pieri-type formula) }
$$
{\mathfrak B}_u\cdot p_m\ =\ \sum {\mathfrak B}_{{\rm end}(\gamma)},
$$
the sum over all chains $\gamma$ in the labeled $0$-Bruhat 
order of ${\mathcal B}_n$
which begin at $u$, have length $m$, and no peaks.
\item[II.] {\rm (Symplectic Pieri-type formula)}
$$
{\mathfrak C}_u\cdot q_m\ =\ \sum {\mathfrak C}_{{\rm end}(\gamma)},
$$
\begin{enumerate}
\item the sum over all chains $\gamma$ in the labeled $0$-Bruhat 
      r\'eseau of ${\mathcal B}_n$
      which begin at $u$, have length $m$, and no descents.
\item  
      the same sum, except with no ascents.
\end{enumerate}
\end{enumerate}
}\medskip

This generalizes Chevalley's formula and the Pieri-type formulas for 
$\Sln{\mathbb C}/B$,
which are expressed in~\cite{Sottile96} 
as a sum of certain labeled chains in the Bruhat order on the symmetric group
${\mathcal S}_n$ with no ascents/no descents.
The duality of these two formulas, one in terms of peaks for an order, and
the other in terms of descents/ascents for an enriched structure on that
order has connections with other dualities in combinatorics.
These include Fomin's duality of graded
graphs~\cite{Fomin_dualityI,Fomin_dualityII} and Stembridge's theory of
enriched $P$-partitions~\cite{stembridge_enriched}, where peak and descent
sets play a complementary role.
These relations are explored in~\cite{BMSvW}, which extends the theory
developed in~\cite{BS_hopf} to the ordered structures of this manuscript.

\begin{ex}\label{ex:pieri}
Represent permutations $w\in{\mathcal B}_3$ by their values
$w(1)\,w(2)\,w(3)$.
Consider the products ${\mathfrak B}_{3\ol{1}2}\cdot p_2$ and
${\mathfrak C}_{3\ol{1}2}\cdot q_2$.
Figure~\ref{fig:height2} shows the part of the 0-Bruhat r\'eseau of
\begin{figure}[htb]
 $$\epsfxsize=3.in \epsfbox{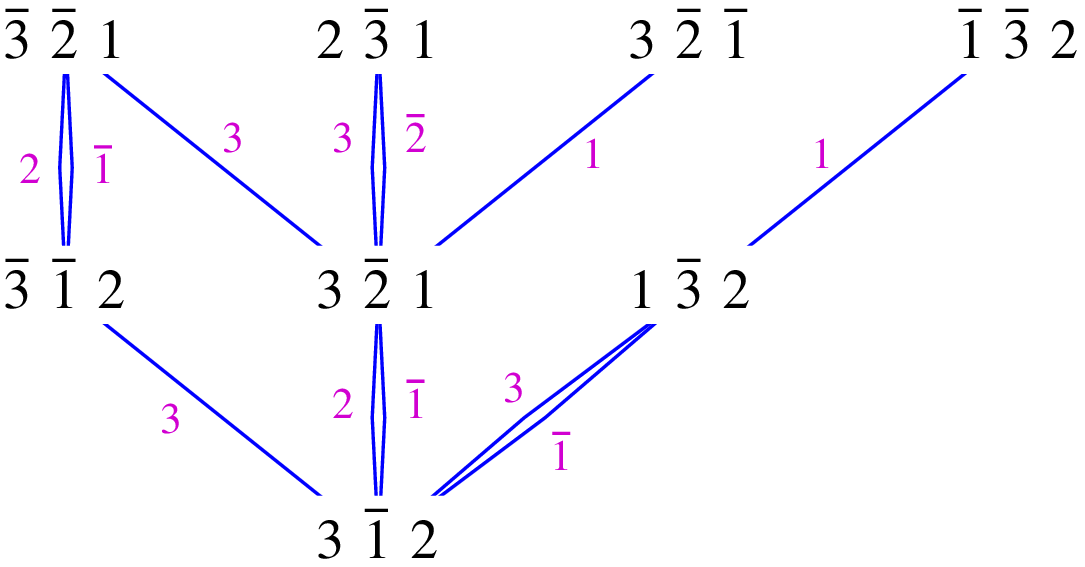}$$
 \caption{Chains above $3\ol{1}2$\label{fig:height2}}
\end{figure}
height 2 above $3\ol{1}2$ in  ${\mathcal B}_3$.
(Erase edges with negative labels to obtain its analog in the 0-Bruhat
order.) 
Chains of length 2 are peakless, so by Theorem~A I, we have
$$
{\mathfrak B}_{3\ol{1}2}\cdot p_2\ =\ 
2\,{\mathfrak B}_{\ol{3}\,\ol{2}1}+{\mathfrak B}_{2\ol{3}1}+
   {\mathfrak B}_{3\ol{2}\,\ol{1}}+{\mathfrak B}_{\ol{1}\,\ol{3}2}.
$$

Every chain in Figure~\ref{fig:height2} with increasing labels may be paired
with a chain with decreasing labels having the {\it same} underlying
permutations, and this pairing exhausts all chains.
Thus, by Theorem~A II, we have
$$
{\mathfrak C}_{3\ol{1}2}\cdot q_2\ =\ 
2\,{\mathfrak C}_{\ol{3}\,\ol{2}1}+2\,{\mathfrak C}_{2\ol{3}1}+
   {\mathfrak C}_{3\ol{2}\,\ol{1}}+{\mathfrak C}_{\ol{1}\,\ol{3}2}.
$$
\end{ex}\medskip

If $\lambda$ is a {\it strict partition} (decreasing integral sequence
$n\geq \lambda_1>\lambda_2>\cdots>\lambda_k>0$), then
$\lambda$ determines a unique 
{\em Grassmannian permutation} $v(\lambda)\in{\mathcal B}_n$ where
$v(i)=\ol{\lambda_i}$ for $i\leq k$ and $0<v(k+1)<\cdots<v(n)$.
If $k=1$ and $\lambda_1=m$, then $v_m=v(\lambda)$.
The Schubert classes $P_\lambda:={\mathfrak B}_{v(\lambda)}$ and 
$Q_\lambda:={\mathfrak C}_{v(\lambda)}$ are pullbacks of Schubert classes
from the Grassmannians of maximal isotropic subspaces 
$\Son/P_0$ and 
$\Spn/P_0$, where $P_0$ is the maximal parabolic
associated to the simple root of exceptional length.

Formulas for products of these $P$- and $Q$-classes are
known~\cite{Stembridge_shifted} as these classes  
are specializations of Schur $P$- and
$Q$-functions~\cite{Jozefiak,Pragacz_S-Q}.
Our proof of Theorem~A uses identities among the
structure constants $b^w_{u\,\lambda}$ and $c^w_{u\,\lambda}$ defined by the
following formulas.
$$
{\mathfrak B}_u\cdot P_\lambda\ =\ \sum_w b^w_{u\,\lambda}\ {\mathfrak B}_w
\qquad\mbox{and}\qquad
{\mathfrak C}_u\cdot Q_\lambda\ =\ \sum_w c^w_{u\,\lambda}\ {\mathfrak C}_w
$$
If $u,w,v(\lambda)\in{\mathcal B}_n$, then these constants do not depend
upon $n$.

Iterating Chevalley's formula~(\ref{eq:Chevalley}) shows that if either of
$b^w_{u\,\lambda}$ or $c^w_{u\,\lambda}$ is non-zero, then
$u<_0w$
and $\ell(w)-\ell(u)$ equals $|\lambda|$, the sum
of the parts of $\lambda$.
In fact the constant $b^w_{u\,\lambda}$ determines and is determined by the 
constant $c^w_{u\,\lambda}$:
Let $s(w)$ count the sign changes ($\{i\mid i>0>w(i)\}$) in $w$.
Then the map ${\mathfrak C}_w\mapsto 2^{s(w)}{\mathfrak B}_w$
embeds $H^*\Spn/B$ into  $H^*\Son/B$
and induces an isomorphism of their rational cohomology rings. 
Thus it suffices to work in $H^*\Spn/B$.
This is fortunate, as a key geometric result,
Theorem~\ref{thm:geom_facts}(2), holds only for $\Spn/B$.

Let $f^w_u$ count the saturated chains in the interval $[u,w]_0$ and
$g^w_u$ count the saturated chains in the r\'eseau $[u,w]_0$.
Iterating Chevalley's formula~(\ref{eq:Chevalley}) with $u=e$, the identity
permutation, we obtain the following expressions.
$$
p_1^m\ =\ \sum_{|\lambda|=m} f^{v(\lambda)}_e\, P_\lambda
\qquad\mbox{and}\qquad
q_1^m\ =\ \sum_{|\lambda|=m} g^{v(\lambda)}_e\, Q_\lambda
$$
Multiplying the first expression by ${\mathfrak B}_u$ and collecting the
coefficients of ${\mathfrak B}_w$ in the resulting expansion (likewise for
the second expression) gives the following proposition.

\begin{prop}\label{prop:chains}
Let $u,w\in{\mathcal B}_n$.
Then
$$
f^w_u\ =\ \sum_{|\lambda|=\ell(w)-\ell(u)}
f^{v(\lambda)}_e\,b^w_{u\,\lambda}
\qquad\mbox{and}\qquad
g^w_u\ =\ \sum_{|\lambda|=\ell(w)-\ell(u)}
g^{v(\lambda)}_e\,c^w_{u\,\lambda}.
$$
\end{prop}
\medskip

Theorem~A and Proposition~\ref{prop:chains}
show a close connection between chains in the
0-Bruhat order/r\'eseau and the structure constants
$b^w_{u\,\lambda}$ and $c^w_{u\,\lambda}$.
This justifies an elucidation of the basic properties of the $0$-Bruhat
order and r\'eseau, which we do in Sections 2 and 6.
These structures have a remarkable property and there are related 
fundamental identities among the structure constants.
\medskip

\noindent{\bf Theorem B. }
{\it
Suppose $u<_0w$ and  $x<_0z$ in ${\mathcal B}_n$ with 
$wu^{-1}=zx^{-1}$.
Then
\begin{enumerate}
\item[(1)] The map $v \mapsto vu^{-1}x$ induces an isomorphism of
  labeled intervals in 
  the 0-Bruhat order and 0-Bruhat r\'eseau 
  $[u,w]_0\stackrel{\sim}{\longrightarrow}[x,z]_0$.
\item[(2)]  For any strict partition $\lambda$, 
  $$
  b^w_{u\,\lambda}\ =\ b^z_{x\,\lambda}
  \qquad\mbox{and}\qquad
  c^w_{u\,\lambda}\ =\ c^z_{x\,\lambda}.
  $$
\end{enumerate}
}
\medskip

We prove Theorem~B(1) in Section 2.1 using combinatorial methods.
Theorem~B(2) is a consequence of a geometric result
(Theorem~\ref{thm:geometric_shapes}) proven in Section 4.
Both parts of Theorem~B are key to our proof of the
Pieri-type formula. 
Interestingly, the Pieri-type formula and Theorem~B(1) together imply
Theorem~B(2):

For any composition $\alpha=(\alpha_1,\ldots,\alpha_s)$ with each
$\alpha_i\geq 0$, let 
$p_\alpha:= p_{\alpha_1}\cdots p_{\alpha_s}$,
$q_\alpha:= q_{\alpha_1}\cdots q_{\alpha_s}$, and
$I(\alpha):=
\{\alpha_1,\alpha_1+\alpha_2,\ldots,\alpha_1+\cdots+\alpha_{s-1}\}$.
The {\it peak set} of a (maximal) chain in a labeled order is the set of
indices of peaks in the chain.
Given a chain in a labeled r\'eseau, its {\it descent set} (respectively
{\it ascent set}) is the set of indices of descents (respectively ascents)
in the chain.

\begin{cor}\label{cor:peak-set}
Let $u,w,x,z\in{\mathcal B}_n$.
\begin{enumerate}
\item[(1)] 
     Let $\alpha$ be any composition. Then the coefficient
     of\/ ${\mathfrak B}_w$ in the product ${\mathfrak B}_u\cdot p_\alpha$ is
     the number of chains in the interval $[u,w]_0$ in the labeled
     $0$-Bruhat order with peak set contained in $I(\alpha)$.     
\item[(2)]
     Let $\alpha$ be any composition. 
     Then the coefficient of\/ ${\mathfrak C}_w$ in
     the product ${\mathfrak C}_u\cdot q_\alpha$ is the number of 
     chains in the interval $[u,w]_0$ in the labeled $0$-Bruhat r\'eseau
     with descent set contained in  $I(\alpha)$.
     This is also the number with ascent set contained in $I(\alpha)$.   
\item[(3)]
     Suppose the Pieri-type formula (Theorem A) holds.
     Then the intervals $[u,w]_0$ and $[x,z]_0$ have the some number of
     chains with peak  set $I(\alpha)$ for every composition
     $\alpha$  if and only if  
     for every strict partition $\lambda$,
     $b^w_{u\,\lambda}=b^z_{x\,\lambda}$.
     The same statement holds for ascent/descent sets for chain in the
     r\'eseaux and the coefficients $c^w_{u\,\lambda},c^z_{x\,\lambda}$.
     In particular, Theorem~B(1) implies Theorem~B(2).
\end{enumerate}
Moreover, the numbers in 1 and 2 depend only upon the multiset
$\{\alpha_1,\ldots,\alpha_s\}$.
\end{cor}

Parts 1 and 2 follow from Theorem A.
For 3, note that the Schur $P$-polynomials (respectively $Q$-polynomials)
are linear combinations of the $p_\alpha$ (respectively the
$q_\alpha$)~\cite[III.8.6]{Macdonald95}. 
This linear combination gives a formula for 
$b^w_{u\,\lambda}$ (respectively $c^w_{u\,\lambda}$) in terms of chains with
given peak sets (respectively, given ascent/descent sets).

Let $\zeta\in{\mathcal B}_n$.
By Theorem~B(1), we may define 
$\eta\preceq\zeta$ if there is a 
$u\in{\mathcal B}_n$ with $u\leq_0\eta u\leq_0\zeta u$
and 
$\Len(\zeta):= \ell(\zeta u)-\ell(u)$
whenever $u\leq_0 \zeta u$.
Then $({\mathcal B}_n,\prec)$ is a graded partial order with rank function
$\Len(\cdot)$.
By the identity of Theorem~B(2), we may define 
$b^\zeta_\lambda:=b^{\zeta u}_{u\,\lambda}$ 
and $c^\zeta_\lambda:=c^{\zeta u}_{u\,\lambda}$ for any $u\in{\mathcal B}_n$
with $u\leq_0\zeta u$ and $|\lambda|=\Len(\zeta)$.

These coefficients satisfy one obvious identity, 
$c^\zeta_\lambda=c^{\zeta^{-1}}_\lambda$, as 
$c^w_{u\, v}=c^{\omega_0w}_{\omega_0u\, v}$
where $\omega_0\in{\mathcal B}_n$ is the longest element.
They also satisfy two others, which we call hidden symmetries.
Let $\rho\in{\mathcal B}_n$ be the permutation defined by
$\rho(i)=i-1-n$ for $1\leq i\leq n$. 
Then $\rho$ is the element with largest rank in $({\mathcal B}_n,\prec)$.
Let $\gamma\in{\mathcal B}_n$ be defined by 
$\gamma(1)=2,\gamma(2)=3,\ldots,\gamma(n)=1$, so that 
$\gamma=(\ol{1},\,\ol{2},\,\ldots,\ol{n})(1,\,2,\,\ldots,\,n)$.
\medskip

\noindent{\bf Theorem C. }
{\it 
For any $\zeta\in{\mathcal B}_n$,
\begin{enumerate}
\item[(1)] $\Len(\zeta)=\Len(\rho \zeta \rho)$ and for 
   any strict partition $\lambda$, we have 
   $b^\zeta_\lambda \ =\ b^{\rho\zeta\rho}_\lambda$ and 
   $c^\zeta_\lambda \ =\ c^{\rho\zeta\rho}_\lambda$.
\item[(2)] If $a\cdot \zeta(a)>0$ for all $a$, then 
   $\Len(\zeta)=\Len(\gamma\zeta\gamma^{-1})$ and for any 
   strict partition $\lambda$, we have 
   $b^\zeta_\lambda \ =\ b^{\gamma\zeta\gamma^{-1}}_\lambda$ 
   and 
   $c^\zeta_\lambda \ =\ c^{\gamma\zeta\gamma^{-1}}_\lambda$.
\end{enumerate}
}\medskip

We prove a strengthening of Theorem~B, relaxing the condition of
equality of $wu^{-1}$ and $zx^{-1}$ to that of {\it shape equivalence}.
Permutations $\eta,\zeta\in{\mathcal B}_n$ are shape equivalent if there
exist sets $I: 0<i_1<\cdots<i_s\leq n$ and 
$J: 0<j_1<\cdots<j_s\leq n$ such that $\eta$ acts as the identity on 
$\{1,\ldots,n\}\setminus I$, $\zeta$ acts as the identity on 
$\{1,\ldots,n\}\setminus J$, and $\eta(i_k)=i_l$ if and only if 
$\zeta(j_k)=j_l$.

Theorems B (the stronger version) and C allow us to determine many of the 
constants $b^w_{u\,\lambda}$ and $c^w_{u\,\lambda}$, showing they equal
certain Littlewood-Richardson coefficients $b^\kappa_{\mu\,\lambda}$ and
$c^\kappa_{\mu\,\lambda}$for Schur $P$- and $Q$-functions.
These are defined by the identities
$$
P_\mu\cdot P_\lambda\ =\ \sum_\kappa b^\kappa_{\mu\,\lambda}\; P_\kappa
\qquad\mbox{and}\qquad
Q_\mu\cdot Q_\lambda\ =\ \sum_\kappa c^\kappa_{\mu\,\lambda}\; Q_\kappa.
$$
A combinatorial formula for these coefficients was
given by Stembridge~\cite{Stembridge_shifted}.

\begin{defn}
Let $\mu,\kappa$ be strict partitions with $\mu\subset \kappa$.
We say that a permutation $\zeta\in{\mathcal B}_n$ has {\it skew shape}
$\kappa/\mu$ if 
\begin{enumerate}
\item[(1)] Either $\zeta$ or $\zeta^\rho$ is shape equivalent to 
      $v(\kappa)v(\mu)^{-1}$, or
\item[(2)] If $a\cdot \zeta(a)>0$ for all $a$, and one of
      $\zeta,\zeta^\gamma,\zeta^{\gamma^2},\ldots,\zeta^{\gamma^{n-1}}$ is
      shape equivalent to $v(\kappa)v(\mu)^{-1}$.
\end{enumerate}
\end{defn}

\begin{cor}
If $u\leq_0 w$ are permutations in ${\mathcal B}_n$ and 
$wu^{-1}$ has a skew shape $\kappa/\mu$, then for any
strict partition $\lambda$ we have
$$
b^w_{u\,\lambda}\ =\ b^\kappa_{\mu\,\lambda}
\qquad\mbox{and}\qquad
c^w_{u\,\lambda}\ =\ c^\kappa_{\mu\,\lambda}.
$$
\end{cor}

We call the partial order $\prec$ the {\em Lagrangian order} and
transfer the labeling from the $0$-Bruhat
order to obtain the labeled Lagrangian order.
In the same fashion, we may transfer the labeling and multiple edges of the
0-Bruhat r\'eseau to $({\mathcal B}_n,\prec)$, obtaining the (labeled) 
{\em Lagrangian r\'eseau}.
By Corollary~\ref{cor:peak-set}(3), Theorem~C has a purely enumerative
corollary.

\begin{cor}\label{cor:hidden-peaks}
For any  $\zeta\in{\mathcal B}_n$,
\begin{enumerate}
\item[(1)] For any subset $S$ of $\{2,\ldots, \Len(\zeta)-1\}$, the intervals 
      $[e,\zeta]_\prec$ and  $[e,\rho\zeta\rho]_\prec$ in the Lagrangian
      order have the same number of chains with peak set $S$. 
\item[(2)] For any subset $S$ of $\{1,\ldots, \Len(\zeta)-1\}$, the intervals
      $[e,\zeta]_\prec$ and  $[e,\rho\zeta\rho]_\prec$ 
      in the Lagrangian r\'eseau order have the same number of
      chains with descent set $S$ and the same number of chains with
      ascent set $S$, and these two numbers are equal.
\item[(3)] If $a\cdot \zeta(a)>0$ for all $a$, then the same is true for 
      $[e,\zeta]_\prec$ and $[e,\gamma\zeta\gamma^{-1}]_\prec$.  
\end{enumerate}
\end{cor}

In general, $[e,\zeta]_\prec\not\simeq [e,\rho\zeta\rho]_\prec$,
(See Figure~\ref{fig:cyclic} in Example~\ref{ex:cyclic})
and if $a\cdot \zeta(a)>0$ for all $a$, then 
in general, 
$[e,\zeta]_\prec\not\simeq [e,\gamma\zeta\gamma^{-1}]_\prec$.
(See Figure~\ref{fig:hidden} in Example~\ref{ex:no-iso}.)

Let $\supp(\zeta):= \{a>0\mid \zeta(a)\neq a\}$, the {\it support} of
$\zeta$. 
A permutation $\zeta\in {\mathcal B}_n$ is {\em reducible} if it 
has a non-trivial factorization
$\zeta=\eta\cdot\xi$ with $\Len(\zeta)=\Len(\eta)+\Len(\xi)$ where $\eta$
and $\xi$ have disjoint supports ($\eta\cdot\xi=\xi\cdot\eta$). 
If $\eta\cdot\xi=\xi\cdot\eta$ with 
$\Len(\xi\cdot\eta)=\Len(\eta)+\Len(\xi)$,
then $\eta\cdot\xi$ is a {\em disjoint product}.
Permutations 
$w\in {\mathcal B}_n$  have unique factorizations into irreducibles.

We represented ${\mathcal B}_n\subset{\mathcal S}_{\pm [n]}$, the group of
permutations of $\{-n,\ldots,-1,1,\ldots,n\}$. 
We also have ${\mathcal S}_n\hookrightarrow {\mathcal B}_n$, and the 
image consists of those $\zeta$ with $a\cdot\zeta(a)>0$ for every $a$:
For $\eta\in {\mathcal S}_n$, let 
$\overline{\eta}\in{\mathcal S}_{-[n]}$ be the permutation such that 
$\overline{\eta}(-i)=-\eta(i)$.
Then $\eta \ol{\eta}\in {\mathcal B}_n$.
For $\zeta\in{\mathcal B}_n$, define $\delta(\zeta)=1$ if $\zeta$ is in
the image of ${\mathcal S}_n$, and  $\delta(\zeta)=0$ otherwise.

In Section~\ref{sec:minimal} we establish the following result.

\begin{lem}
Let $\zeta\in {\mathcal B}_n$ and suppose $\zeta$ is irreducible.
Then $\Len(\zeta)\geq\# \supp(\zeta)-\delta(\zeta)$.
If $\Len(\zeta)=\# \supp(\zeta)-\delta(\zeta)$, then either
\begin{enumerate}
\item[(1)] 
     we have $\delta(\zeta)=1$ and there exists a cycle 
     $\eta\in {\mathcal S}_n$ with $\zeta=\eta \ol{\eta}$, or
\item[(2)]
     we have $\delta(\zeta)=0$ and $\zeta$ is a single cycle in 
     ${\mathcal S}_{\pm[n]}$.
\end{enumerate}
\end{lem}

\begin{defn}
If every irreducible factor $\eta$ of $\zeta$ satisfies 
$\Len(\eta)=\# \supp(\eta)-\delta(\eta)$, then we say that $\zeta$ is 
{\it minimal}.

If $\zeta\in{\mathcal B}_n$ is minimal, then set
\begin{eqnarray*}
\theta(\zeta)&:=& 
      2^{\#\{\mbox{\scriptsize irreducible factors of $\zeta$}\} -1},\\
\chi(\zeta)&:=& 
      2^{\#\{\mbox{\scriptsize irreducible factors $\eta$ of 
        $\zeta$ with $\delta(\eta)=1$}\}}.
\end{eqnarray*}
If $\zeta$ is not minimal, then set $\theta(\zeta)=\chi(\zeta)=0$.
\end{defn}

We state the Pieri-type formula in terms of the permutation $wu^{-1}$.
\medskip

\noindent{\bf Theorem D. }
{\em 
Let $u,w\in {\mathcal B}_n$ and $m\leq n$.
Then 
$$
{\mathfrak B}_u\cdot p_m\ =\ \sum \theta(wu^{-1}){\mathfrak B}_w
\qquad\mbox{and}\qquad
{\mathfrak C}_u\cdot q_m\ =\ \sum \chi(wu^{-1}){\mathfrak C}_w,
$$
the sum over all $w\in {\mathcal B}_n$ with $u\leq_0 w$ and 
$\ell(w)-\ell(u)=m$.
}\medskip

This is similar to the form of the Pieri-type
formula for $\Sln{\mathbb C}/B$ in~\cite{LS82a}, which is in
terms of the cycle structure of the permutation $wu^{-1}$.
This also generalizes the form of the the Pieri formula for 
Grassmannians of maximal isotropic subspaces~\cite{Hiller_Boe}.
Our proof for $Sp_{2n}{\mathbb C}/B$ 
shows these multiplicities to arise from the
intersection of a linear subspace of ${\mathbb P}^{2n-1}$ with a collection of
quadrics, one for each irreducible factor $\eta$ of $wu^{-1}$ with
$\delta(\eta)=1$,
similar to the proof of the Pieri-type formula for maximal isotropic
Grassmannians in~\cite{sottile_maximal}.
We relate the two formulations (Theorems A and D) of the Pieri-type formula
in sections~\ref{sec:lagr-or} and~\ref{sec:lagr-re}.

\begin{ex}
We return to Example~\ref{ex:pieri}.
For $u=3\ol{1}2$ and $w$ equal to each of\/ 
$\ol{3}\,\ol{2}1$, $2\ol{3}1$, $3\ol{2}\,\ol{1}$, and $\ol{1}\,\ol{3}2$ in
turn, $wu^{-1}$ is the permutation in ${\mathcal S}_{\pm[3]}$:
$$
(12)(\ol{1}\,\ol{2})\,(3\ol{3}),\quad
(132)(\ol{1}\,\ol{3}\,\ol{2}),\quad
(12\ol{1}\,\ol{2}),\quad \mbox{and}\quad
(13\ol{1}\,\ol{3}).
$$
As permutations in ${\mathcal B}_3$, the first has 2 irreducible factors,
for these we have $\delta((12)(\ol{1}\,\ol{2}))=1$ and $\delta((3\ol{3}))=0$.
The other permutations are irreducible with 
$\delta((132)(\ol{1}\,\ol{3}\,\ol{2}))=1$ and 
$\delta((12\ol{1}\,\ol{2}))=\delta((13\ol{1}\,\ol{3}))=0$.
Thus the values of $\theta$ are $2,1,1,1$ and of $\chi$ are $2,2,1,1$, which
shows the two forms of the Pieri-type formula agree on this example.

We remark that the last two permutations, $(12\ol{1}\,\ol{2})$ and 
$(13\ol{1}\,\ol{3})$, are shape equivalent.
\end{ex}

While we use the cohomology rings of complex varieties,
our results and methods are valid for the Chow
rings~\cite{Fulton_intersection} and $l$-adic ({\'e}tale)
cohomology~\cite{Deligne_SGA4.5} of these same varieties over
any field not of characteristic 2. 

This paper is organized as follows:
Section 2 contains basic combinatorial definitions and 
properties of the Bruhat order on
${\mathcal B}_\infty$ analogous to those of the
symmetric group  established
in~\cite{BS98a,BS_monoid}.
Section 3 contains the basic geometric 
definitions.
In Section 4, we use geometry to establish the main identity, 
Theorem~B(2). 
In Section 5, we establish additional geometric results and prove Theorem~C.
In Section 6, we establish further combinatorial properties of the
Lagrangian order and r\'eseau needed for the proof of the Pieri-type formula,
which is given in Section 7.

\section{Orders on ${\mathcal B}_\infty$}\label{sec:orders}

We derive the basic properties of the 0-Bruhat order on ${\mathcal B}_\infty$
analogous to properties of the $k$-Bruhat order on
${\mathcal S}_\infty$.
Further properties are developed in Section 6. 

Let $\#S$ be the cardinality of a finite set $S$.
For an integer $j$, its absolute value is $|j|$ and
let $\ol{\jmath}:=-j$.
Likewise, for a set $P$ of integers, 
define $\ol{P}:=\{\ol{\jmath}\mid j\in P\}$
and $\pm P:= P\cup \ol{P}$.
Set $[n]:=\{1,\ldots,n\}$
and let  ${\mathcal S}_{\pm[n]}$ be the group of permutations of $\pm[n]$.
Let $e$ be the identity permutation in ${\mathcal S}_{\pm[n]}$ and 
$\omega_0$ the longest element in ${\mathcal S}_{\pm[n]}$:
$\omega_0(i)=\ol{\imath}$.
Then ${\mathcal B}_n$ is the subgroup of ${\mathcal S}_{\pm[n]}$ 
for which $\omega_0 w \omega_0 = w$ and 
$\omega_0\in {\mathcal B}_n$.
We also have ${\mathcal B}_n\subset{\mathcal S}_{[\ol{n},n]}$, the symmetric
group on $[\ol{n},n]:=\pm[n]\cup\{0\}$.
We refer to elements of these groups as permutations.
Permutations $w\in{\mathcal B}_n$ are often represented by
their values 
$w(1)\,w(2)\,\ldots\,w(n)$.
For example, $2\,4\,\ol{3}\,\ol{1}\in{\mathcal B}_4$.
The {\em length} $\ell(w)$ of $w\in {\mathcal B}_n$ is
$$
\ell(w)\quad =\quad \#\{0<i<j\mid w(i)>w(j)\}\ +\ 
\sum_{i>0> w(i)}|w(i)|.
$$
Thus $\ell(2\,4\,\ol{3}\,\ol{1})=4+4=8$.
Note that $\omega_0$ is the longest element in ${\mathcal B}_n$.

An important class of permutations are the {\em Grassmannian
permutations}, those $v\in {\mathcal B}_n$ for which
$v(1)<v(2)<\cdots<v(n)$.
Such a permutation is determined by its initial negative values.
If $v(k)<0<v(k+1)$, define $\lambda(v)$ to be the decreasing sequence
$\ol{v(1)}>\ol{v(2)}>\cdots>\ol{v(k)}$.
Note that $\ell(v)=\ol{v(1)}+\cdots+\ol{v(k)} =: |\lambda(v)|$.
Likewise, given a decreasing sequence $\mu$ of positive integers 
(a {\it strict partition}) with 
$n\geq\mu_1$, let $v(\mu)$ be the Grassmannian permutation with 
$\lambda(v(\mu))=\mu$.
We write $\mu\subset\lambda$ for strict partitions $\mu,\lambda$
if $\mu_i\leq\lambda_i$ for all $i$, equivalently, if
$v(\mu)\leq_0v(\lambda)$.

The inclusion $\pm[n]\hookrightarrow\pm[n{+}1]$ induces inclusions
${\mathcal B}_n\hookrightarrow{\mathcal B}_{n+1}$ and 
${\mathcal S}_{\pm[n]}\hookrightarrow{\mathcal S}_{\pm[n+1]}$.
Define ${\mathcal B}_\infty:= \bigcup_n{\mathcal B}_n$
and ${\mathcal S}_{\pm\infty}:= \bigcup_n{\mathcal S}_{\pm[n]}$.

\begin{figure}[htb]
$$\epsfxsize=3.5in \epsfbox{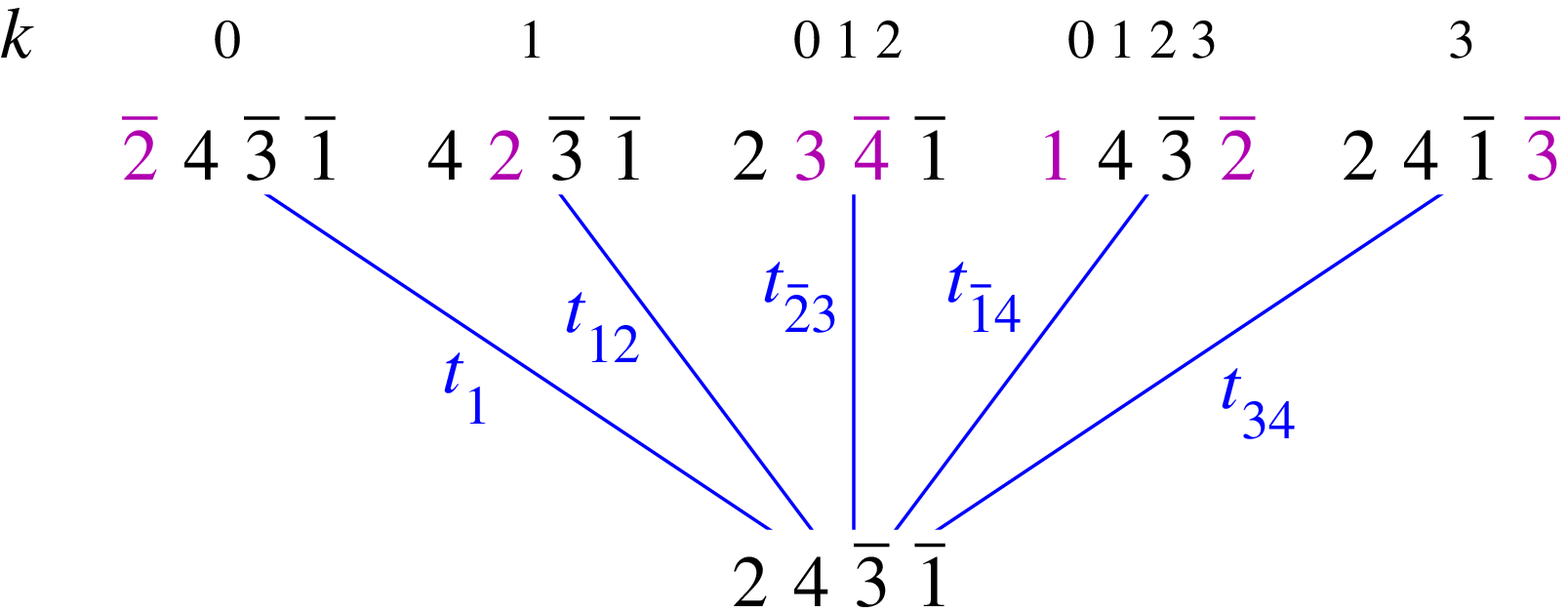}$$
\caption{Covers of $2\,4\,\ol{3}\,\ol{1}$\label{fig:covers}}
\end{figure}

In this representation, ${\mathcal B}_\infty$ has three types of
reflections, which are, as elements of ${\mathcal S}_{\pm\infty}$:
$$
\begin{array}{lcl}
t_{i\,j}&:=& (\ol{\jmath},\,\ol{\imath})(i,\,j)\\
t_{j}   &:=& (\ol{\jmath},\,j)\\
t_{\ol{\imath}\,j}
        &:=& (\ol{\jmath},\,i)(\ol{\imath},\,j)
\end{array}\qquad 
\mbox{for } 0<i<j.
$$
These reflections act on positions on the right and on values on the left.
The {\em Bruhat order} on ${\mathcal B}_\infty$ is defined by its covers:
$u\lessdot w$ if $\ell(u)+1=\ell(w)$ and $u^{-1}w$ is a reflection.
For each $k=0,1,\ldots$, define the {\em $k$-Bruhat order} 
(on ${\mathcal B}_n$ or ${\mathcal B}_\infty$) by its covers:
Set $u\lessdot_k w$ if $u\lessdot w$ and
$$
u^{-1}w\quad \mbox{is one of}\quad\left\{
\begin{array}{ll}
t_{i\,j}                 &\mbox{ with }i\leq k<j, \ \mbox{ \ or}\\
t_{j}                    &\mbox{ with }k<j, \ \mbox{ \ or}\\
t_{\ol{\imath}\,j} &\mbox{ with }k<j.
\end{array}\right.
$$
For example, Figure~\ref{fig:covers} shows all covers $w\in{\mathcal B}_4$
of $u=2\,4\,\ol{3}\,\ol{1}$,
the reflection $u^{-1}w$, and for which $k$ this is a cover in the
$k$-Bruhat order.

\subsection{The 0-Bruhat order}\label{sec:0-Bruhat}

While these orders are analogous to the $k$-Bruhat orders on 
${\mathcal S}_\infty$~\cite{BS98a,BS_monoid,LS83,Sottile96},
only the 0-Bruhat order on $B_\infty$ enjoys most 
properties of the $k$-Bruhat orders on ${\mathcal S}_\infty$.
This is because the 0-Bruhat order is an induced suborder of the
$0$-Bruhat order on ${\mathcal S}_{\pm\infty}$.

The {\em length} $l(w)$ of a permutation $w\in {\mathcal S}_{\pm\infty}$ 
counts the inversions of $w$:
$$
l(w)\ :=\ \#\{ i<j \mid w(i)>w(j)\}.
$$
The Bruhat order ($\vartriangleleft$) on ${\mathcal S}_{\pm\infty}$ is defined
by its covers:  
$u\tridot w$ if and only if 
$wu^{-1}$ is a transposition and $l(w)=l(u)+1$.
If $k\in{\mathbb Z}$, this is a cover (written $\tridotk$\!\!) in the k-Bruhat
order  ($\vartriangleleft_k$) on 
${\mathcal S}_{\pm\infty}$ (or ${\mathcal S}_{\pm[n]}$) if
$wu^{-1}=(a,\,b)$ with $a<k<b$.
The k-Bruhat order has a non-recursive characterization, needed below: 

\begin{prop}[\cite{BS98a}, Theorem~A]\label{prop:one}
Let $u,w\in {\mathcal S}_{\pm\infty}$ and $k\in{\mathbb Z}$. 
Then $u\vartriangleleft_k w$ if and only if
\begin{enumerate}
\item[(1)] $a<k < b$ implies $u(a)\leq w(a)$ and $u(b)\geq w(b)$.
\item[(2)] If $a<b$,  $u(a)<u(b)$,  and $w(a)>w(b)$, then $a<k<b$.
\end{enumerate}
\end{prop}

For the remainder of this section, we will be concerned with the case $k=0$.

\begin{thm}\label{thm:induced_order}
The 0-Bruhat order on ${\mathcal B}_\infty$ is the order induced from the
0-Bruhat order on ${\mathcal S}_{\pm\infty}$ by the inclusion
${\mathcal B}_\infty\hookrightarrow{\mathcal S}_{\pm\infty}$.
\end{thm}

\noindent{\bf Proof. }
For $u,w\in {\mathcal B}_\infty$, it is straightforward to verify
$$
u\lessdot_0 u t_j \ \Longleftrightarrow\ 
u\tridoto u(\ol{\jmath},\,j)
$$
and 
$$
u\lessdot_0 u t_{\ol{\imath}\,j}\ 
\Longleftrightarrow\ 
u\tridoto u(\ol{\jmath},\,i)
\tridoto u(\ol{\jmath},\,i)(\ol{\imath},\,j).
$$
Thus $u<_0 w\Rightarrow u\vartriangleleft_0 w$, and so $<_0$ is a suborder
of $\vartriangleleft_0$.

To show this suborder is induced, suppose that 
$u\vartriangleleft_0 w$ with $u,w\in {\mathcal B}_\infty$,
and argue by induction on $l(w)-l(u)$.
Suppose 
$u\tridoto v\vartriangleleft_0 w$.
If $v=u (\ol{\jmath},\,j) =u t_j$, then $v\in{\mathcal B}_n$ and we
are done by induction.

Suppose now that $v=u(\ol{\jmath},\,i)\not\in {\mathcal B}_n$.
By the involution $x\mapsto\omega_0 x\omega_0$ of 
$({\mathcal S}_{\pm\infty},\vartriangleleft_0)$, 
$\omega_0v\omega_0=u(\ol{\imath},\,j)$ also satisfies
$u\tridoto\omega_0v\omega_0\vartriangleleft_0 w$.
The criteria of Proposition~\ref{prop:one} show that 
either $0<u(i)$ and $0< u(j)$ so that 
$u\tridoto u(\ol{\jmath},\,j)\vartriangleleft_0 w$, 
or else $u(j)\cdot u(i) <0$ and so
$v\tridoto v(\ol{\imath},\,j)\vartriangleleft_0w$.
In the first case, $u(\ol{\jmath},\,j)=u t_j \in{\mathcal B}_n$,
and in the second, 
$v(\ol{\imath},\,j)=u t_{\ol{\imath},\,j}\in{\mathcal B}_n$,
which completes the proof.
\QED

\begin{rem}\label{rem:covers}
If $u(\ol{\jmath},\,i)\tridoto u$, then either 
$u t_{ij}\lessdot_0 u$ or else both $u t_i\lessdot_0 u$ and 
$u t_j\lessdot_0 u$.
\end{rem}

\begin{ex}
We illustrate Theorem~\ref{thm:induced_order} in
Figure~\ref{fig:obtain}.
There, the elements of ${\mathcal B}_3$ are boxed.
\begin{figure}[htb]
$$\epsfxsize=5.4in \epsfbox{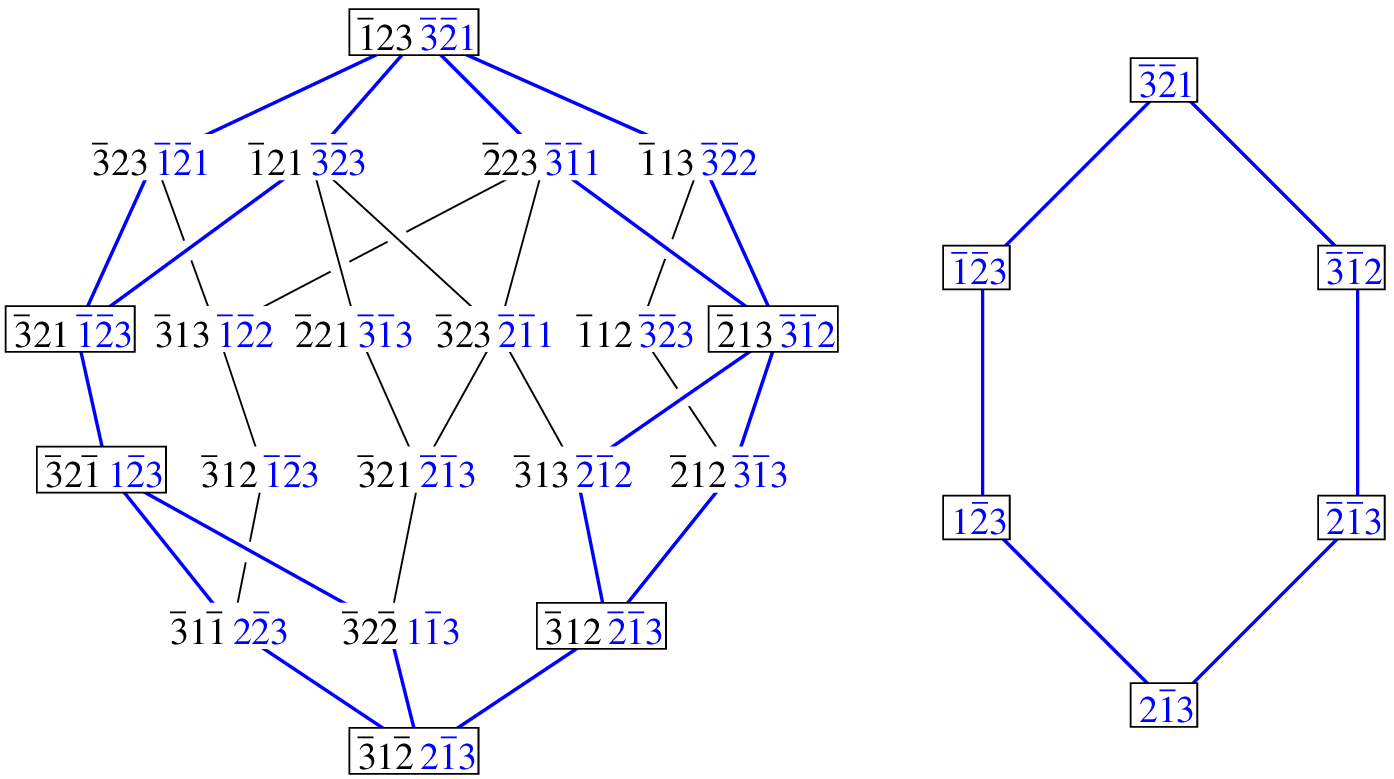}$$
\caption{The intervals $[\ol{3}\,1\,\ol{2}\:2\,\ol{1}\,3,\;
   \ol{1}\,2\,3\:\ol{3}\,\ol{2}\,1]_{\vartriangleleft_0}$
   and $[2\,\ol{1}\,3,\;\ol{3}\,\ol{2}\,1]_0$.
   \label{fig:obtain}}
\end{figure}
\end{ex}

This relation between the two partial orders $({\mathcal B}_\infty,<_0)$
and $({\mathcal S}_{\pm\infty},\vartriangleleft_0)$ makes many 
properties of  $({\mathcal B}_\infty,<_0)$ easy corollaries of the analogous
results for $({\mathcal S}_{\pm\infty},\vartriangleleft_0)$
(established in~\cite{BS98a,BS_monoid}).
We discuss these properties in the remainder of this section, omitting 
proofs. 
\smallskip

The 0-Bruhat order has a non-recursive characterization:

\begin{prop}\label{prop:thmA}
 Let $u,w\in {\mathcal B}_\infty$.
 Then $u\leq_0 w$ if and only if
 \begin{enumerate}
  \item[(1)] $0<i\Longrightarrow u(i)\geq w(i)$, \ and 
  \item[(2)] $0<i<j$ and $u(i)<u(j)\Longrightarrow w(i)<w(j)$.
 \end{enumerate}
\end{prop}

For $P\subset \{1,2,\ldots\}={\mathbb N}$, let 
$\#P\in{\mathbb N}\cup \{\infty\}$ be the cardinality of $P$.
Then the inclusion
$P\hookrightarrow {\mathbb N}$ induces compatible inclusions
$$
\begin{picture}(122,57)(-5,0)
\put(30,2){\vector(1,0){63}}	\put(30,4){\oval(4,4)[l]}
\put(30,47){\vector(1,0){65}}	\put(30,49){\oval(4,4)[l]}
\put(5,37){\vector(0,-1){24}}	\put(7,37){\oval(4,4)[t]}
\put(105,37){\vector(0,-1){24}}	\put(107,37){\oval(4,4)[t]}
\put(60,8){$\varepsilon_P$}	\put(60,52){$\varepsilon_P$}
\put(0,45){${\mathcal B}_{\#P}$}
\put(-5,1){${\mathcal S}_{\pm\#P}$}
\put(100,45){${\mathcal B}_\infty$}
\put(98,1){${\mathcal S}_{\pm\infty}$}
\end{picture}
$$
{\em Shape equivalence} is the equivalence relation on ${\mathcal B}_\infty$
induced by $u \sim \varepsilon_P(u)$ for $P\subset {\mathbb N}$ and 
$u\in{\mathcal B}_{\#P}$.
Let $[u,w]_0:=\{v\mid u\leq_0v\leq_0w\}$ denote the interval in the 0-Bruhat
order between $u$ and $w$, a finite graded poset.
A corollary of Theorems~\ref{thm:induced_order} and Theorem~E(i)
of~\cite{BS98a} is the following fundamental result about the 0-Bruhat order
on ${\mathcal B}_\infty$.
\medskip

\noindent{\bf Theorem~B(1)}
{\em 
Suppose $u,w,x,z\in {\mathcal B}_\infty$ with $wu^{-1}$
shape equivalent to $zx^{-1}$.
Then $[u,w]_0\simeq [x,z]_0$.
If $\varepsilon_P(wu^{-1})=zx^{-1}$, then this isomorphism is given by
$$
[u,w]_0\ni v\ \longmapsto\ \varepsilon_P(vu^{-1})x \in [x,z]_0.
$$
}\medskip

This property is not shared by the $k$-Bruhat order on 
${\mathcal B}_\infty$, for any $k>0$.

\begin{ex}\label{ex:not_nice}
Consider the following two intervals in the 1-Bruhat order on 
${\mathcal B}_4$:
$$\epsfxsize=2.2in \epsfbox{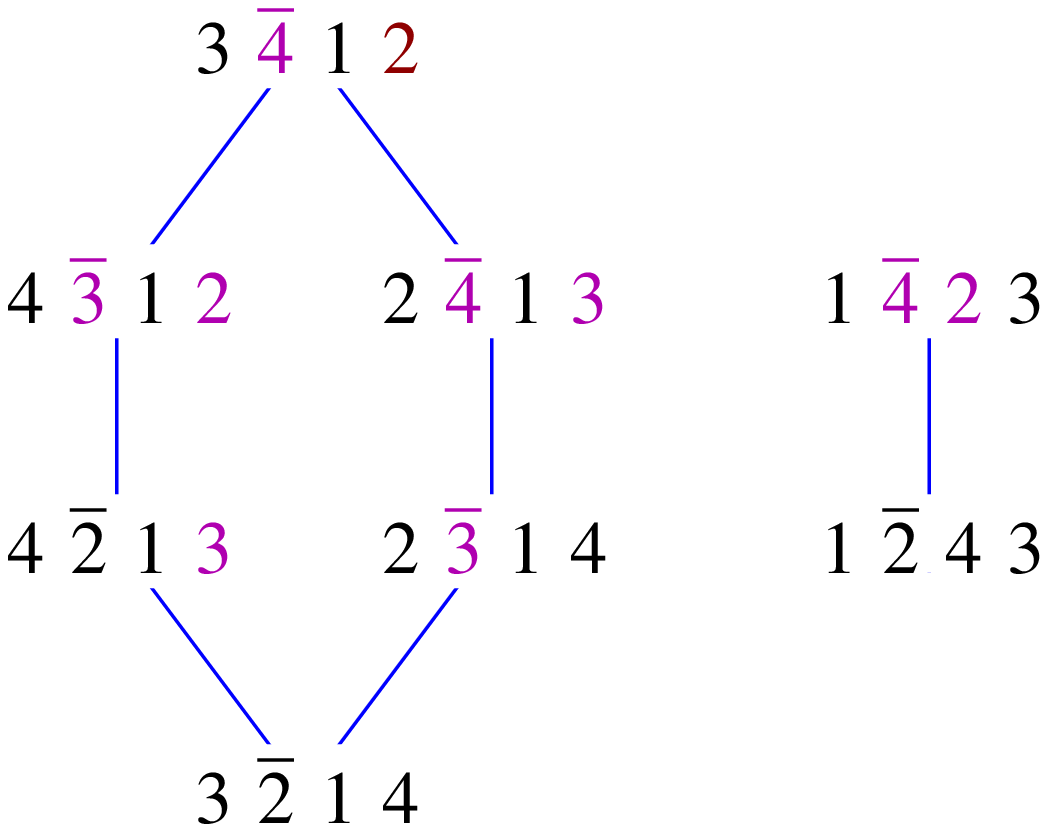}$$
Note that 
$(3\ol{2}14)^{-1}\cdot(3\ol{4}12) = (1\ol{2}43)^{-1}\cdot(1\ol{4}23)
= (\ol{2},\,\ol{4})(2\,,4)$ and the two intervals are not isomorphic.
\end{ex}

\begin{rem}\label{rem:X}
For any $\zeta\in{\mathcal B}_\infty$, there is a 
$u\in{\mathcal B}_\infty$ with $u\leq_0\zeta u$:
Suppose $\{a\in \pm{\mathbb N}\mid a>\zeta(a)\}=\{a_1,\ldots,a_m\}$
with $\zeta(a_1)<\zeta(a_2)<\cdots<\zeta(a_m)$.
Further let 
$\{a_{m+1}<a_{m+2}<\cdots\}= {\mathbb N}\setminus
\{|a_1|,|a_2|,\ldots,|a_m|\}$.
If we define $u(j)=a_j$, then Proposition~\ref{prop:thmA} implies that
$u\leq_0\zeta u$.
Note that if $[m]=\{a>0\mid \zeta(a)\neq a\}$, then $\zeta u$ is a
Grassmannian permutation. 
\end{rem}

By Theorem~B(1), we may define a new
partial order on ${\mathcal B}_\infty$, which we call the {\em Lagrangian
order:}
For $\eta,\zeta\in{\mathcal B}_\infty$, set $\eta\preceq \zeta$ if there is
a $u\in{\mathcal B}_\infty$ with $u\leq_0\eta u\leq_0\zeta u$.
By Remark~\ref{rem:X}, it has a unique minimal element $e$. 
This order is graded by the rank, $\Len(\zeta)$, 
where $\Len(\zeta):=\ell(\zeta u)-\ell(u)$ whenever $u\leq_0\zeta u$.
These notions have definitions independent of $\leq_0$:

\begin{defth}[{\it cf.}\/ Definition 3.2.2~\cite{BS98a}]\label{defth:2.6}
Let $\eta,\zeta\in{\mathcal B}_\infty$.
\begin{enumerate}
\item[(1)]
  Then $\eta\preceq\zeta$ if and only if
   \begin{enumerate}
   \item[(i)] $a\in \pm {\mathbb N}$ with $a>\eta(a)
    \Longrightarrow \eta(a)\geq \zeta(a)$, and 
   \item[(ii)]
    $a,b\in \pm {\mathbb N}$ with $a<b$, 
    $a>\zeta(a),b>\zeta(b)$, and $\zeta(a)<\zeta(b)
    \Longrightarrow \eta(a)<\eta(b)$.\vspace{-10pt}
   \end{enumerate}
   \vspace{15pt}
\item[(2)] 
    ${\displaystyle \Len(\zeta)\ =\ \sum_{a,0>\zeta(a)} |\zeta(a)|
      \;\ -\;\ \#\{(a,b)\mid 0<a<b,  a=\zeta(a), a>\zeta(b)\}}$
	\vspace{5pt}

   ${\displaystyle \qquad\ 
      -\;\ \#\{(a,b)\mid a<b, a>\zeta(a),b>\zeta(b),\zeta(a)>\zeta(b)\}
	\;\ -\  \sum_{0>a>\zeta(a)} |a|}$.
\end{enumerate}
\end{defth}

\noindent{\bf Proof. }
Let $u$ be the permutation with $u\leq_0 \zeta u$ constructed from $\zeta$
in Remark~\ref{rem:X}.
If $u\leq_0 \eta u\leq_0\zeta u$, then $\eta$ satisfies the
conditions in ({\em i}), and conversely.

For ({\em ii}), consider the difference $\ell(\zeta u)-\ell(u)$.
The length of $\zeta u$ is the first sum, plus the number of
inversions of the form $0<i\leq n<j$ with $\zeta u(i)>\zeta u(j)=u(j)$.
In the construction of $u$, each of these is also an inversion in $u$
involving positions $0<i\leq n<j$, and so are canceled in the difference.
The second term counts the remaining inversions of this type in $u$,
the third term counts the inversions with $0<i<j\leq n$ in $u$,
and the fourth term is $\sum_{i>0>u(i)} |u(i)|$.
\QED

The Lagrangian order is the ${\mathcal B}_\infty$-counterpart of the
Grassmann-Bruhat order $\prcc$ on ${\mathcal S}_\infty$
of~\cite{BS98a,BS_monoid}. 
This is defined as follows: 
Let $\eta,\zeta\in{\mathcal S}_\infty$.
Then $\eta\prcc\zeta$ if and only if 
there is a $u\in{\mathcal S}_\infty$ with 
$u\vartriangleleft_0 \eta u\vartriangleleft_0 \zeta u$.
The Grassmann-Bruhat order is ranked with 
$||\eta||:=l(\eta u)-l(u)$ whenever $u\vartriangleleft_0\eta u$.
Let $s(\zeta)$ count the sign changes $\{a>0\mid 0>\zeta(a)\}$ in $\zeta$.
We have the following relation between these two orders. 

\begin{cor}\label{cor:compare_orders}
\mbox{ }
\begin{enumerate}
\item[(1)] $({\mathcal B}_\infty,\prec)$ is an induced suborder of 
      $({\mathcal S}_{\pm\infty},\prcc)$.
\item[(2)] For $\zeta\in{\mathcal B}_\infty(\subset{\mathcal S}_{\pm\infty})$,
      we have $\Len(\zeta)=(||\zeta|| + s(\zeta))/2$.
\end{enumerate}
\end{cor}

\noindent{\bf Proof. }
The first statement is a consequence of
Theorem~\ref{thm:induced_order}.  
For the second statement, consider any maximal chain in $[e,\zeta]_\prec$
(in ${\mathcal B}_\infty$).
By Theorem~\ref{thm:induced_order}, this gives a maximal chain in
$[e,\zeta]_\prcs$ (in ${\mathcal S}_{\pm\infty}$),
where covers of the form $\eta\precdot t_{ab}\,\eta$ are replaced by 
$\eta\pcdot (a,b)\eta\pcdot (a,b)(\ol{a},\ol{b})\eta$.
Thus $||\zeta||=\Len(\zeta)+\tau$, where $\tau$ counts the covers in that
chain if the form $\eta\precdot t_{ab}\,\eta$.
Since only covers of the form $\eta\precdot t_b\eta$ contribute to
$s(\zeta)$, we have   $||\zeta||=2\Len(\zeta)-s(\zeta)$.
\QED

Let $\eta,\zeta\in{\mathcal B}_\infty$.
If $\zeta\cdot\eta=\eta\cdot\zeta$ with 
$\Len(\eta\cdot\zeta)=\Len(\eta)+\Len(\zeta)$,
and neither of $\zeta$ or $\eta$ is the identity, then 
$\eta\cdot\zeta$ is the {\em disjoint product} of $\eta$ and $\zeta$.
(In general $\Len(\eta\cdot\zeta)\leq\Len(\eta)+\Len(\zeta)$.) 
If a permutation cannot be factored in this way, it is {\em irreducible}.
Permutations $\zeta\in{\mathcal B}_\infty$ factor uniquely into
irreducibles.
This is described in terms of
non-crossing partitions~\cite{Kreweras}:
(A non-crossing partition of $\pm{\mathbb N}$ is a set partition such that
if $a<c<b<d$ with $a,b$ in a part $\pi$ and $c,d$ in a part $\pi'$,
then $\pi=\pi'$, as otherwise $\pi,\pi'$ are crossing.)

First, consider $\zeta$ as an element of ${\mathcal S}_{\pm\infty}$.
Let $\Pi$ be the finest non-crossing partition of $\pm{\mathbb N}$
which is refined by the partition given by the cycles of $\zeta$.
For each non-singleton part $\pi$ of $\Pi$, let $\zeta_\pi$ be the product of
the cycles of $\zeta$ which partition $\pi$.
(These $\zeta_\pi$ are the irreducible factors of $\zeta$, as an element of 
${\mathcal S}_{\pm\infty}$.)
Since $\zeta\in{\mathcal B}_\infty$, for each such part $\pi$ of $\Pi$, either
$\pi = \ol{\pi}$ or else one of $\pi,\ol{\pi}$ consists solely
of positive integers.
In the first case, $\zeta_\pi$ is an irreducible factor of $\zeta$ 
(as an element of ${\mathcal B}_\infty$), and in the second,
$\zeta_\pi\zeta_{\ol{\pi}}$ is an irreducible factor of $\zeta$.

The main result concerning this disjointness is the following:

\begin{thm}[{\it cf.}\/ Theorem G (i) \cite{BS98a}]\label{thm:disj-prod}
Suppose $\zeta=\zeta_1\cdots\zeta_s$ is the factorization of 
$\zeta\in{\mathcal B}_\infty$ into irreducibles.
Then the map
$(\eta_1,\ldots,\eta_s)\mapsto \eta_1\cdots\eta_s$
induces an isomorphism
$$
[e,\zeta_1]_\prec\times\cdots\times[e,\zeta_s]_\prec\ 
\stackrel{\sim}{\relbar\joinrel\longrightarrow}\ [e,\zeta]_\prec.
$$
\end{thm}

\begin{figure}[htb]
$$\epsfxsize=1.6in \epsfbox{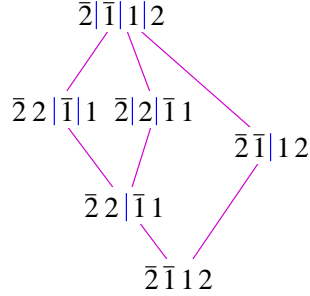}$$
\caption{Non-crossing partitions in $\pm[2]$\label{fig:three}}
\end{figure}

This factorization into irreducibles suggests defining a type $B$
non-crossing partition to be a non-crossing partition $\Pi$ whose blocks
$\pi$ are either stable under $\omega_0$ ($\pi=\ol{\pi}$), or else 
$\pi, \ol{\pi}$ are distinct, with one consisting solely of positive
integers.
These differ from the non-crossing partitions of type $B$
introduced by Reiner~\cite{Reiner_97}, which form a graded lattice.
Figure~\ref{fig:three} shows the partitions of $\pm[2]$ defined here.

We summarize some properties of $({\mathcal B}_\infty,\prec)$.

\begin{thm}[{\it cf.}~Theorem 3.2.3 of~\cite{BS98a}]\label{thm:prpoerties}
\mbox{ } 
\begin{enumerate}
\item[(1)] $({\mathcal B}_\infty,\prec)$ is a graded poset with minimal element
      $e$ and rank function $\Len(\cdot)$.
\item[(2)] The map $\lambda\mapsto v(\lambda)$ exhibits the lattice of strict
      partitions as an induced suborder of $({\mathcal B}_\infty,\prec)$.
\item[(3)] If $u\leq_0\zeta u$, then $\eta\mapsto \eta u$ induces an
      isomorphism  
      $[e,\zeta]_\prec\stackrel{\sim}{\longrightarrow}[u,\zeta u]_0$.
\item[(4)] If $\eta\preceq \zeta$, then $\xi\mapsto \xi\eta^{-1}$ induces an
      isomorphism $[\eta,\zeta]_\prec
                   \stackrel{\sim}{\longrightarrow}[e,\zeta\eta^{-1}]_\prec$.
\item[(5)] For every infinite set $P\subset {\mathbb N}$, the map
      $\varepsilon_P:{\mathcal B}_\infty\rightarrow{\mathcal B}_\infty$ is
      an injection of graded posets.
      Thus if $\eta,\zeta\in{\mathcal B}_\infty$ are shape equivalent, then 
      $[e,\zeta]_\prec\simeq [e,\eta]_\prec$.
\item [(6)]The map $\eta\mapsto \eta\zeta^{-1}$ induces an order-reversing
      isomorphism between $[e,\zeta]_\prec$ and $[e,\zeta^{-1}]_\prec$.
\end{enumerate}
\end{thm}

\begin{figure}[htb]
 $$
  \epsfxsize=6in \epsfbox{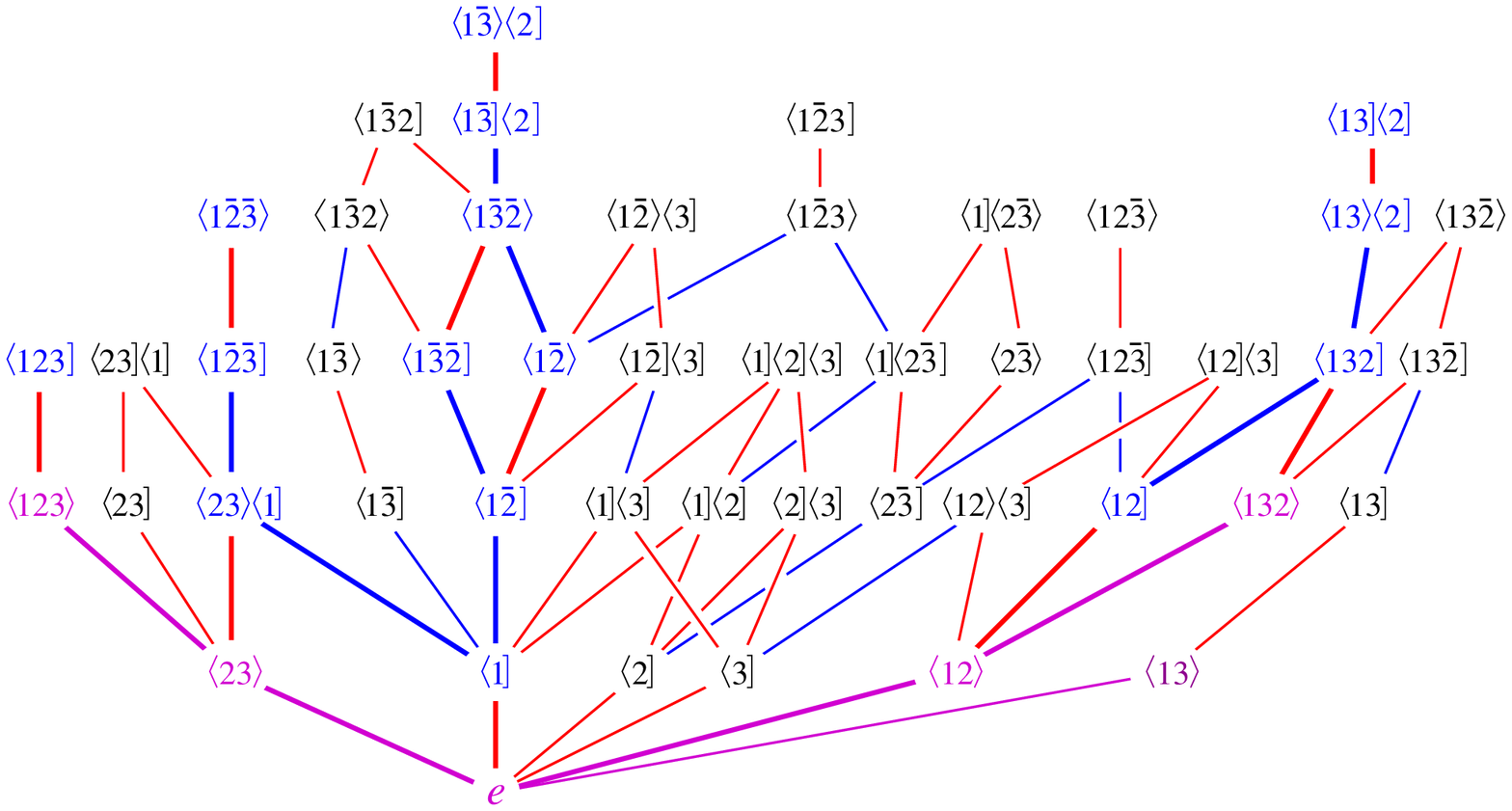}
 $$
 \caption{The Lagrangian order on ${\mathcal B}_3$ \label{fig:new_order}}
\end{figure}

If $\zeta\in{\mathcal B}_\infty$ is factored into into disjoint cycles in 
${\mathcal S}_{\pm\infty}$, the resulting cycles have one of two forms:
$$
(a,\,b,\,\ldots,\,c)\qquad\mbox{or}\qquad
(a,\,b,\,\ldots,\,c,\,\ol{a},\,\ol{b},\,\ldots,\,\ol{c})
$$
with $|a|,|b|,\ldots,|c|$ distinct.
Furthermore, every cycle $\eta=(a,\,b,\,\ldots,\,c)$ of the first type is
paired with another, 
$\ol{\eta}:=(\ol{a},\,\ol{b},\,\ldots,\,\ol{c})$, also of the first type.
This motivates a `cycle notation' for permutations 
$\zeta\in{\mathcal B}_\infty$.
Write $\langle a,\,b,\,\ldots,\,c\rangle$ for the product
$(a,\,b,\,\ldots,\,c)\cdot(\ol{a},\,\ol{b},\,\ldots,\,\ol{c})$
and $\langle a,\,b,\,\ldots,\,c]$ for cycles 
$(a,\,b,\,\ldots,\,c,\,\ol{a},\,\ol{b},\,\ldots,\,\ol{c})$ of the second
type.
Call either of these {\it cycles} in ${\mathcal B}_\infty$.
In some examples and figures, the commas may be omitted.
With this notation, Figure~\ref{fig:new_order}
shows the Lagrangian order on ${\mathcal B}_3$.
The thickened lines are between skew Grassmannian permutations
$v(\lambda)v(\mu)^{-1}$ for $\mu\subset\lambda$. 

\subsection{A monoid for the Lagrangian order}

The `Schubert vs.~Schur' structure constants $b^w_{u\,\lambda}$
and  $c^w_{u\,\lambda}$
are related to the enumeration of (saturated) chains
in the 0-Bruhat order on ${\mathcal B}_\infty$, and hence to the enumeration
of chains in the Lagrangian order on ${\mathcal B}_\infty$.
We develop the elementary theory of chains in these orders
along the lines of~\cite{BS_monoid}.

A chain in either $[u,\zeta u]_0$ or $[e,\zeta]_\prec$ is a particular
factorization of $\zeta$ into transpositions 
$t_b$ and $t_{a\,b}$.
We give an algorithm for finding a chain in $[e,\zeta]_\prec$.
For this, set $t_{\ol{b}\,b}=t_{b\,\ol{b}}=t_b$.

\begin{alg}[{\it cf.}\/ Remark 3.1.2 \cite{BS98a}]\label{alg:chain}
\mbox{ }

\noindent{\tt input: }A permutation $\zeta\in {\mathcal B}_\infty$.

\noindent{\tt output: }Permutations
$\zeta,\zeta_1,\ldots,\zeta_m=e$ such that  
$$
e\ \precdot\ \zeta_{m-1}\ \precdot\ \cdots\ 
\precdot\ \zeta_1\ \precdot\  \zeta
$$
is a saturated chain in the Lagrangian order.

Output $\zeta$.
While $\zeta\neq e$, do
\begin{enumerate}
\item[1] Choose $b\in {\mathbb N}$ maximal subject to 
    $b>\zeta(b)$.
\item[2] Choose $a$ minimal subject to 
    $a\leq \zeta(b)<\zeta(a)$. 
\item[3] $\zeta:=\zeta t_{a\,b}$,  output $\zeta$.
\end{enumerate}

Before every execution of \/ {\rm 3}, 
$\zeta t_{a\,b}\precdot \zeta$.
Moreover, this algorithm terminates in $\Len(\zeta)$ iterations and the
reverse of the sequence produced is a chain in $[e,\zeta]_\prec$.
\end{alg}

The Lagrangian order has a theory of reduced decompositions,
analogous to the usual theory for the Coxeter group $B_\infty$ with respect
to the weak order. 
(Also analogous to that for the Grassmann-Bruhat order on  
${\mathcal S}_\infty$~\cite{BS_monoid}.) 
We express this in the context of monoids.

Define a monoid ${\mathcal M}$ with 0 and generators
${\bf t}_{a\, b}, {\bf t}_b$ for integers $0<a<b$, one for each
reflection of the Weyl group ${\mathcal B}_\infty$. 
To simplify the list of relations these satisfy, set
${\bf t}_{b\,b}:= {\bf t}_b$.
Also, if ${\bf w}\equiv{\bf u}$ is a relation between words ${\bf w}$ and 
${\bf u}$ in these generators, then 
${\bf w}^{\rm op}\equiv{\bf u}^{\rm op}$, where ${\bf w}^{\rm op}$ is the
word ${\bf w}$ read backwards.
That said, these generators are subject to the following relations:

\begin{equation}\label{eq:relations}
\begin{array}{llrclll}
({\rm i})&&{\bf t}_{a\,c}{\bf t}_a{\bf t}_{a\,b}&\equiv&
     {\bf t}_{a\,b}{\bf t}_{b\,c}{\bf t}_b&&
     \mbox{if }a<b<c,\\
({\rm ii})&&{\bf t}_{b\,c}{\bf t}_{c\,d}{\bf t}_{a\,c}&\equiv&
     {\bf t}_{b\,d}{\bf t}_{a\,b}{\bf t}_{b\,c}&&
     \mbox{if }a<b<c<d\\
({\rm iii})&&{\bf t}_{a\,b}{\bf t}_{c\,d}&\equiv&
     {\bf t}_{c\,d}{\bf t}_{a\,b}&&
     \mbox{if $a\leq b<c\leq d$ or $a<c<d<b$},\\
({\rm iv})&&{\bf t}_{a\,c}{\bf t}_b&\equiv& 
     {\bf t}_b{\bf t}_{a\,c}\ \equiv\ 0&&
     \mbox{if }a<b\leq c\mbox{ or }a=b=c,\\
({\rm v})&&{\bf t}_{a\,c}{\bf t}_{b\,d}&\equiv& 
     {\bf t}_{b\,d}{\bf t}_{a\,c}\ \equiv\ 0&&
     \mbox{if }a\leq b<c\leq d,\\
({\rm vi})&&{\bf t}_{b\,c}{\bf t}_{a\,b}{\bf t}_{b\,c}&\equiv&
 {\bf t}_{a\,b}{\bf t}_{b\,c}{\bf t}_{a\,b}\ \equiv\ 0,&&
        \mbox{if }a<b<c.
\end{array} 
\end{equation}

These hold because ${\mathcal M}$ is a sub monoid of the monoid
for the Grassmann-Bruhat order~\cite{BS_monoid} in the same way the
Lagrangian order is an induced 
suborder of the Grassmann-Bruhat order on ${\mathcal S}_\infty$.

The relation between ${\mathcal M}$ and the Lagrangian order on 
${\mathcal B}_\infty$ is obtained via a faithful representation of
${\mathcal M}$ as linear operators on ${\mathbb Q}{\mathcal B}_\infty$:
Define linear operators $\wh{\bf t}_{a\,b}$ 
and $\wh{\bf t}_b$ on ${\mathbb Q}{\mathcal B}_\infty$ by 
\begin{eqnarray*}
  \wh{\bf t}_{a\,b}.\zeta &:=&
  \left\{\begin{array}{ll}
  t_{a\,b}\zeta&\mbox{ if } \Len(\zeta)+1=|t_{a\,b}\zeta|\\
  0&\mbox{ otherwise.}\end{array}\right.\\
  \wh{\bf t}_b.\zeta &:=&
  \left\{\begin{array}{ll}
  t_a\zeta&\mbox{ if } \Len(\zeta)+1=|t_b\zeta|\\
  0&\mbox{ otherwise.}\end{array}\right.
\end{eqnarray*}

To simplify the following statement, let the index $\alpha$ represent
either of $b$ or $a\,b$.

\begin{thm}[{\it cf.} Theorem~1.1~\cite{BS_monoid}] \mbox{ }
\begin{enumerate}
\item[(1)] 
    The operators $\wh{\bf t}_\alpha$ satisfy the relations (1.1), and a
    composition of operators is characterized by its value at the identity. 
    That is,
    $\wh{\bf t}_{\alpha'_m}\cdots\wh{\bf t}_{\alpha'_1}
    =\wh{\bf t}_{\alpha_n}\cdots\wh{\bf t}_{\alpha_1}$ 
    if and only if\/
    $\wh{\bf t}_{\alpha'_m}\cdots\wh{\bf t}_{\alpha'_1}e
    =\wh{\bf t}_{\alpha_n}\cdots\wh{\bf t}_{\alpha_1}e$.
\item[(2)] 
    For 
    ${\bf x}={\bf t}_{\alpha_n}\cdots{\bf t}_{\alpha_2}{\bf t}_{\alpha_1}
    \in {\mathcal M}$, the map 
    ${\bf x}\mapsto \wh{x}:=
    \wh{\bf t}_{\alpha_n}\cdots\wh{\bf t}_{\alpha_2}\wh{\bf t}_{\alpha_1}$
    is a faithful representation of ${\mathcal M}$.
\item[(3)] 
    The following map is a well defined bijection:
    \begin{eqnarray*}
    {\mathcal M}&\longrightarrow& {\mathcal B}_\infty\cup\{0\},\\
    {\bf x}&\longmapsto& \wh{\bf x} e.
    \end{eqnarray*}
\item[(4)] 
    For $\zeta,\eta\in {\mathcal B}_\infty$, 
    $\eta\preceq \zeta$ if and only if there exists 
    ${\bf x}\in {\mathcal M}$ such that $\zeta=\wh{\bf x}\eta$.  
\item[(5)] 
    The set 
    ${\bf R}(\zeta)=\{\wh{\bf x}:\wh{\bf x}e=\zeta\}$
    corresponds to the set of all maximal chains in $[e,\zeta]_\prec$.
\end{enumerate}
\end{thm}

We call the elements of ${\bf R}(\zeta)$ the {\it $\prec$-reduced
decompositions} of $\zeta$.

\section{Isotropic flag manifolds and maximal Grassmannians}

Let $V$ denote either ${\mathbb C}^{2n+1}$
equipped with non-degenerate symmetric bilinear form  or
${\mathbb C}^{2n}$ equipped with a
non-degenerate alternating bilinear form.
In the first case, $V$ an {\it odd orthogonal (vector) space}, and in the
second, a {\it symplectic (vector) space}.
A linear subspace $K$ of $V$  is {\em isotropic} if the restriction of 
the form to $K$ is identically zero.
Isotropic subspaces have dimension at most $n$.
An {\em isotropic flag} in $V$ is a sequence  $\Edot$ of isotropic
subspaces:
$$
\Edot\ :\ E_{\ol{n}}\ \subset\ E_{\ol{n-1}}\ \subset\ \cdots\ 
\subset\ E_{\ol{1}},
$$
where $\dim E_{\ol{\imath}} = n+1-i$.
Let $K^\perp$ be the annihilator of a subset $K$ of $V$.
Given an isotropic flag $\Edot$ in $V$, 
we obtain a canonical complete flag  in $V$  (also written $\Edot$) 
by defining $E_i:=E_{\ol{i+1}}^\perp$ for $i=1,\ldots,n$, and in the odd
orthogonal case,  $E_0:=E_{\ol{1}}^\perp$.
Henceforth, flags will always be complete, although we may only specify 
$E_{\ol{n}},\ldots,E_{\ol{1}}$.

The group $G$ of linear transformations of $V$ which preserve the given form
acts transitively on the set of isotropic
flags in $V$.
Since the stabilizer of an isotropic flag is a Borel subgroup $B$ of $G$,
this exhibits the set of isotropic flags as the homogeneous space $G/B$.
Here, $G$ is either $\Son$ (odd orthogonal) or $\Spn$ (symplectic).
Similarly, $G$ acts transitively on the set of maximal isotropic subspaces 
of $V$, exhibiting it as the homogeneous space $G/P_0$.
Here $P_0$ is the stabilizer of a maximal isotropic subspace, a maximal
parabolic subgroup associated to the simple root of exceptional length.
Let $\pi: G/B \twoheadrightarrow G/P_0$ be the projection map.

The rational cohomology rings~\cite{Borel_cohomology} of $G/B$ 
for both the symplectic and odd-orthogonal flag manifolds are isomorphic to 
$$
{\mathbb Q}[x_1,\ldots,x_n]/ 
\langle e_i(x_1^2,\ldots,x_n^2) , i=1,\ldots,n\rangle,
$$
where $e_i(a_1,\ldots,a_n)$ is the $i$th elementary symmetric polynomial in
$a_1,\ldots,a_n$. 
However, their integral cohomology rings differ~\cite{EdGr}:
\begin{eqnarray*}
  H^*Sp_{2n}{\mathbb C}/B&\simeq& {\mathbb Z}[x_1,\ldots,x_n]/ 
      \left\langle e_i(x_1^2,\ldots,x_n^2) \right\rangle.\\
  H^*So_{2n+1}{\mathbb C}/B&\simeq& 
   {\mathbb Z}[x_1,\ldots,x_n, c_1,\ldots,c_n]/ I,\\
  I &=& \left\langle e_i(x_1^2,\ldots,x_n^2),\, 2c_i-e_i(x_1,\ldots,x_n), \,
c_{2i}=(-1)^ic_i^2\right\rangle.
\end{eqnarray*}
These rings have another description in terms of Schubert classes.

\subsection{Schubert varieties}
Since an isotropic flag $\Edot\in G/B$ is also a complete flag in $V$, 
we have a canonical embedding $G/B \hookrightarrow {\mathbb F}\ell(V)$, the
manifold of complete flags in $V$. 
Similarly, there is an embedding $G/P_0 \hookrightarrow {\bf G}_n(V)$,
the Grassmannian of $n$-dimensional subspaces of $V$.
We use these maps to understand some structures of $G/B$.

Given $w\in{\mathcal B}_n$ and an isotropic (complete) flag $\Edot\in G/B$,
the {\it Schubert variety} $Y_w\Edot$ of $G/B$ is the collection of all
flags $\Fdot\in G/B$ satisfying
\begin{equation}\label{eq:flag_def}
  \dim E_i \cap F_{\ol{\jmath}}\ \ \geq\ \ \#\{n\geq l\geq j\mid w(l)\leq i\}
\end{equation}
for each $j\in[n]$ and $-n\leq i\leq n$ ($i\neq 0$ in the symplectic case).
This has codimension $\ell(w)$ in $G/B$.
Also, $Y_w\Edot\subset Y_u\Edot$ if and only if $u\leq w$ in the 
Bruhat order.
The Schubert cell $Y^\circ_u\Edot$ is the set of flags $\Fdot$
for which equality holds in~(\ref{eq:flag_def}).
These are the flags in $Y_u\Edot$ which are not in any sub-Schubert variety
($Y_w\Edot$ with $u\leq w$). 

If now $\Edot\in{\mathbb F}\ell(V)$ and $w\in{\mathcal S}_{[\ol{n},n]}$
(${\mathcal S}_{\pm[n]}$ in the symplectic case), 
then the Schubert variety $X_w\Edot$ of ${\mathbb F}\ell(V)$ is the
collection of flags $\Fdot\in{\mathbb F}\ell(V)$
satisfying~\ref{eq:flag_def}) for all $\ol{n}\leq i,j\leq n$ 
($i,j\neq 0$ in the symplectic case).
Furthermore, if $w\in{\mathcal B}_n$ and
$\Edot\in G/B$, then 
$$
Y_w\Edot\  =\ G/B\bigcap X_w\Edot
$$
and  this is a scheme-theoretic equality~\cite{LS_geom_78}.

The Schubert cells constitute a cellular decomposition of $G/B$.
Thus {\em Schubert classes}, the cohomology classes Poincar\'e dual to the
fundamental cycles of Schubert varieties, form ${\mathbb Z}$-bases
for these cohomology rings.
Write ${\mathfrak B}_w$ for the class $[Y_w\Edot]$ in 
$H^*\Son/B$ Poincar\'e dual to the fundamental cycle of 
$Y_w\Edot$ of $\Son/B$ and ${\mathfrak C}_w$ for the 
corresponding class in $H^*\Spn/B$.
Since these are bases, there are integral structure constants 
$b^w_{u\,v}$ and $c^w_{u\,v}$ for $u,w,v\in{\mathcal B}_n$ defined by the
identities:
$$
{\mathfrak B}_u \cdot{\mathfrak B}_v \ =\ \sum_w b^w_{u\,v}{\mathfrak B}_w
\qquad\mbox{and}\qquad
{\mathfrak C}_u\cdot {\mathfrak C}_v \ =\ \sum_w c^w_{u\,v}{\mathfrak C}_w.
$$

Let $s(w)$ count the number of sign changes in the permutation $w$.
Then the isomorphism of rational rings 
is induced by the map~\cite{BH}:
$$
{\mathfrak C}_w\ \longmapsto\ 2^{s(w)}{\mathfrak B}_w.
$$
Thus
\begin{equation}\label{eq:coeff_compare}
    2^{s(u)+s(v)}b^w_{u\,v}\ =\ 2^{s(w)}c^w_{u\,v}.
\end{equation}
It suffices to establish identities and formulas for 
$\Spn/B$.
We do this, because a crucial geometric result
(Theorem~\ref{thm:geom_facts}(2)) does not hold for $\Son/B$.
\medskip

Two flags $\Edot, \Epdot$ are {\it opposite} if 
$\dim(E_i\bigcap E'_{\ol{\imath}})=1$ for all $i$.
In what follows, $Y_u$ and $Y'_v$, will always denote Schubert
varieties defined by fixed, but arbitrary, opposite isotropic flags.
A consequence of Kleiman's theorem on the transversality of a general
translate~\cite{Kleiman}, results in~\cite{Deodhar}, and some combinatorics,
is the following proposition.

\begin{prop}\label{prop:transverse}
Let $u,w\in{\mathcal B}_n$.
Then $Y_u\bigcap Y'_{\omega_0w}\neq 0$ if and only if $u\leq w$ in the Bruhat
order.
If $u\leq w$, then $Y_u, Y'_{\omega_0w}$ meet generically transversally, and
the intersection cycle is irreducible of dimension $\ell(w)-\ell(u)$.
\end{prop}

The top-dimensional component of $H^* G/B$ is generated by the class of a
point $[\mbox{\rm pt}]={\mathfrak B}_{\omega_0}$ or 
${\mathfrak C}_{\omega_0}$.
The map $\deg:H^*G/B \rightarrow {\mathbb Z}$ selects the coefficient of
[pt] in  a cohomology class.
The {\it intersection pairing} on $H^*G/B$ is the composition
$$
\beta, \gamma\in H^*G/B \ \longmapsto\ \deg (\beta\cdot\gamma).
$$
By Proposition~\ref{prop:transverse}, the product
$[Y_u]\cdot [Y_v]$ is the cohomology class 
$[Y_u\bigcap Y'_v]$.
In particular, when $v=\omega_0 u$, these intersections are single, reduced
points, so that $[Y_u]$ and $[Y_{\omega_0 u}]$ are dual under the 
intersection pairing.
Thus 
$$
c^w_{u\,v}\ =\ \deg ({\mathfrak C}_u\cdot {\mathfrak C}_{\omega_0 w}
\cdot {\mathfrak C}_v),
$$
which is also the number of points in the intersection
$$
Y_u\bigcap Y'_{\omega_0 w}\bigcap Y''_v,
$$
where $Y''_v$ is defined by a flag ${E_{\,\DOT}''}$ opposite to both 
$\Edot$ and $\Epdot$ (which define $Y_u$ and $Y'_{\omega_0 w}$).

We derive a useful description of flags in 
$Y^\circ_u\bigcap {Y'}^\circ_{\omega_0w}$
when $u\leq_0 w$.
For $S\subset {\mathbb C}^m$, let $\Span{S}$ be the linear
span of $S$. 

\begin{lem}\label{lem:param}
Suppose that $u\leq_0 w$ and $\Edot, \Epdot$ are opposite isotropic flags
in $V$.
Then there are algebraic functions
$g_{\ol{\jmath}}:Y^\circ_u\Edot\bigcap 
Y^\circ_{\omega_0w}\Epdot\rightarrow V$ for $1\leq j\leq n$ 
such that for each flag 
$\Fdot\in Y^\circ_u\Edot\bigcap Y^\circ_{\omega_0w}\Epdot$
and each $1\leq j\leq n$, 
\begin{enumerate}
\item[(1)] $F_{\ol{\jmath}}=
	\left\langle g_{\ol{n}}(\Fdot),\ldots,
           g_{\ol{\jmath}}(\Fdot)\right\rangle$, \ 
   and  
\item[(2)] $g_{\ol{\jmath}}(\Fdot)\in E_{u(j)}\bigcap
      E'_{\ol{w(j)}}$.
\end{enumerate}
\end{lem}

\noindent{\bf Proof. }
The representation of Schubert cells via parameterized matrices give
$V$-valued functions $f_{\ol{\jmath}}$ defined on the Schubert
cell  
$Y^\circ_u\Edot$ such that if $\Fdot$ is a flag in that cell, then
$F_{\ol{\jmath}}=\Span{f_{\ol{n}}(\Fdot),\ldots,f_{\ol{\jmath}}(\Fdot)}$,
and $f_{\ol{\jmath}}\in E_{u(j)}$.

Construct the functions $g_{\ol{\jmath}}$ inductively.
First, set $g_{\ol{n}}(\Fdot):= f_{\ol{n}}(\Fdot)$ for 
$\Fdot\in Y^\circ_u\Edot\bigcap Y^\circ_{\omega_0w}\Epdot$.
Since $F_{\ol{n}}\subset E_{u(n)}\bigcap E'_{\ol{w(n)}}$, conditions 1 and 2
are satisfied for $g_{\ol{n}}$.
Suppose we have constructed $g_{\ol{\imath}}$ for $n\geq i>j$.
Let $g_{\ol{\jmath}}(\Fdot)$ be the intersection of $E'_{\ol{w(j)}}$ 
with the affine space
$$
W_j\ :=\ f_{\ol{\jmath}}(\Fdot)\ +\ 
\Span{g_{\ol{\imath}}(\Fdot) \mid i>j\mbox{ and } w(i)<w(j)}.
$$
There is a unique point of intersection:
Since $\Fdot\in Y^\circ_{\omega_0w}\Epdot$, 
$$
\dim E'_{\ol{w(j)}}\bigcap F_{\ol{\jmath}}\ =\ 
\#\{i\mid i\geq j\mbox{ and } w(i)\geq w(j)\}. 
$$
Since $u\leq_0 w$, if $i>j$ and $w(i)<w(j)$, then necessarily
$u(i)<u(j)$, by condition 2 of Proposition~\ref{prop:thmA}.
Hence $W_j\subset E_{u(j)}$ and so 
$g_{\ol{\jmath}}(\Fdot)\in E_{u(j)}\bigcap E'_{\ol{w(j)}}$.   
\QED

Schubert varieties $\Upsilon_\lambda$ of $G/P_0$ are indexed by strict
partitions $\lambda$: decreasing
sequences $n\geq \lambda_1>\cdots>\lambda_k$ of positive integers.
The projection map $\pi: G/B \twoheadrightarrow G/P_0$ maps Schubert
varieties to Schubert varieties, with $\pi Y_u = \Upsilon_\lambda$, where
$\lambda$ consists of the positive numbers among 
$\{\ol{u(1)},\ol{u(2)},\ldots,\ol{u(n)}\}$ arranged in decreasing order.
For a strict partition $\lambda$, let $v(\lambda)$ be the Grassmannian
permutation whose (initial) negative values are 
$\ol{\lambda_1}<\ol{\lambda_2}<\cdots<\ol{\lambda_k}$.
Let $\lambda^c$ be the decreasing sequence obtained from 
the integers in $[n]$ which do not appear in $\lambda$.
Then  $Y_{v(\lambda)}=\pi^{-1}\Upsilon_\lambda$ and
$$
\pi\ :\ Y_{\omega_0v(\lambda)} \ \longrightarrow\  \Upsilon_{\lambda^c}
$$
is birational.
One may see this by considering typical elements of their Schubert cells.

Set $P_\lambda := [\Upsilon_\lambda]$ in $H^*\Son/P_0$
(equivalently  $P_\lambda := [Y_{v(\lambda)}]$ in $H^*\Son/B$)
and let $Q_\lambda$ be the corresponding class in the symplectic case.
For $1\leq m\leq n$, the {\em special Schubert variety}
$\Upsilon_{(m)}$ is the collection of all
maximal isotropic subspaces which meet a fixed $(n+1-m)$-dimensional
isotropic subspace. 
Let $p_m$ (respectively $q_m$) denote the class $P_{(m)}$ in
either $H^*\Son/B$ or $H^*\Son/P_0$ 
(respectively the class $Q_{(m)}$ in either $H^*\Spn/B$ or $H^*\Spn/P_0$).

We are particularly interested in the constants 
$b^w_{u\,\lambda}:=b^w_{u\,v(\lambda)}$ and 
$c^w_{u\,\lambda}:=c^w_{u\,v(\lambda)}$ which give the 
structure of the cohomology of $G/B$ as a module over the
cohomology of $G/P_0$.
Using the intersection pairing and the  projection
formula ({\it cf.}~\cite[8.1.7]{Fulton_intersection}),
we have 
\begin{eqnarray*}
c^w_{u\,\lambda} &=&
   \deg({\mathfrak C}_u\cdot{\mathfrak C}_{\omega_0w}\cdot \pi^*(Q_\lambda))\\
&=&\deg (\pi_*({\mathfrak C}_u\cdot{\mathfrak C}_{\omega_0w})\cdot Q_\lambda),
\end{eqnarray*}
and a similar formula for $b^w_{u\,\lambda}$.
Our main technique will be to find formulas for
$\pi_*({\mathfrak C}_u\cdot{\mathfrak C}_{\omega_0w})$ by studying the effect
of the map $\pi$ on the cycle 
$Y_u\cap Y'_{\omega_0w}$.

To that end, define 
${\mathcal Y}^w_u:=\pi(Y_u\cap Y'_{\omega_0w})$.
These cycles ${\mathcal Y}^w_u$ are, like Schubert varieties, defined
up to translation by $G$.
In the theorems below, write ${\mathcal Y}^w_u={\mathcal Y}^z_x$ to mean
that the cycles may be carried onto each other by an  element of $G$.
(The proofs are more explicit.)

The main result of Section~\ref{sec:identity} concerns these cycles.

\begin{thm}\label{thm:geometric_shapes}
Let $u,w\in{\mathcal B}_n$ with $u\leq_0w$.
Then
\begin{enumerate}
\item[(1)] The map $\pi : Y_u\cap Y'_{\omega_0w} \rightarrow {\mathcal Y}^w_u$
     has degree 1.
\item[(2)] If we have $x,z\in{\mathcal B}_n$ with $z\leq_0z$ and $wu^{-1}$ shape
     equivalent to $zx^{-1}$, then 
     ${\mathcal Y}^w_u={\mathcal Y}^z_x$.
\end{enumerate}
\end{thm}
This implies the identity of Theorem~B(2).
As a consequence of Theorem~\ref{thm:geometric_shapes}(2), define
${\mathcal Y}_\zeta:={\mathcal Y}^{\zeta u}_u$ for any 
$\zeta,u\in{\mathcal B}_n$ with $u\leq_0\zeta u$.

These cycles satisfy more identities.
Let $W$ be a $2m$-dimensional symplectic space.
We study the manifolds of maximal isotropic subspaces of $V,
W$, and $V\oplus W$, together with a map
$$
\Xi\ :\ \Spn/P_0 \times \mbox{\it Sp}_{2m}{\mathbb C}/P_0
\ \longrightarrow\  {\it Sp}_{2n+2m}{\mathbb C}/P_0,
$$
defined by $(H,K)\mapsto H\oplus K$.
Also define $\rho\in{\mathcal B}_n$ by
$\rho(i)=i{-}1{-}n$ for $1\leq i\leq n$
and $\gamma\in{\mathcal B}_n$ by
$\gamma(i)=i{+}1$ for $1\leq i< n$ and $\gamma(n)=1$.

\begin{thm}\label{thm:geom_facts}
\mbox{ }
\begin{enumerate}
\item[(1)] For any $\zeta\in{\mathcal B}_n$, 
      ${\mathcal Y}_{\zeta}\ =\ {\mathcal Y}_{\zeta^{-1}}$.
\item[(2)] If $\eta\cdot\zeta$ is a disjoint product in ${\mathcal B}_{n+m}$
     with $\eta'\in{\mathcal B}_n$ shape equivalent to $\eta$ and 
     $\zeta'\in{\mathcal B}_m$ shape equivalent to $\zeta$, then
     $$
     \Xi( {\mathcal Y}_{\eta'}\times{\mathcal Y}_{\zeta'})\ =\ 
     {\mathcal Y}_{\eta\cdot\zeta}.
     $$
\item[(3)] For any $\zeta\in{\mathcal B}_n$,
     $$
     {\mathcal Y}_{\zeta}\ =\ {\mathcal Y}_{\rho\zeta\rho}.
     $$
\item[(4)] For any $\zeta\in{\mathcal B}_n$ with 
     $a\cdot\zeta(a)>0$ for every $a$,
     $$
     {\mathcal Y}_{\zeta}\ =\ {\mathcal Y}_{\gamma\zeta\gamma^{-1}}.
     $$
\end{enumerate}
\end{thm}

The first statement follows from the observation that
$Y_u\bigcap Y_{\omega_0\zeta u}= 
Y_{\omega_0\zeta u}\bigcap Y_{\omega_0\zeta^{-1}(\omega_0\zeta u)}$.
Since $\Len(\zeta)=\dim {\mathcal Y}_\zeta$, Theorem~C(1) follows from part
(3), by the projection formula.
Similarly, Theorem~C(2) follows from part (4).
We remark that part (2) is true only for the symplectic case, while
(1), (3), and (4) hold also for the odd orthogonal flag manifold.
Statements (2), (3), and (4) are proven in Section~\ref{sec:more_identities}.

\section{Identities of structure constants}
\label{sec:identity}
We establish Theorem~\ref{thm:geometric_shapes} which implies 
Theorem~B(2).
As in Section~\ref{sec:orders}, many results and methods are similar to
those  of~\cite{BS98a} for analogous results about $\Sln{\mathbb C}/B$. 
Our discussions are therefore brief.
The results here hold for both 
$\Son/B$ and $\Spn/B$, with nearly identical proofs.
We only provide justification for $\Spn/B$.

Let $H_2=\Span{h,\ol{h}}\simeq {\mathbb C}^2$ be a symplectic space of
dimension 2.
Then the orthogonal direct sum $V\oplus H_2$ is a symplectic space of
dimension $2n+2$.
For each $1\leq p\leq n+1$, define embeddings
$\psi_p,\psi_{\ol{p}}\ :\ \Spn/B \ \hookrightarrow\ 
\mbox{\it Sp}_{2n+2}{\mathbb C}/B$,
the space of isotropic flags in $V\oplus H_2$,
by
$$
\left(\psi_p \Edot\right)_j\ = \ 
\left\{\begin{array}{ccl}
E_{j+1}&& j\leq \ol{p}\: (<0)\\
\Span{E_j,h}&& \ol{p} <j<0
\end{array}\right..
$$
Define $\psi_{\ol{p}}$ by replacing $h$ with $\ol{h}$ in the definition
above. 
We find the effect of these maps on cohomology by determining
the image of a Schubert variety under $\psi_p$.

First, define two maps between ${\mathcal B}_n$ and 
${\mathcal B}_{n+1}$.
For every $1\leq p\leq n+1$ and $q\in \pm[n+1]$, define the injection
$\varepsilon_{p,q}:{\mathcal B}_n\hookrightarrow{\mathcal B}_{n+1}$ by:
$$
\varepsilon_{p,q} (w)(j) \quad =\quad \left\{\begin{array}{lcl}
w(j)     && j < p \mbox{ and } |w(j)|  <  |q|\\
w(j)-1   && j < p \mbox{ and } w(j) \leq -|q|\\
w(j)+1   && j < p \mbox{ and } w(j) \geq  |q|\\
q        && j = p\\
w(j-1)   && j > p \mbox{ and } |w(j)|  <  |q|\\
w(j-1)-1 && j > p \mbox{ and } w(j) \leq -|q|\\
w(j-1)+1 && j > p \mbox{ and } w(j) \geq  |q|
\end{array}\right..
$$
Let $/_p:{\mathcal B}_{n+1}\rightarrow{\mathcal B}_n$ be the left inverse of
$\varepsilon_{p,q}$, defined by $\varepsilon_{p,w(p)}(w/_p)=w$.
If we represent permutations as permutation matrices, then the effect of
$/_p$ on $w\in{\mathcal B}_{n+1}$ is to erase the $p$th and $\ol{p}$th
columns and the $w(p)$th and $\ol{w(p)}$th rows.
The effect of $\varepsilon_{p,q}$ on $w$ is to expand its permutation matrix
with new $p$th, $\ol{p}$th columns and $q$th, $\ol{q}$th rows filled with
zeroes, except for 1s at positions $(q,p)$ and $(\ol{q},\ol{p})$.
For example:
$$
\begin{array}{c}
\varepsilon_{3,\ol{2}}(\ol{2}\,3\,4\,1)\ =\ \ol{3}\,4\,\ol{2}\,5\,1\\
\varepsilon_{3,2}(\ol{2}\,3\,4\,1)\ =\ \ol{3}\,4\,2\,5\,1
\end{array} 
\qquad\mbox{and}\qquad
4\,\ol{1}\,5\,\ol{2}\,3/_4\ =\ 3\,\ol{1}\,4\,2.
$$

These definitions imply the following proposition
({\it cf.}~\cite{Sottile96}, Lemma 12).

\begin{prop}\label{prop:subsets}
Let $w\in{\mathcal B}_n$, $1\leq p, |q|\leq n+1$, 
and $\Edot$ any isotropic flag.
Then
$$
\psi_p\, Y_w \Edot \ \subset\ Y_{\varepsilon_{p,q}(w)}\psi_{\ol{q}}\Edot.
$$
\end{prop}

Recall that $e$ is the identity permutation and $Y_e$ is the flag manifold
$G/B$. 

\begin{cor}\label{cor:intersection}
Let $\Edot,\Epdot$ be opposite isotropic flags.
Then for any $q\in \pm[n+1]$,
$\psi_q\Edot, \psi_{\ol{q}}\Epdot$ are opposite flags,
and for any $1\leq p\leq n+1$,
$\psi_p\, Y_w\Edot$ equals either of the following cycles
$$
Y_{\varepsilon_{p,\ol{n+1}}(w)}\psi_{n+1}\Edot \bigcap
Y_{\varepsilon_{p,n+1}(e)}\psi_{\ol{n+1}}\Epdot
\qquad\quad 
Y_{\varepsilon_{p,n+1}(w)}\psi_{\ol{n+1}}\Edot \bigcap
Y_{\varepsilon_{p,\ol{n+1}}(e)}\psi_{n+1}\Epdot.
$$
\end{cor}

\noindent{\bf Proof. }
It is straightforward to check that the flags are opposite.
Moreover, by Proposition~\ref{prop:subsets}, 
$\psi_p\, Y_w\Edot$ is a subset of either intersection, as 
$Y_e=\Spn/B$.
Since
$$
\ell(\varepsilon_{p,\ol{n+1}}(w))\ =\ \ell(w)+n+p,
\qquad\quad
\ell(\varepsilon_{p,n+1}(w))\ =\ \ell(w) +n+1-p,
$$
and $\dim \Spn/B=n^2$, Proposition~\ref{prop:transverse} implies
that all three cycles are irreducible with the same dimension, proving their
equality.
\QED

\begin{cor}\label{cor:push-forward}
For any $w\in{\mathcal B}_n$ and $1\leq p\leq n$, we have
$$
(\psi_p)_*{\mathfrak C}_w\ =\ 
{\mathfrak C}_{\varepsilon_{p,\ol{n+1}}(w)}\cdot
        {\mathfrak C}_{\varepsilon_{p,n+1}(e)}\ =\ 
{\mathfrak C}_{\varepsilon_{p,n+1}(w)}\cdot
        {\mathfrak C}_{\varepsilon_{p,\ol{n+1}}(e)}.
$$
\end{cor}

\begin{rem}
As in Section~4.2 of~\cite{BS98a}, the Schubert classes
${\mathfrak C}_{\varepsilon_{p,n+1}(e)}$ and 
${\mathfrak C}_{\varepsilon_{p,\ol{n+1}}(e)}$ are certain special Schubert
classes 
from Grassmannian projections.
Were the corresponding Pieri-type formulas known, we could deduce
formulas for $(\psi_p)_*{\mathfrak C}_w$.
Only one of these classes is a special Schubert class from $G/P_0$:
$$
{\mathfrak B}_{\varepsilon_{1,\ol{n+1}}(e)}\ =\ p_{n+1}
\qquad\mbox{and}\qquad
{\mathfrak C}_{\varepsilon_{1,\ol{n+1}}(e)}\ =\ q_{n+1}.
$$
We deduce formulas for $(\psi_1)_*{\mathfrak B}_w$ and
$(\psi_1)_*{\mathfrak C}_w$ from this.

Recall that $H^*\Spn/B$ is generated by $x_1,\ldots,x_n$.
These classes are Chern classes of certain line bundles on $\Spn/B$:
Let ${\mathcal E}_{\DOT}\rightarrow\Spn/B$ be the flag of bundles whose
fibre at $\Edot$ is $\Edot$.
Then $x_i=-c_1({\mathcal E}_{\ol{i+1}}/{\mathcal E}_{\ol{\imath}})$.
Thus
$$
\psi_p^*(x_i)\ =\ \left\{\begin{array}{lcl}
                                x_i&\ &i<p\\
                                0  && i=p\\
                                x_{i-1}&&i>p
                         \end{array}\right..
$$
\end{rem}

\begin{thm}
Let $v\in{\mathcal B}_{n+1}$.
Then
\begin{enumerate}
\item[(1)] In $H^*\Son/B$,  $\psi_1^*{\mathfrak B}_v=
   \sum \theta(\varepsilon_{1,\ol{n+1}}(y)v^{-1})\ {\mathfrak B}_y$, the
   sum over all $y\in{\mathcal B}_n$ with 
   $\ell(v)=\ell(y)$.
  
\item[(2)] In $H^*\Spn/B$, $\psi_1^*{\mathfrak C}_v=
   \sum \chi(\varepsilon_{1,\ol{n+1}}(y)v^{-1})\ {\mathfrak C}_y$, the
   sum over all $y\in{\mathcal B}_n$ with 
   $\ell(v)=\ell(y)$.
\end{enumerate}
\end{thm}

\noindent{\bf Proof. }
These are consequences of the projection formula and
Corollary~\ref{cor:push-forward}.
For the first, 
$$ 
\psi_1^*{\mathfrak B}_v\ =\ \sum_{y\in{\mathcal B}_n}
          \deg(\psi_1^*{\mathfrak B}_v\cdot{\mathfrak B}_{\omega_0 y})
	   {\mathfrak B}_y.
$$
But
$$
\deg(\psi_1^*{\mathfrak B}_v\cdot{\mathfrak B}_{\omega_0 y})\ =\ 
\deg({\mathfrak B}_v\cdot(\psi_1)_*({\mathfrak B}_{\omega_0 y}))\ =\ 
\deg({\mathfrak B}_v\cdot{\mathfrak B}_{\varepsilon_{1,n+1}(\omega_0y)}
      \cdot q_{n+1})
$$
Since
$\varepsilon_{1,n+1}(\omega_0y)=\omega_0(\varepsilon_{1,\ol{n+1}}(y))$
and $\ell(\varepsilon_{1,\ol{n+1}}(y))=\ell(y)+n+1$, 
the result follows by the Pieri-type formula (Theorem~D).
\QED

\begin{lem}\label{lem:reduction}
Suppose $u<_0 w$ in ${\mathcal B}_{n+1}$ and $u(p)=w(p)=q$ for some 
$1\leq p\leq n+1$. 
Then
\begin{enumerate}
\item[(1)] $u/_p<_0 w/_p$ and $\ell(w)-\ell(u) = \ell(w/_p)-\ell(u/_p)$.
\item[(2)] For any opposite isotropic flags $\Edot,\Epdot$ in $V$,
    $$
     \psi_p\left( Y_{u/_p}\Edot \bigcap Y_{\omega_0 w/_p}\Epdot\right)
     \ =\ Y_u\psi_{\ol{q}}\Edot \bigcap Y_{\omega_0 w}\psi_q\Epdot.
    $$
\end{enumerate}
\end{lem}

\noindent{\bf Proof. }
Since $u<_0w$ and $u(p)=w(p)$, Proposition~\ref{prop:thmA} implies that
$u/_p<_0 w/_p$.
Moreover, $wu^{-1}$ is shape equivalent to $w/_p(u/_p)^{-1}$, so the first
statement follows from Theorem~B(1).

For (2), Proposition~\ref{prop:subsets} gives the inclusion $\subset$.
By Corollary~\ref{cor:intersection},
$\psi_{\ol{q}}\Edot$ and $\psi_q\Epdot$ are opposite flags.
Thus, by Proposition~\ref{prop:transverse} and (1), both sides are
irreducible and have the same dimension, proving their equality.
\QED

\begin{thm}\label{thm:elem_equiv}
Suppose $u<_0 w$ in ${\mathcal B}_{n+1}$ and $u(p)=w(p)=q$ for some 
$1\leq p\leq n+1$. 
Then for any strict partition $\lambda$, we have
$$
b^w_{u\,\lambda}\ =\ b^{w/_p}_{u/_p\,\lambda}
\qquad\mbox{and}\qquad
c^w_{u\,\lambda}\ =\ c^{w/_p}_{u/_p\,\lambda}.
$$
\end{thm}

\noindent{\bf Proof. }
We first study the map 
$\Psi : \Spn/P_0 \hookrightarrow \mbox{\it Sp}_{2n+2}{\mathbb C}/P_0$, 
defined by $K \mapsto \Span{K,h}$.
Here, ${\it Sp}_{2n+2}{\mathbb C}/P_0$ is the Grassmannian of 
maximal isotropic subspaces of $V\oplus H_2$.
If $\Edot,\Epdot$ are opposite isotropic flags in $V$, the analog of 
Corollary~\ref{cor:intersection} is
$$
\Psi (\Upsilon_\lambda\Edot)\ =\ 
\Upsilon_\lambda \psi_{n+1}\Edot \bigcap \Upsilon_{(n+1)}\psi_{\ol{n+1}}\Epdot,
$$
where $(n{+}1)$ is a decreasing sequence of length 1.
We leave this to the reader.
As this intersection is generically transverse, 
$\Psi_* Q_\lambda = Q_\lambda \cdot q_{n+1}$.
Thus 
$$
\Psi^* Q_\lambda \ =\ \left\{\begin{array}{lcl}
 Q_\lambda   &\quad&\lambda_1<n+1\\
0&&\lambda_1=n+1
\end{array}\right..
$$
To see this, note that $\Psi^* Q_\lambda=\sum_\mu d^\mu_\lambda Q_\mu$,
where
\begin{eqnarray*}
d^\mu_\lambda &:=& \deg((\Psi^* Q_\lambda)\cdot  Q_{\mu^c})\\
&=& \deg(Q_\lambda\cdot \Psi_*(Q_{\mu^c}))\\
&=& \deg(Q_\lambda\cdot q_{n+1} \cdot  Q_{\mu^c})
\quad=\quad \delta^\mu_\lambda,
\end{eqnarray*}
the Kronecker delta, by the Pieri formula for isotropic
Grassmannians~\cite{Hiller_Boe}.

Consider the commutative diagram:
$$
\begin{picture}(305,67)
\put(4,5){${\mathcal Y}^{w/_p}_{u/_p}\ =\ 
    \pi(Y_{u/_p}\bigcap Y'_{\omega_0{w/_p}})$}
\put(63.8,53){$Y_{u/_p}\bigcap Y'_{\omega_0{w/_p}}$}
\put(200,5){$\pi(Y_u\bigcap Y'_{\omega_0w})\ =\ {\mathcal Y}^w_u$}
\put(211.52,53){$Y_u\bigcap Y'_{\omega_0w}$}
\put(140,7.5){\vector(1,0){50}}    	\put(137,56){\vector(1,0){60}}
\put(160,61){$\psi_p$}			\put(160,10){$\Psi$}
\put(91.5,46){\vector(0,-1){26}}     	\put(230.5,46){\vector(0,-1){26}}
\put(81,32){$\pi$}  			\put(235,32){$\pi$}  
\end{picture}
$$
Thus $[{\mathcal Y}^w_u]=\Psi_*[{\mathcal Y}^{w/_p}_{u/_p}]$
and the maps $\pi$ have the same degree $\delta$ as the horizontal maps are 
isomorphisms.

Let $\lambda$ be a strict partition. 
Then
$c^w_{u\,\lambda}= \delta\cdot \deg (Q_{\lambda}\cdot [{\mathcal Y}^w_u])$
which is 
$$
\delta\cdot \deg (Q_{\lambda}\cdot\Psi_*[{\mathcal Y}^{w/_p}_{u/_p}])
= \delta\cdot \deg (\Psi^*(Q_{\lambda})\cdot 
           [{\mathcal Y}^{w/_p}_{u/_p}])
= \delta\cdot\deg(Q_{\lambda}\cdot[{\mathcal Y}^{w/_p}_{u/_p}])
= c^{w/_p}_{u/_p\,\lambda}.\qquad
\QED
$$

\begin{lem}\label{lem:identity}
Suppose $u<_0w$ and $x<_0 z$ in ${\mathcal B}_n$ with $w u^{-1}=z x^{-1}$
and $u(i)\neq w(i)$ for all $i\in[n]$.
Then, if\/ $Y_u,Y_x$, (respectively $Y'_{\omega_0w}, Y'_{\omega_0{z}}$)
are defined with the same
flags, there is commutative diagram
$$
\begin{picture}(165,70)(8.5,0)
\put(14,5){$\pi(Y_u\bigcap Y'_{\omega_0w})$}
\put(8.5,53){$Y_u\bigcap Y'_{\omega_0w}$}
\put(99.2,5){$\pi(Y_x\bigcap Y'_{\omega_0{z}})$}
\put(128,53){$Y_x\bigcap Y'_{\omega_0{z}}$}
\put(85,5){$=$}
\put(63,56){\line(1,0){6}}    \put(75,56){\line(1,0){6}}
\put(87,56){\line(1,0){6}}    \put(99,56){\line(1,0){6}}    
\put(111,56){\vector(1,0){8}}
\put(87.5,61){$f$}
\put(20,28){$\pi$}	\put(145,28){$\pi$}
\put(27.8,46){\vector(1,-2){13}} \put(146,46){\vector(-1,-2){13}}
\end{picture}
$$
where $f$ is an isomorpism between Zariski open subsets of
$Y_u\bigcap Y'_{\omega_0w}$ and $Y_x\bigcap Y'_{\omega_0{z}}$, 
and the maps $\pi$ have degree 1.
\end{lem}

\noindent{\bf Proof of Theorem~B(2). }
Let $u,w,x,z\in{\mathcal B}_n$ with $w u^{-1}$ shape equivalent to 
$z x^{-1}$ and $\lambda$ a strict partition.
If $u(p)=w(p)$, then $u/_p<_0w/_p$ with 
$w/_p(u/_p)^{-1}$ shape equivalent to $w u^{-1}$.
Thus we may assume that the hypotheses of
Lemma~\ref{lem:identity} hold, by
Theorem~\ref{thm:elem_equiv}.
But then 
$\pi_*[Y_u\bigcap Y'_{\omega_0 w}]=
\pi_*[Y_x\bigcap Y'_{\omega_0 z}]$, hence
$\pi_*({\mathfrak C}_u\cdot{\mathfrak C}_{\omega_0 w})=
\pi_*({\mathfrak C}_x\cdot{\mathfrak C}_{\omega_0 z})$,
showing $c^w_{u\,\lambda}=c^z_{x\,\lambda}$.
\QED

\noindent{\bf Proof of Lemma~\ref{lem:identity}. }
Let $g_{\ol{\jmath}}$ for $1\leq j\leq n$ be the functions of
Lemma~\ref{lem:param}, defined for all 
$\Fdot\in Y^\circ_u\bigcap {Y'}^\circ_{\omega_0w}$.
For such $\Fdot$, let $f(\Fdot)$ be the flag whose $\ol{\jmath}$-th subspace
is 
$$
\Span{g_{u^{-1}x({\ol{n}})}(\Fdot),\ldots,g_{u^{-1}x({\ol{\jmath}})}(\Fdot)}.
$$
Then $f$ is an isomorphism between the intersections of Schubert cells,
which are Zariski dense in the intersections of Schubert varieties,
and the diagram is commutative.
Since we may assume $w$ is a Grassmannian permutation, 
and in this case, the map 
$\pi: Y_{\omega_0w}\rightarrow \pi(Y_{\omega_0w})$
has degree 1, it follows that the maps $\pi$
(which have the same degree) have degree 1.
\QED

\section{More identities}\label{sec:more_identities}

\subsection{Product decomposition}
Suppose that $W$ is a $2m$-dimensional complex symplectic vector space
and consider the map
$$
\Xi\ : \ \mbox{\it Sp}_{2m}{\mathbb C}/P_0  \times \Spn/P_0\ 
\longrightarrow\ 
\mbox{\it Sp}_{2m+2n}{\mathbb C}/P_0
$$
defined by
$$ 
\Xi\ : \ (H,K)\ \longmapsto \ H\oplus K,
$$
where  $H\subset W$ and $K\subset V$ are maximal isotropic (Lagrangian) 
subspaces.

\begin{thm}
Let $\eta,\zeta\in{\mathcal B}_{n+m}$ with $\eta\cdot\zeta$ a disjoint
product and $\#\supp(\eta)\leq m$, 
$\#\supp(\zeta)\leq n$.  
Then for any $\eta'\in{\mathcal B}_m$ and $\zeta'\in{\mathcal B}_n$
with $\eta\sim \eta'$ and  $\zeta\sim \zeta'$  ($\sim$ is shape
equivalence), there is an element $g$ of 
$\mbox{\it Sp}_{2m+2n}{\mathbb C}$ such that 
$$
\Xi ({\mathcal Y}_{\eta'} \times {\mathcal Y}_{\zeta'})\ =\ 
g ({\mathcal Y}_{\eta\cdot\zeta}).
$$
\end{thm}

\noindent{\bf Proof. }
This is a consequence of Lemma~5.2.1 of~\cite{BS98a}, 
the analogous fact for the classical flag manifold and Grassmannian.
Restricting that result to the symplectic flag manifold and Grassmannian of
maximal isotropic subspaces proves the theorem.
\QED

\noindent
{\bf Remark. } This does not hold for the odd orthogonal case.
In fact, even the map $\Xi$ cannot be defined:
If $V,W$ are odd-orthogonal spaces, then $V\oplus W$ is an even-dimensional
space.

\subsection{Some hidden symmetries}\label{sec:cyclic}
Define the permutation $\rho\in {\mathcal B}_n$ by
$\rho(i) = i{-}1{-}n$.
For example, in ${\mathcal B}_4$, we have 
$\rho=\Span{1,\,\ol{4}}\Span{2,\,\ol{3}}$, in
the cycle notation of Section~\ref{sec:orders}.
In general, $\rho=\Span{1,\,\ol{n}}\Span{2,\;\ol{n{-}1}}\cdots
\langle \lfloor\frac{n{+}1}{2}\rfloor,\ 
    \lfloor\ol{\frac{n{+}2}{2}}\rfloor  \rangle$.
For a permutation $\zeta\in{\mathcal B}_n$, let 
$\zeta^\rho:=\rho\zeta\rho$ denote the conjugation of $\zeta$ by $\rho$.

A (saturated) chain in either the Lagrangian order or Lagrangian r\'eseau
gives a sequence $a_1,\ldots,a_m$ of labels.
The {\it peak set} of that chain consists of those $j\in\{2,3,\ldots,m-1\}$
with $a_{j-1}<a_j>a_{j+1}$, its {\it descent set} of those
$i\in\{1,2,\ldots,m-1\}$ with $a_i>a_{i+1}$, and 
its {\it ascent set}  those $i$ with $a_i<a_{i+1}$.
\medskip

\noindent{\bf Theorem~C(1). }
{\it 
Let $\zeta\in{\mathcal B}_n$.
Then $\Len(\zeta)=\Len(\zeta^\rho)$ and if $\lambda$ 
is a strict partition with $|\lambda|=\Len(\zeta)$,
then 
$$
   b^\zeta_\lambda = b^{\zeta^\rho}_\lambda\qquad
   \mbox{and}\qquad
   c^\zeta_\lambda = c^{\zeta^\rho}_\lambda.
$$
}\medskip

\noindent{\bf Proof. }
The first statement follows from Corollary~\ref{cor:compare_orders}
as $||\zeta||=||\zeta^\rho||$ in 
$({\mathcal S}_{\pm[n]},\prcc)$ and both $\zeta$ and $\zeta^\rho$ have
the same number of sign changes.
The last statement is a consequence of Lemma~\ref{lem:cyclic} below,
which shows $\pi_*({\mathfrak C}_u\cdot {\mathfrak C}_{\omega_0 \zeta u})
=\pi_*({\mathfrak C}_x\cdot {\mathfrak C}_{\omega_0\zeta^\rho x})$
whenever $u<_0\zeta u$ and $x<_0 \zeta^\rho x$.
\QED

\begin{ex}\label{ex:cyclic}
In ${\mathcal B}_4$, let $\zeta=\Span{134}\Perm{2}$.
Then $\zeta^\rho=\Span{142}\Perm{3}$.
Consider the intervals $[e,\zeta]_\prec$ and 
$[e,\zeta^\rho]_\prec$ in the labeled Lagrangian r\'eseau 
displayed in Figure~\ref{fig:cyclic}.
\begin{figure}[htb]
$$
\epsfxsize=4.7in \epsfbox{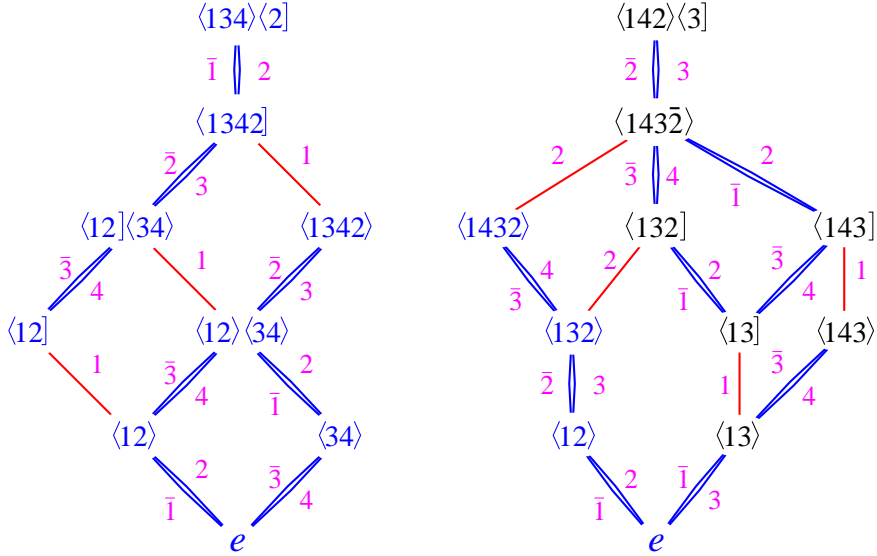}
$$
\caption{Conjugation by $\rho$ on labeled intervals\label{fig:cyclic}}
\end{figure}
While they are not isomorphic, they have the same rank,
the same number of maximal chains, $80$, and the underlying orders each
have 5 chains.
Moreover, they each have 2 chains with peak set $\{3\}$, and one each with
peak sets $\{2\},\{4\}$, and $\{2,4\}$.
The r\'eseaux have the same number of chains with fixed descent sets.
The $j$th component of the following vector  records the number of chains
with descent set equal to the position of the 1's in the binary
representation of $j{-}1$: 
$$
 (0,2,6,4,6,12,8,2,2,8,12,6,4,6,2,0)
$$
\end{ex}

\begin{defn}\label{def:Etdot}
Let $\Edot,\Epdot$ be opposite isotropic flags in $V$.
Define a flag $\Etdot$ by:
$$
\Etdot\ :\ 
E_1\bigcap E'_{\ol{1}}\ \subset\  \cdots\
 \subset\  E_n\bigcap E'_{\ol{1}}\ \subset\  
(E_{\ol{n}}+E'_{\ol{1}})\ \subset\ \cdots\ \subset\ 
(E_{\ol{2}}+E'_{\ol{1}})\ \subset\  V.
$$
Define $\Etpdot$ the same way, but with the roles of $\Edot$ and $\Epdot$
reversed. 
This gives opposite flags $\Etdot, \Etpdot$, and since 
$(A\bigcap B)^\perp =(A^\perp+B^\perp)$, they are isotropic.
\end{defn}

\begin{lem}\label{lem:cyclic}
Suppose $u,w,x,z\in{\mathcal B}_n$ with 
$\rho u^{-1}w\rho=x^{-1}z$, and $u(j)\neq w(j)$ for $1\leq j\leq n$.
Then, for any opposite isotropic flags $\Edot,\Epdot$ in 
$V$, there is a commutative diagram
$$
\begin{picture}(208,72)(-1.5,0)
\put(4,5){$\pi(Y_u\Edot\bigcap Y_{\omega_0w}\Epdot)$}
\put(-1.5,53){$Y_u\Edot\bigcap Y_{\omega_0w}\Epdot$}
\put(119.2,5){$\pi(Y_x\Etdot\bigcap Y_{\omega_0{z}}\Etpdot)$}
\put(146,53){$Y_x\Etdot\bigcap Y_{\omega_0{z}}\Etpdot$}
\put(102,5){$=$}
\put(81,56){\line(1,0){5}}    \put(91,56){\line(1,0){5}}    
\put(101,56){\line(1,0){5}}   \put(111,56){\line(1,0){5}}    
\put(121,56){\line(1,0){5}}    \put(131,56){\vector(1,0){8}}
\put(105,61){$f$}
\put(23,28){$\pi$}	\put(175,28){$\pi$}
\put(30.8,46){\vector(1,-2){13}} \put(176,46){\vector(-1,-2){13}}
\end{picture}
$$
with $f$ an isomorpism of Zariski open subsets of\/
$Y_u\Edot\bigcap Y_{\omega_0w}\Epdot$ and 
$Y_x\Etdot\bigcap Y_{\omega_0{z}}\Etpdot$.
\end{lem}

\noindent{\bf Proof. }
Let $\Gdot,\Gpdot$ be oppsite (not necessarily isotropic)
flags in $V$.
Define $\Gdot\!^+$ to be
$$
\Gdot\!^+\ :\ G_{\ol{n-1}}\bigcap G'_{n-1}\  \subset \ 
G_{\ol{n-2}}\bigcap G'_{n-1}\  \subset\  \cdots\  \subset \ 
G'_{n-1}\  \subset\  V.
$$
Define $\Gpdot^+$ to be
$$
\Gpdot^+\ :\ G_{\ol{n}}\ \subset\  (G_{\ol{n}}+G'_{\ol{n}})\ \subset\ 
\cdots\  \subset\  (G_{\ol{n}}+G'_{n-1})\ \subset\  V.
$$
For $\zeta\in {\mathcal S}_{\pm [n]}$, let $\zeta^+$ be
the conjugation of
$\zeta$ by the cycle $(\ol{n},\ldots,\ol{1},1,\ldots,n)$.
In Section~5.3 of \cite{BS98a}, the following proposition is proven:

\begin{prop}\label{prop:cyclic}
Let $u,w,x,z\in{\mathcal S}_{\pm[n]}$ with $u\vartriangleleft_0 w$,
$x\vartriangleleft_0 z$, $(u^{-1}w)^+=x^{-1}z$, and $w$ is a Grassmannian
permutation with descent $0$, 
($w(\ol{n})<\cdots<w(\ol{1})$ and $w(1)<\cdots<w(n)$).
If $\pi:{\mathbb F}\ell V\twoheadrightarrow {\bf G}_n(V)$ is the projection,
then there is a commutative diagram:
$$
\begin{picture}(226,70)(-3.5,0)
\put(3,5){$\pi(X_u\Gdot\bigcap X_{\omega_0w}\Gpdot)$}
\put(-3.5,53){$X_u\Gdot\bigcap X_{\omega_0w}\Gpdot$}
\put(119.2,5){$\pi(X_x\Gdot\!^+\bigcap X_{\omega_0{z}}\Gpdot^+)$}
\put(148,53){$X_x\Gdot\!^+\bigcap X_{\omega_0{z}}\Gpdot^+$}
\put(105,5){$=$}
\put( 83,56){\line(1,0){5}}    \put( 93,56){\line(1,0){5}}   
\put(103,56){\line(1,0){5}}    \put(113,56){\line(1,0){5}}    
\put(123,56){\line(1,0){5}}    \put(133,56){\vector(1,0){8}}
\put(107,61){$f$}
\put(26,28){$\pi$}	\put(183,28){$\pi$}
\put(31.8,46){\vector(1,-2){13}} \put(186,46){\vector(-1,-2){13}}
\end{picture}
$$
with $f$ an isomorpism of Zariski open subsets of\/
$X_u\Gdot\bigcap X_{\omega_0w}\Gpdot$ and 
$X_x\Gdot\!^+\bigcap X_{\omega_0{z}}\Gpdot^+$.
\end{prop}

It suffices to prove Lemma~\ref{lem:cyclic} with $w$ a Grassmannian
permutation, by Lemma~\ref{lem:identity}. 
Observe that $(\Etdot,\Etpdot)$ is the result of $n$ applications of the
map $(\Edot,\Epdot)\mapsto (\Edot^+,\Epdot^+)$.
Similarly, 
$\rho=(\ol{n},\ldots,\ol{2},\ol{1},1,\ldots,n)^n$.
Thus, iterating Proposition~\ref{prop:cyclic} $n$ times
gives the commutative diagram in ${\mathbb F}\ell V$ and 
${\bf G}_n V$:
$$
\begin{picture}(212,70)(-3,0)
\put(2.5,5){$\pi(X_u\Edot\bigcap X_{\omega_0w}\Epdot)$}
\put(-3,53){$X_u\Edot\bigcap X_{\omega_0w}\Epdot$}
\put(119.2,5){$\pi(X_x\Etdot\bigcap X_{\omega_0{z}}\Etpdot)$}
\put(146,53){$X_x\Etdot\bigcap X_{\omega_0{z}}\Etpdot$}
\put(105,5){$=$}
\put(83,56){\line(1,0){5}}    \put(93,56){\line(1,0){5}}    
\put(103,56){\line(1,0){5}}   \put(113,56){\line(1,0){5}}    
\put(123,56){\line(1,0){5}}    \put(133,56){\vector(1,0){8}}
\put(110,61){$f$}
\put(25.5,28){$\pi$}	\put(178,28){$\pi$}
\put(32.3,46){\vector(1,-2){13}} \put(179,46){\vector(-1,-2){13}}
\end{picture}
$$
Restricting this to the subset of isotropic flags gives 
the diagram of the lemma.
\QED
 
These same arguments prove the analog of
Lemma~\ref{lem:cyclic} for $\Son/C$.

\subsection{More hidden symmetries}\label{sec:hidden}
Until now, we have deduced identities in $H^*\Spn/B$ by rectricting
constructions involving Schubert subvarieties of ${\mathbb F}\ell(V)$
to those in $\Spn/B$ via the embedding 
$\Spn/B\hookrightarrow{\mathbb F}\ell(V)$.
This is the geometric counterpart of the embedding 
${\mathcal B}_n\hookrightarrow{\mathcal S}_{\pm[n]}$ studied in Section~2.
Here, we explore the geometry of the map 
$\iota:{\mathcal S}_n\hookrightarrow{\mathcal B}_n$, where a permutation 
$w\in{\mathcal S}_n$ is extended to act on $-[n]$ by $w(-i)=-w(i)$.
An immediate consequence of Definition-Theorem~\ref{defth:2.6} and its analog
for $\prcc$ is the following lemma.

\begin{lem}\label{lemma:iota}
The map $\iota$ is an embedding of Bruhat orders
$({\mathcal S}_n,\vartriangleleft)\hookrightarrow({\mathcal B}_n,\leq)$
and it respects the length functions in each order.
Furthermore, $\iota({\mathcal S}_n)$ consists of those permutations 
$\zeta\in{\mathcal B}_n$ with $a\cdot \zeta(a)>0$ for all $a$.
\end{lem}

Let $L, L^\perp$ be complementary Lagrangian subspaces in $V$.
The pairing $(x,y)\in L\oplus L^\perp\mapsto \beta(x,y)$,
where $\beta$ is the alternating form, identifies them as linear duals.
Given a subspace $H$ of $L$, let $H^\perp\subset L^\perp$ denote its
annihilator in $L^\perp$.
Then $H+H^\perp$ is a Lagrangian subspace of $V$.

Let ${\mathbb F}\ell(L)$ be the space of complete flags 
$\Fdot:=F_1\subset F_2\subset \cdots \subset F_n=L$ in $L$.
Note that here $\dim F_i=i$.
For each $k=0,1,\ldots,n$, define an injective map
$$
\varphi_k\ :\ {\mathbb F}\ell(L)\  \longrightarrow\  \Spn/B
$$
by
\begin{equation}\label{eq:fln-embedd}
  (\varphi_k\Fdot)_{\ol{\jmath}}\ =\ \left\{\begin{array}{lcl}
     F_{n+1-j}&\quad&j\geq n-k+1\\
     F_k+F_{k+j-1}^\perp&& j\leq n-k+1
  \end{array}\right..
\end{equation}
Then $(\varphi_k\Fdot)_{\ol{1}}=(F_k+F_k^\perp)$ is Lagrangian, showing
that $\varphi_k\Fdot$ is an isotropic flag.

For $w\in {\mathcal S}_n$ 
the Schubert variety $X_w\Edot$ of ${\mathbb F}\ell(L)$ 
consists of those flags $\Fdot\in{\mathbb F}\ell(L)$ satisfying
\begin{equation}\label{eq:small-flags}
         \dim E_a \bigcap F_b =\#\{b\geq l \mid w(l)+a\geq n+1\} .
\end{equation}
We determine the image of Schubert varieties of ${\mathbb F}\ell(L)$
under these maps $\varphi_k$.

Define $\epsilon_k:{\mathcal S}_n\rightarrow{\mathcal B}_n$ by
\begin{equation}\label{eq:perm-embedd}
   (\epsilon_kw)(j)\ =\ \left\{\begin{array}{lcl}
     w(j+k)&\quad& 1\leq j\leq n-k\\
     \ol{w(n{+}1{-}j)}&& n-k<j\leq n \rule{0pt}{17pt}
    \end{array}\right.
\end{equation}
Note that $\epsilon_0w = \iota(w)$.

\begin{lem}\label{lem:epsilon_k}
Let $u,w\in{\mathcal S}_n$ with $u\vartriangleleft_k w$.
Then $\epsilon_k$ induces an isomorphism of graded posets
$[u,w]_{\vartriangleleft_k} \stackrel{\sim}{\longrightarrow}[u,w]_0$
and $\iota(wu^{-1})= \epsilon_k w(\epsilon_k u)^{-1}$.
\end{lem}

\noindent{\bf Proof. }
A consequence of the definitions is that, for $u,w\in{\mathcal S}_n$,
$$
u\vartriangleleft_kw \ \Longleftrightarrow\ 
\epsilon_k u <_0 \epsilon_k w,
$$
and $\iota(wu^{-1})= \epsilon_k w(\epsilon_k u)^{-1}$.
The Lemma follows immediately from these observations.
\QED

\begin{cor}
The map $\iota$ is an embedding of ranked orders 
$\iota:({\mathcal S}_\infty,\prcc)\ \hookrightarrow\ 
   ({\mathcal B}_\infty,\prec)$.
\end{cor}

Let $w^\vee$ be defined by $w^\vee(j)=n+1-w(j)$.
Then $X_w\Edot$ and $X_{w^\vee}\Epdot$ are dual under the intersection
pairing, where $\Edot,\Epdot$ are opposite flags.

\begin{lem}\label{lem:hidden_subset}
With these definitions, $\varphi_k X_w\Edot$ is a subset of either
$$
Y_{\epsilon_kw}\varphi_n\Edot \quad\mbox{\rm or}\quad
Y_{\omega_0\epsilon_kw^\vee}\varphi_0\Edot .
$$
\end{lem}

\noindent{\bf Proof. }
Let $\Fdot\in X_w\Edot$.
We show $\varphi_k\Fdot\in Y_{\epsilon_kw}\varphi_n\Edot$,
that is, for each $-n\leq i\leq n$ and $1\leq j\leq n$ ($i\neq 0$),
\begin{equation}\label{eq:key-identity}
  \dim \left(\varphi_n\Edot\right)_i \bigcap
  \left(\varphi_k\Fdot\right)_{\ol{\jmath}}\ \geq\ 
  \#\{ n\geq l\geq j\mid i\geq \epsilon_kw(l)\}.
\end{equation}

Suppose that $j>n-k+1$. 
Then $(\varphi_k\Fdot)_{\ol{\jmath}}=F_{n+1-j}\subset L=
      (\varphi_n \Edot)_{\ol{1}}$.
If $n\geq l\geq j$, then $(\epsilon_k w)(l)=\ol{w(n{+}1{-}l)}<0$.
Thus if $i>0$,~(\ref{eq:key-identity}) holds as both sides equal $n+1-j$.
Suppose $i<0$.
Then $(\varphi_n\Edot)_i=E_{n+1-\ol{\imath}}$ and so the left side
of~\ref{eq:key-identity}) is 
\begin{eqnarray*}
 \dim E_{n+1-\ol{\imath}}\bigcap F_{n+1-j}&\geq&
 \#\{ m\leq n+1-j\mid w(m)+n+1-\ol{\imath} \geq n+1\}\\
 &=&
 \#\{n\geq l\geq j\mid i\geq \ol{w(n{+}1{-}l)} = \epsilon_kw(l)\}.
\end{eqnarray*}

Now suppose that $j\leq n-k+1$.
Then $(\varphi_k\Fdot)_{\ol{\jmath}}=F_k+F_{k+j-1}^\perp$.
Thus the left side of~(\ref{eq:key-identity}) is
\begin{equation}\label{eq:new-dim}
 \dim(\varphi_n\Edot)_i\bigcap F_k \ +\ 
 \dim(\varphi_n\Edot)_i\bigcap F_{k+j-1}^\perp.
\end{equation}
If $i<0$, then $(\varphi_n\Edot)_i\subset L$ and only the first term 
of~(\ref{eq:new-dim}) contributes.
By the previous paragraph, this is
$$
\dim(\varphi_n\Edot)_i\bigcap F_k \ \geq\ 
\#\{n\geq l\geq k\mid i\geq \epsilon_kw(l)\}.
$$
If $k\geq l$, then $\epsilon_kw(l)>0>i$, showing this equals the right side
of~(\ref{eq:key-identity}).

If now $i>0$, then $(\varphi_n\Edot)_i = L+E_{n-i}^\perp$.
Thus~(\ref{eq:new-dim}) is
$$
k+\dim E_{n-i}^\perp\bigcap F_{k+j-1}^\perp\ =\ 
k+n- \dim(E_{n-1}+F_{k+j-1}).
$$
But this is 
$k+n-\dim E_{n-1} -\dim F_{k+j-1}+\dim E_{n-1}\bigcap F_{k+j-1}$,
which is at least
\medskip

\noindent
${\displaystyle
i-j+1+\#\{k+j-1\geq m\geq1\mid w(m)+n-1\geq n+1\}
}$\vspace{-5pt}
\begin{eqnarray*}
\hspace{2.2in}&=&
n-j+1-\#\{n\geq m\geq k+j\mid w(m)\geq i+1\}\\
&=& k+\#\{n-k\geq l\geq j\mid w(l+k)\leq i\}\\
&=& k+\#\{n-k\geq l\geq j\mid \epsilon_kw(l)\leq i\}.
\end{eqnarray*}
This equals the right side of~(\ref{eq:key-identity}) 
since $l>n-k$ implies $\epsilon_kw(l)<0<i$.

Similar arguments show 
$\varphi_k\Fdot\in Y_{\omega_0\epsilon_kw^\vee}\varphi_0\Edot$.
\QED

\begin{cor}
Let $u,w\in{\mathcal S}_n$ with $u\vartriangleleft_k w$ and
$\Edot,\Epdot\in {\mathbb F}\ell(L)$ be opposite flags.
Then $\varphi_n\Edot,\varphi_0\Epdot$ are opposite isotropic flags, and
$$
\varphi_k\left(  X_u\Edot \bigcap X_{w^\vee}\Epdot\right)\ = \
Y_{\epsilon_k u}\varphi_n\Edot 
    \bigcap Y_{\omega_0\epsilon_k w}\varphi_0\Epdot.
$$
\end{cor}

\noindent{\bf Proof. }
Lemma~\ref{lem:hidden_subset} gives the inclusion $\subset$
and it is easy to see that $\varphi_n\Edot,\varphi_0\Epdot$ are opposite.
By Lemma~\ref{lem:epsilon_k}, both sides have the same dimension,
proving equality.
\QED

For each $k=1,2,\ldots,n$, let 
$\pi_k:{\mathbb F}\ell(L)\rightarrow {\bf G}_k(L)$ be the projection induced
by $\Edot \mapsto E_k$.
As in Lemma~\ref{lem:identity}, if 
$u\vartriangleleft_k w$ in ${\mathcal S}_n$ and $\Edot,\Epdot$ are opposite
flags in ${\mathbb F}\ell(L)$, then the intersection
$X_u\Edot\bigcap X_{w^\vee}\Epdot$ is mapped birationally onto its image 
$\pi_k(X_u\Edot\bigcap X_{w^\vee}\Epdot)$ in ${\bf G}_k(L)$.
Furthermore, the image cycle depends only upon $\eta:=wu^{-1}$.
Denote it by ${\mathcal X}_\eta$.

Define $\Phi_k: {\bf G}_k(L) \rightarrow \Spn/B$ by
$H\mapsto (H+H^\perp)$.
Then $\Phi_k\circ\pi_k = \pi\circ \varphi_k$ and we have the following
corollary.

\begin{cor}\label{cor:hidden-identity}
For any $\eta\in{\mathcal S}_n$, 
${\mathcal Y}_{\iota(\eta)}=\Phi_k({\mathcal X}_\eta)$, where 
$k=\#\{a\mid a<\eta(a)\}$.
\end{cor}

Recall that $\gamma:= \iota (1,2,\ldots,n) \in {\mathcal B}_n$.
\medskip

\noindent{\bf Theorem~C(2). }
{\it 
Let $\zeta\in \iota({\mathcal S}_n)$.
Then $\Len(\zeta)=\Len(\zeta^\gamma)$ and if $\lambda$ is a strict partition
with $|\lambda|=\Len(\zeta)$, then 
$$
  b^\zeta_\lambda = b^{\zeta^\gamma}_\lambda\qquad
  \mbox{and}\qquad
  c^\zeta_\lambda = c^{\zeta^\gamma}_\lambda.
$$
}\medskip

\noindent{\bf Proof. }
The first statement follows from Lemma~\ref{lemma:iota} as 
$||\eta||=||\eta^\gamma||$ for $\eta\in{\mathcal S}_n$.
The last statement is a consequence of the identity 
${\mathcal X}_\eta={\mathcal X}_{\eta^\gamma}$
(Proposition~\ref{prop:cyclic}) and Corollary~\ref{cor:hidden-identity}.
\QED

\begin{ex}\label{ex:no-iso}
Let $\eta=(1,2,4,3)$.  
Then $\zeta=\iota(\eta)= \Span{1,2,4,3}$ and $\zeta^\gamma=\Span{1,4,2,3}$.
The labeled intervals $[e,\eta]_{\prcs}$ and $[e,\zeta]_\prec$ are
isomorphic.
Consider the intervals $[e,\zeta]_\prec$ and $[e,\zeta^\gamma]_\prec$
in the labeled Lagrangian r\'eseau displayed in Figure~\ref{fig:hidden}.
While they are not isomorphic, they have the same rank, the same number of
maximal chains, 16, and the underlying orders each have 2 maximal chains.
Moreover, they each have a peakless chain and one with peak set $\{2\}$.
The r\'eseaux each have 2 increasing
chains, 2 decreasing chains, 6 with descent set $\{1\}$, and 6 with descent
set $\{2\}$. 
\begin{figure}[htb]
$$\epsfxsize=4.5in \epsfbox{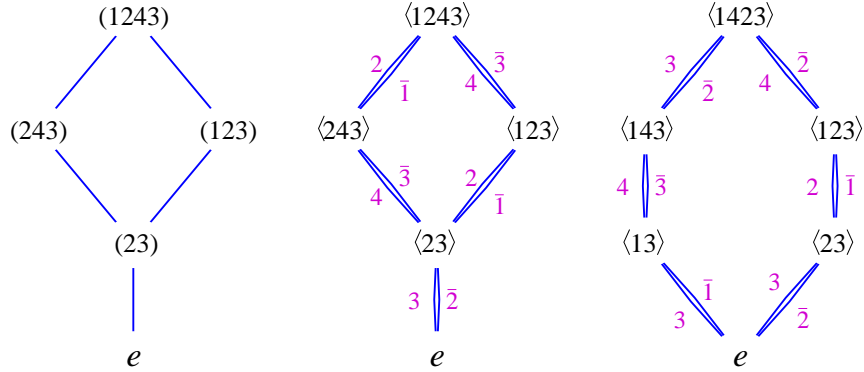}$$
\caption{Conjugation by $\gamma$ on labeled intervals\label{fig:hidden}}
\end{figure}
\end{ex}

\section{Minimal permutations and labeled r\'eseaux}

\subsection{Minimal permutations}\label{sec:minimal}

For a cycle $\zeta\in{\mathcal B}_\infty$, let 
$\delta(\zeta)=1$ if $\zeta$ has the form  $\iota(\eta)$
for $\eta\in {\mathcal S}_\infty$ and
$\delta(\zeta)=0$ otherwise.
Note that $\delta(\zeta)=1$ if and only if $a>0$ implies $\zeta(a)>0$.
\medskip

\begin{lem}\label{lem:minimal_cycle}
Let $\zeta\in{\mathcal B}_\infty$ be a cycle.
Then $\Len(\zeta)\ \geq\ \#\supp(\zeta) - \delta(\zeta)$.
\end{lem}

Recall that $s(\zeta)$ counts the number of sign changes in $\zeta$.

\noindent{\bf Proof. }
A saturated chain in $[e,\zeta]_\prcs$ (in ${\mathcal S}_{\pm\infty}$) gives
a factorization of $\zeta$ into transpositions.
If $\zeta$ consists of two cycles in ${\mathcal S}_{\pm\infty}$, then 
$||\zeta||\ \geq \ 2(\#\supp(\zeta) -1)$
and by Corollary~\ref{cor:compare_orders}(2),
$$
\Len(\zeta)\ \geq\ \#\supp(\zeta) -1 + s(\zeta) \ \geq\ 
\#\supp(\zeta) -\delta(\zeta),
$$
with equality only if $s(\zeta)=0$, that is, only if $\delta(\zeta)=1$.

Similarly, if $\zeta$ is a single cycle in ${\mathcal S}_{\pm\infty}$, then 
$||\zeta||\ \geq \ 2\#\supp(\zeta) -1$.
Since $\delta(\zeta)=0$ and $s(\zeta)\geq 1$,  
Corollary~\ref{cor:compare_orders}(2) gives
$$
\Len(\zeta)\ \geq\ \#\supp(\zeta) -\delta(\zeta).
\qquad \QED
$$

\begin{cor}
If $\zeta\in{\mathcal B}_\infty$ is irreducible and 
$\Len(\zeta)=\#\supp(\zeta) - \delta(\zeta)$, then $\zeta$ is a single cycle.
\end{cor}

\noindent{\bf Proof. }
Recall that if $\eta,\xi\in{\mathcal B}_\infty$ have disjoint supports, then
$\Len(\eta\cdot\xi)\geq\Len(\eta)+\Len(\xi)$, with equality only when 
$\eta\cdot\xi$ is a disjoint product.
Thus, by Lemma~\ref{lem:minimal_cycle}, 
$$
\Len(\zeta)\ \geq\ \#\supp(\zeta) - \delta(\zeta),
$$
with equality only when $\zeta$ is a single cycle.
\QED

\begin{defn}
A {\it minimal cycle} is a cycle  $\zeta\in{\mathcal B}_\infty$ for which 
$\Len(\zeta)=\#\supp(\zeta) - \delta(\zeta)$.
For minimal cycles, $s(\zeta)+\delta(\zeta)=1$.
A permutation $\zeta\in{\mathcal B}_\infty$ is {\it minimal} if each of its
irreducible factors are minimal cycles.
\end{defn}

\begin{cor} 
If $\eta,\zeta\in{\mathcal B}_\infty$ with $\eta\prec \zeta$ and $\zeta$ is
minimal, then so is $\eta$
\end{cor}

\noindent{\bf Proof. }
By Theorem~\ref{thm:disj-prod}, we may assume $\zeta$ is irreducible.
Then $\zeta$ is a single cycle and the result follows by 
induction on $\Len(\zeta)$, similar to the proof of
Lemma~\ref{lem:minimal_cycle}.
\QED

\begin{cor}\label{cor:minimal_perm}
Let $\zeta\in{\mathcal B}_\infty$.
If $\zeta=\zeta_1\cdots\zeta_s$ is the factorization of
$\zeta$ into irreducibles, then 
$$
\Len(\zeta)\ \geq\ \#\supp(\zeta) - \sum_i \delta(\zeta_i),
$$
with equality only if $\zeta$ is  minimal.
\end{cor}

\begin{lem}\label{lem:sign-cover}
If\/ $\zeta\in{\mathcal B}_\infty$ is a minimal cycle with
$\delta(\zeta)=0$, then 
there is a unique $a\in\supp(\zeta)$ with $a>0>\zeta(a)=:\ol{\alpha}$.
Furthermore,
$$
\begin{array}{ccl}
1\ \precdot\ t_\alpha\ \preceq\ \zeta&&\mbox{if }\ a\leq\alpha\\
1\ \preceq t_\alpha\zeta\ \precdot\ \zeta&&\mbox{if }\ a\geq\alpha.
\end{array}
$$
\end{lem}

\noindent{\bf Proof. }
Let $\zeta\in{\mathcal B}_\infty$ be a minimal cycle with $\delta(\zeta)=0$.
Then $s(\zeta)=1$, so there is a unique $a>0$ with $\ol{\alpha}:=\zeta(a)$.
We prove the lemma when $a>\alpha$:
If $a=\alpha$, then $\zeta=t_a$, as $\zeta$ is irreducible
and if $a<\alpha$, then replacing $\zeta$ by $\zeta^{-1}$, reduces to the 
case $a>\alpha$.

Suppose $a>\alpha$.
We claim that if $b>a$, then $\zeta(b)>\alpha$.
The lemma follows from this claim.
Indeed, then condition (ii) of 
Definition~\ref{defth:2.6}(1) is satisfied, and hence 
$t_{\alpha}\zeta\prec\zeta$.
Since $a<\alpha$, we have $\supp(t_{\alpha}\zeta)=\supp(\zeta)$.
As $\delta(t_{\alpha}\zeta)=1$, it follows that 
$\Len(t_{\alpha}\zeta)\geq \#\supp(\zeta)-1 = \Len(\zeta)-1$
and thus $t_{\alpha}\zeta\precdot\zeta$.

Let $\zeta$ be a permutation for which the claim does not hold with
$\#\supp(\zeta)$ minimal.
By Algorithm~\ref{alg:chain}, $\eta:=\zeta t_{y\,x}\precdot\zeta$,
where $x$ is maximal in $\supp(\zeta)$ and $y$ is minimal subject to 
$y\leq \zeta(x)<\zeta(y)\leq x$.
Note that $y\neq a$.
Since $a<x$, $\delta(\eta)=0$.
As $\zeta$ is minimal and $\eta\precdot\zeta$, either $\eta$ is irreducible
and $\supp(\eta)\subsetneq\supp(\zeta)$ or else $\eta$ is the disjoint
product of two minimal cycles with $x$ in the support of one and $|y|$ in
the support of the other.

If $y<0$, then $y=\ol{a}$ as $\zeta(y)>0$.
Then $\eta(a)=\ol{\zeta(x)}$ and $\eta(x)=\alpha$, so 
$\supp(\eta)\subset\supp(\zeta)$ and so $\eta$ is irreducible with $x$ in
the support of one component and $a$ in the other.
But $x>a>\alpha=\eta(x)>\ol{\zeta(x)}=\eta(a)$, contradicting disjointness.

Suppose now that $y>0$. 
Then $\zeta(x)>\alpha$, for otherwise $y=\ol{a}<0$.
Thus $x>b$.
Since $b>a>\alpha>\eta(b)$ with $\eta(a)=\ol{\alpha}$,
the component $\eta'$ of $\eta$ whose support contains $a$ also contains $b$
and $\supp(\eta')\subsetneq \supp(\zeta)$ with $\delta(\eta')=0$, 
contradicting the minimality of $\#\supp(\zeta)$.
\QED

\subsection{The Grassmann-Bruhat order on ${\mathcal S}_\infty$}

We develop some additional combinatorics
for the symmetric group ${\mathcal S}_\infty$.
For $k\in{\mathbb N}$, the $k$-Bruhat order  ($\vartriangleleft_k$)
on ${\mathcal S}_{\infty}$ is 
the analog of the $k$-Bruhat order on ${\mathcal B}_{\pm\infty}$.
The interval $[u,w]_k$ in the
$k$-Bruhat order depends only upon $wu^{-1}$, so we define the 
{\it Grassmann-Bruhat  order} on ${\mathcal S}_{\infty}$ 
(the Lagrangian order on ${\mathcal B}_\infty$ is its analog) by
$\eta\prcc\zeta$ if there is a $k\in{\mathbb N}$ and 
$u\in{\mathcal S}_\infty$ with 
$u\vartriangleleft_k\eta u\vartriangleleft_k\zeta u$.
This order is ranked by $||\zeta||=l(\zeta u)-l(u)$ for 
$u\vartriangleleft_k\zeta u$, and it has an independent description:
\medskip

\noindent{\bf Definition 3.2.2 of~\cite{BS98a}}
{\sl
Let $\eta,\zeta\in{\mathcal S}_\infty$.
Then $\eta\prcc\zeta$ if and only if 
\begin{enumerate}
\item[(1)] $a<\eta(a) \Rightarrow \eta(a)\leq\zeta(a)$
\item[(2)] $a>\eta(a) \Rightarrow \eta(a)\geq\zeta(a)$
\item[(3)] If $a<\zeta(a)$ and $b<\zeta(b)$  (respectively 
     $a>\zeta(a)$ and $b>\zeta(b)$)  with $a<b$ and $\zeta(a)<\zeta(b)$,
     then $\eta(a)<\eta(b)$.
\end{enumerate}
}
\medskip

Covers $\eta\pcdot\zeta$ correspond to transpositions
$(\alpha,\,\beta)=\zeta\eta^{-1}$ and we construct a labeled Hasse diagram
for $({\mathcal S}_\infty,\prcc)$, labeling such a cover with the greater of
$\alpha,\beta$.
By Theorem 3.2.3 of~\cite{BS98a}, the map $\eta\mapsto \zeta\eta^{-1}$
induces an order-reversing isomorphism between $[e,\zeta]_\prcs$ and 
$[e,\zeta^{-1}]_\prcs$, preserving the edge labels.
Also, if $P=\{p_1,p_2,\ldots\}\subset {\mathbb N}$ and 
$\varepsilon_P:{\mathcal S}_P \hookrightarrow {\mathcal S}_\infty$ is the
map induced by the inclusion $P\hookrightarrow {\mathbb N}$
(these maps induce shape equivalence), then $\varepsilon_P$ induces an
isomorphism 
$[e,\zeta]_\prcs\stackrel{\sim}{\longrightarrow}
[e,\varepsilon_P(\zeta)]_\prcs$,
preserving the relative order of the edge labels.
Specifically, an edge label $i$ of $[e,\zeta]_\prcs$
is mapped to the label $p_i$ of $[e,\varepsilon_P(\zeta)]_\prcs$.
Lastly, we remark that Algorithm~\ref{alg:chain} restricted to 
${\mathcal S}_\infty$, and with $t_{a\,b}$ replaced by the transposition
$(a,\,b)$, gives a chain in the $\prcc$-order on ${\mathcal S}_\infty$
from $e$ to $\zeta$.

\begin{lem}\label{lem:prec_cover}
Let $\zeta\in{\mathcal S}_\infty$ and suppose that $x$ is maximal 
subject to $x\neq \zeta(x)$.
Then, for any $\alpha$
\begin{enumerate}
\item[(1)]  $(\alpha,\,x)\zeta\prcc \zeta\ \Longrightarrow\ 
	\zeta^{-1}(x)\leq \zeta(x)=\alpha$.
\item[(2)] $(\alpha,\,x)\prcc \zeta\ \Longrightarrow\ 
	\zeta^{-1}(x)\geq \zeta(x)=\alpha$.
\end{enumerate}
\end{lem}

\noindent{\bf Proof. }
For 1, let $\eta:=(\alpha,\,x)\zeta\prcc \zeta$ and set 
$a=\zeta^{-1}(x)$ and $b=\zeta^{-1}(\alpha)$.
Note that $a\neq b$ and $\eta(b)=x$.
We claim that $b=x$ and $a<\alpha$, which will establish 1.

Suppose $b\neq x=\eta(b)$.
Then, by the maximality of $x$, $b<\eta(b)$
and so the definition of $\prcc$ implies
$\eta(b)\leq \zeta(b)=\alpha$.
Since $\alpha<x$, this implies $x<x$, a contradiction.
Suppose now that $a>\alpha=\eta(a)$.
By the definition of $\prcc$, this implies that $\eta(a)\geq\zeta(a)=x$,
and so $a>x$, contradicting the maximality of $x$.

The second assertion follows from the first by applying the anti-isomorphism 
$\eta \mapsto \eta\zeta^{-1}$ between $[e,\zeta]_\prcs$ and 
$[e,\zeta^{-1}]_\prcs$,
$$
e\ \prcc\ (\alpha,\,x)\ \prcc\ \zeta\ 
\Longleftrightarrow\  (\alpha,\,x)\zeta^{-1}\ \prcc\ \zeta^{-1}.
\qquad\QED
$$

A cycle $\zeta\in{\mathcal S}_\infty$ is {\it minimal} if
$||\zeta||+1=\#\supp(\zeta)$. 
A permutation is {\it minimal} if it is the disjoint product of minimal
cycles. 
A maximal chain in an interval $[e,\zeta]_\prcs$ is {\it peakless} if we do
not have $a_{i-1}<a_i>a_{i+1}$ for any $i=2,\ldots,||\zeta||{-}1$, where 
$a_1,\ldots,a_{||\zeta||}$ is the sequence of labels in that chain.

\begin{lem}\label{lem:unique-vee}
Suppose $\zeta\in{\mathcal S}_\infty$ is a minimal cycle.
Then there is a unique peakless chain in the labeled interval 
     $[e,\,\zeta]_\prcs$.
     If $\beta$ is the smallest label in such a chain, then the
     transposition of that cover is $(\alpha,\beta)$ where $\alpha<\beta$
     are the two smallest elements of $\supp(\zeta)$.
\end{lem}

\noindent{\bf Proof. } 
We argue by induction on $||\zeta||$, which we assume is at least 2,
as the case $||\zeta||=1$ is immediate.
Replacing $\zeta$ by a shape equivalent permutation 
if necessary, we may assume that $\supp(\zeta)=[n]$, so that
$||\zeta||=m=n-1$.

Replacing $\zeta$ by $\zeta^{-1}$ would only reverse such a chain, so we may
assume that  
$a:= \zeta^{-1}(n) < b:= \zeta(n)$.
We claim that $(b,\,n)\zeta=\zeta(a,\,n) \pcdot\zeta$.
Given this, the conclusion of the lemma follows.
Indeed, let $\eta:=(b,\,n)\zeta$.
Since $\eta(n)=n$, this is an irreducible minimal permutation in 
${\mathcal S}_{n-1}$.
By the inductive hypothesis, $[e,\eta]_\prcs$ has a unique chain with
labels $\beta_1>\cdots>\beta_k<\cdots<\beta_{n-2}$, and each $\beta_i<n$.
The unique extension of this to a chain in $[e,\zeta]_\prcs$ has
$\beta_{n-2}<\beta_{n-1}=n$.
This is the unique such chain in $[e,\zeta]_\prcs$ as $\eta\pcdot\zeta$ is
the unique terminal cover in $[e,\zeta]_\prcs$ with edge label $n$,
by Lemma~\ref{lem:prec_cover}.
Also note that unless $n=2$, $1\leq a<b$, which proves the second part of 1.

By Algorithm~\ref{alg:chain}, if $y$ is chosen minimal so that
$y\leq \zeta(n)=b<\zeta(y)$, then $\zeta(y,\,n)\pcdot\zeta$.
We show that $y=a$, which will establish the claim and complete the proof.

Suppose $y\neq a$.
Since $a<b<n=\zeta(a)$, the minimality of $y$ implies that $y<a$.
But then $\zeta(y,\,n)$ consists of two cycles $\eta$ and $\eta'$
and we have $\eta(a)=n$ and $\eta'(y)=b$.
Since $y<a<b<n$, these cycles are not disjoint, so we have
\begin{eqnarray*}
n-2\ =\  ||\eta\cdot\eta'||&>& ||\eta||+||\eta'||\\
&\geq& \#\supp(\eta)-1 + \#\supp(\eta')-1\ =\ n-2,
\end{eqnarray*}
a contradiction.
\QED

\subsection{The labeled Lagrangian order}\label{sec:lagr-or}
The labeled Lagrangian and 0-Bruhat orders on 
${\mathcal B}_\infty$ are obtained from the Hasse diagrams of the underlying
orders by labeling each cover with the integer $\beta$, where that cover is
either $\zeta\precdot t_\beta \zeta$ ($u \precdot t_\beta u$) or 
$\zeta\precdot t_{\alpha\,\beta}\zeta$ ($u\precdot t_{\alpha\,\beta}u$).
Recall that the map 
$\iota:{\mathcal S}_\infty\rightarrow{\mathcal B}_\infty$
maps the labeled Grassmann-Bruhat order on ${\mathcal S}_\infty$
isomorphically onto its image in the labeled Lagrangian order, preserving
edge labels.
The Pieri-type formula for $\Son/B$ has two formulations (Theorems A and D),
which we relate here. 
We say that a chain in $[e,\zeta]_\prec$ is {\it peakless} if in its
sequence $\beta_1,\ldots,\beta_m$ of labels, we do not have 
$\beta_{i-1}<\beta_i>\beta_{i+1}$, for any $i=2,\ldots,m-1$.

\begin{lem}\label{lem:minimal-peak}
Let $\zeta\in{\mathcal B}_\infty$ be a minimal cycle.
Then there is a unique peakless chain in the labeled interval 
$[e,\zeta]_\prec$.
If $\delta(\zeta)=0$, then the minimal label corresponds to the cover whose
reflection is of the form $t_a$.
\end{lem}

\noindent{\bf Proof. }
If $\delta(\zeta)=1$ this is an immediate consequence of
Lemma~\ref{lem:unique-vee}. 
Suppose that $\delta(\zeta)=0$.
Replacing $\zeta$ by a shape equivalent permutation if necessary, 
we may assume that $\supp(\zeta)=[n]$ and $n>1$.
Replacing $\zeta$ by $\zeta^{-1}$ if necessary, we may assume that
$a:=\zeta^{-1}(n) < b:=\zeta(n)$.
As in the proof of Lemma~\ref{lem:unique-vee},
$(b,n)\zeta\pcdot\zeta$ in the Grassmann-Bruhat order on 
${\mathcal S}_{\pm\infty}$.
By Remark~\ref{rem:covers},
either $t_{b\,n}\zeta\precdot \zeta$ or else we have 
both $t_b\zeta\precdot \zeta$ and $t_n\zeta \precdot \zeta$.
The second case implies $s(\zeta)>1$, contradicting the minimality of
$\zeta$. 
Thus $\eta:=t_{b\,n}\zeta\precdot \zeta$.

Then $\eta$ is a minimal cycle with $\delta(\eta)=0$ and 
$\supp(\eta)=[n{-}1]$.
Appending the cover $\eta\stackrel{n}{\longrightarrow}\zeta$ to the unique
peakless chain in $[e,\eta]_\prec$ gives a peakless chain in
$[e,\zeta]_\prec$. 
Moreover, $\eta$ is the unique permutation with $\eta\precdot\zeta$ and
$\eta\stackrel{n}{\longrightarrow}\zeta$, showing the uniqueness of this
chain.
\QED

For any $\zeta\in{\mathcal B}_\infty$, let $\Pi(\zeta)$ be the number of
peakless chains in $[e,\zeta]_\prec$.

\begin{lem}\label{lem:peak-prod}
If $\eta,\zeta\in {\mathcal B}_\infty$ are disjoint, then
$$
\Pi(\eta\cdot\zeta)\ =\ 2\, \Pi(\eta)\cdot\Pi(\zeta).
$$
\end{lem}

\noindent{\bf Proof. }
For $\xi\in{\mathcal B}_\infty$, let $W(\xi)$ be the multiset of words
formed from labels of maximal chains in $[e,\xi]_\prec$.
The alphabet of these words is a subset of $\supp(\xi)$.
Thus $W(\eta)$ and $W(\zeta)$ have disjoint alphabets.
Note that $W(\eta\cdot\zeta)$ consists of all pairs of words in 
$W(\eta)\times W(\zeta)$.
The lemma follows from Lemma~\ref{lem:combin}, a combinatorial
result concerning peakless words and shuffles.

For a set $A$ of words in an ordered alphabet ${\mathcal A}$, let peak$(A)$
be the subset of peakless words from $A$.  
Suppose that $A'$ is another set of words with a different alphabet
${\mathcal A}'$ and fix some total order on the disjoint union 
${\mathcal A}\coprod {\mathcal A}'$ which extends the given orders on each
of ${\mathcal A},{\mathcal A}'$.
Let sh$(A,A')$ be all shuffles of pairs of words in $A\times A'$.

\begin{lem}\label{lem:combin}
The natural restriction map sh$(A,A') \twoheadrightarrow A\times A'$
induces a 2 to 1 map
$$
\mbox{\rm peak}(\mbox{\rm sh}(A,A'))\ \relbar\joinrel\twoheadrightarrow\ 
\mbox{\rm peak}(A)\times\mbox{\rm peak}(A').
$$
\end{lem}

\noindent{\bf Proof. }
It is clear that the restriction map takes a peakless word in 
sh$(A,A')$ to a pair of peakless words in $A\times A'$.
Given a pair of peakless words $(\omega,\omega')\in A\times A'$,
there are exactly two shuffles of $\omega,\omega'$ which are peakless:
Suppose the minimal letter $a$ in $\omega$ is greater than the minimal
letter in $\omega'$.
Then these two shuffles differ only in their subwords consisting of $a$ and 
$u'$, where $u'$ is that subword of $\omega'$ consisting of all letters less
than $a$.
Then $u'$ is a segment of $\omega'$, as
$\omega'$ is peakless.
The two subwords of peakless shuffles are $a.u'$ and $u'.a$.
\QED

\begin{lem}\label{lem:peakless-unique}
Let $\zeta\in{\mathcal B}_\infty$ and suppose there is a peakless chain in
$[e,\zeta]_\prec$.
Then $\zeta$ is minimal.
\end{lem}

\noindent{\bf Proof. }
Suppose by way of contradiction that $\zeta\in{\mathcal B}_\infty$ is
irreducible and not minimal, but $\Pi(\zeta)\neq 0$.
We may further assume that among all such permutations, $\zeta$ has minimal
rank, and that $\supp(\zeta)=[n]$.
Let $\beta_1>\cdots>\beta_k<\cdots<\beta_m$ be the labels in a peakless
chain in $[e,\zeta]_\prec$.
Replacing $\zeta$ by $\zeta^{-1}$ if necessary (which merely reverses the
chain), we may assume that $\beta_m=n$ and so $\beta_1\neq n$, by
Lemma~\ref{lem:prec_cover} and Theorem~\ref{thm:induced_order}.
Let $\eta$ be the penultimate member of this chain.
Then $\Pi(\eta)\neq 0$, as the initial segment of this chain gives a
peakless chain in $[e,\eta]_\prec$.
Thus $\eta$ is a minimal permutation, by our assumption on $\zeta$, and so 
$$
|\eta|\ \leq\  \#\supp(\eta)-\delta(\eta) \ \leq\ 
 n-\delta(\zeta)\ <\ \Len(\zeta)\ =\ |\eta|+1,
$$
as $\zeta$ is not minimal and $\eta\precdot\zeta$ so
$\delta(\eta)\geq\delta(\zeta)$. 
Therefore the weak inequalities must be equalities, so that
$\supp(\eta)=[n]$. 
Since $\beta_1>\cdots>\beta_k<\cdots<\beta_{m-1}$ are the labels of a chain
in $[e,\eta]_\prec$ and $\beta_{m-1}<n$, we must have $\beta_1=n$, as 
$\supp(\eta)=[n]$.
But this contradicts our earlier observation about $\beta_1$.
\QED

We relate the two formulations of the Pieri formula in the odd-orthogonal
case. 
For $\zeta\in{\mathcal B}_\infty$, define 
$$
\theta(\zeta)\ =\ \left\{ \begin{array}{lcl}
2^{\#\left\{\mbox{\scriptsize irreducible factors of $\zeta$}\right\} -1}&\ &
\mbox{\rm if $\zeta$ is minimal}\\
0&& \mbox{\rm otherwise}\end{array}\right.
$$

\begin{cor}\label{cor:theta-peak}
For $\zeta\in{\mathcal B}_\infty$,
$\Pi(\zeta)=\theta(\zeta)$.
\end{cor}

\noindent{\bf Proof. }
This is clear if $\zeta$ is minimal as both $\Pi(\zeta)$ and $\theta(\zeta)$
satisfy the same recursion, by Lemmas~\ref{lem:minimal-peak}
and~\ref{lem:peak-prod}. 
\QED

\subsection{The Lagrangian r\'eseau}\label{sec:lagr-re}
The enumerative significance of the structure constants $c^w_{u\, v}$ is
best expressed in terms of 
maximal chains in certain directed multigraphs associated to intervals in
the Bruhat order, which we call labeled r\'eseaux.
A cover $\eta\precdot\zeta$ in the Lagrangian order corresponds to a
reflection $\eta\zeta^{-1}$, which is either of the form $t_{a\,b}$
or of the form $t_a$.
The {\em Lagrangian r\'eseau} on ${\mathcal B}_\infty$ is the labeled 
directed multigraph where a cover $\eta\precdot\zeta$ in the Lagrangian
order with $\eta\zeta^{-1}=t_{a}$ is given a single edge
$\eta\stackrel{a}{\longrightarrow}\zeta$ and a cover with 
$\eta\zeta^{-1}=t_{a\,b}$ is given two edges
$\eta\stackrel{\ol{a}}{\longrightarrow}\zeta$ and 
$\eta\stackrel{b}{\longrightarrow}\zeta$.
We obtain the labeled Lagrangian order from this r\'eseau by erasing
those edges whose negative labels.

In the Grassmann-Bruhat order on ${\mathcal S}_\infty$, there are two
conventions for labeling a cover $\eta \pcdot \zeta$:
This cover gives a transposition $(\alpha,\beta):=\eta\zeta^{-1}$
with $\alpha<\beta$, and we may choose either $\alpha$ or $\beta$.
For want of a better term, we call the consistent choice of $\alpha$ the
lower convention, and the consistent choice of $\beta$ the upper convention.
We make use of the following fact.

\begin{prop}\label{prop:unique-increasing}
Let $\eta\in{\mathcal S}_\infty$.
If there is a chain in $[e,\eta]_\prcs$ with decreasing labels in the lower
convention, then there is a chain in $[e,\eta]_\prcs$
with decreasing labels in the upper convention, and these chains are
unique.
The same is true for chains with increasing labels, and in either case
$\eta$ is minimal.
\end{prop}

A chain with increasing labels is an {\it increasing chain} and one
with decreasing labels is a {\it decreasing chain}.

\begin{lem}\label{lem:inc-exists}
Let $\zeta\in{\mathcal B}_\infty$ be a minimal cycle.
Then the r\'eseau $[e,\zeta]_\prec$ has an increasing chain.
If $\delta(\zeta)=1$, then there are at least 2 increasing chains.
\end{lem}

\noindent{\bf Proof. }
Consider the peakless chain in the labeled order $[e,\zeta]_\prec$:
\begin{equation}\label{eq:pkls_chain}
  e\ \stackrel{\alpha_1}{\longrightarrow}\ \zeta_1\ 
  \stackrel{\alpha_2}{\longrightarrow}\ \cdots\ 
  \stackrel{\alpha_m}{\relbar\joinrel\longrightarrow}\ \zeta_m\ =\ \zeta
\end{equation}
Let $\alpha_k$ be the minimal label in this chain.
Then 
$\zeta_{k-1}=\iota(\eta_{k-1})$ for some $\eta_{k-1}\in{\mathcal S}_\infty$.
To see this, if $\delta(\zeta)=0$, then by Lemma~\ref{lem:minimal-peak},
the label $\alpha_k$ corresponds to the only cover whose 
reflection is not in $\iota({\mathcal S}_\infty)$, and so  
$\delta(\zeta_{k-1})=1$.

The pullback of the initial segment of this chain to 
$[e,\eta_{k-1}]_\prcs$
gives a decreasing chain (with labels $\alpha_1,\ldots,\alpha_{k-1}$)
in the upper convention.
Consider the unique decreasing chain
$$
  e\ \stackrel{\beta_1}{\longrightarrow}\ \eta_1\ 
  \stackrel{\beta_2}{\longrightarrow}\ \cdots\ 
  \stackrel{\beta_{k-1}}{\relbar\joinrel\relbar\joinrel\longrightarrow}\ 
  \eta_{k-1}
$$
in the lower convention.
Then
$$
  e\ \stackrel{\ol{\beta_1}}{\relbar\joinrel\longrightarrow}\ \iota(\eta_1)\ 
  \stackrel{\ol{\beta_2}}{\relbar\joinrel\longrightarrow}\ \cdots\ 
  \stackrel{\ol{\beta_{k-1}}}%
     {\relbar\joinrel\relbar\joinrel\longrightarrow}\ 
  \iota(\eta_{k-1})\ =\ \zeta_{k-1}
$$
is an increasing chain in the r\'eseau $[e,\zeta_{k-1}]_\prec$.
Concatenating the end of the peakless chain~(\ref{eq:pkls_chain})
onto this gives an increasing chain
$$
  e\ \stackrel{\ol{\beta_1}}{\relbar\joinrel\longrightarrow}\ \iota(\eta_1)\ 
  \stackrel{\ol{\beta_2}}{\relbar\joinrel\longrightarrow}\ \cdots\ 
  \stackrel{\ol{\beta_{k-1}}}%
     {\relbar\joinrel\relbar\joinrel\longrightarrow}\ 
  \iota(\eta_{k-1})
  \stackrel{\alpha_k}{\longrightarrow}\ \cdots\ 
  \stackrel{\alpha_m}{\relbar\joinrel\longrightarrow}\  \zeta
$$
in the r\'eseau $e,\zeta]_\prec$.

Suppose $\delta(\zeta)=1$ and consider the middle portion of this increasing
chain:
$$
\iota(\eta_{k-2})\stackrel{\ol{\beta_{k-1}}}%
     {\relbar\joinrel\relbar\joinrel\longrightarrow}\ 
\iota(\eta_{k-1})=\zeta_{k-1}\ 
\stackrel{\alpha_k}{\longrightarrow}\ 
\zeta_k.
$$
Let $\ol{b}$ be the label of the other edge between $\zeta_{k-1}$ and
$\zeta_k$. 
Then we claim that $\ol{\beta_{k-1}}<\ol{b}$ so that replacing 
$\zeta_{k-1} \ \stackrel{\alpha_k}{\longrightarrow}\ \zeta_k$
by 
$\zeta_{k-1} \ \stackrel{\ol{b}}{\longrightarrow}\ \zeta_k$
gives a second increasing chain in the r\'eseau.

To see this, first note that $\beta_{k-1}=b$ is impossible as these are
consecutive covers in the Lagrangian order 
(see relation (iv) of Equation~(\ref{eq:relations})). 
Define $\eta$ by $\iota(\eta)=\zeta$ and pull this chain back to 
$[e,\eta]_\prcs$.
It is the unique peakless chain in $[e,\eta]_\prcs$ and $\alpha_k$ is the
minimal label.
By Lemma~\ref{lem:unique-vee}(1), $b=\min(\supp(\eta))$ and so 
$\beta_{k-1}\geq b$.
\QED

\begin{lem}
Let $\zeta\in{\mathcal B}_\infty$ and suppose there is an increasing chain
in the r\'eseau $[e,\zeta]_\prec$.
Then $\zeta$ is minimal.
If $\zeta$ is a minimal cycle, then there are precisely $2^{\delta(\zeta)}$
such chains. 
\end{lem}

\noindent{\bf Proof. }
Let 
\begin{equation}\label{eq:inc_chain}
  e\ \stackrel{\beta_1}{\longrightarrow}\ \zeta_1\ 
  \stackrel{\beta_2}{\longrightarrow}\ \cdots\ 
  \stackrel{\beta_m}{\relbar\joinrel\longrightarrow}\ \zeta_m\ =\ \zeta
\end{equation}
be an increasing chain in $[e,\zeta]_\prec$.
Suppose that $\beta_{k-1}<0<\beta_k$.
Then for $i<k$, $\delta(\zeta_i)=1$.
Define $\eta_i\in{\mathcal S}_\infty$ by $\iota(\eta_i)=\zeta_i$ for $i<k$.
Then
$$
  e\ \stackrel{\ol{\beta_1}}{\longrightarrow}\ \eta_1\ 
  \stackrel{\ol{\beta_2}}{\relbar\joinrel\longrightarrow}\ \cdots\ 
  \stackrel{\ol{\beta_{k-1}}}%
  {\relbar\joinrel\relbar\joinrel\longrightarrow}\ \eta_{k-1}
$$
is a decreasing chain in $[e,\eta]_\prcs$, with the lower labeling
convention.
Let 
$$
   e\ \stackrel{\alpha_1}{\longrightarrow}\ \xi_1\ 
  \stackrel{\alpha_2}{\longrightarrow}\ \cdots\ 
  \stackrel{\alpha_{k-1}}{\relbar\joinrel\relbar\joinrel\longrightarrow}\ 
  \xi_{k-1}\ =\ \eta_{k-1}
$$
be the unique decreasing chain in the upper labeling convention.
Concatenating  the image of this chain in $[e,\zeta]_\prec$ with the end of
the chain~(\ref{eq:inc_chain}) gives a peakless chain
\begin{equation}\label{eq:peakless}
   e\ \stackrel{\alpha_1}{\longrightarrow}\ \iota(\xi_1)\ 
  \stackrel{\alpha_2}{\longrightarrow}\  \cdots\ 
  \stackrel{\alpha_{k-1}}{\relbar\joinrel\relbar\joinrel\longrightarrow}\ 
  \iota(\xi_{k-1})=\zeta_{k-1}
  \stackrel{\beta_k}{\longrightarrow}\ \cdots\ 
  \stackrel{\beta_m}{\relbar\joinrel\longrightarrow}\  \zeta
\end{equation}
in the interval $[e,\zeta]_\prec$ in the Lagrangian order.
By Lemma~\ref{lem:peakless-unique}, $\zeta$ is necessarily minimal.
\smallskip

Suppose now that $\zeta$ is a minimal cycle, then the r\'eseau
$[e,\zeta]_\prec$ has an increasing chain and the peakless
chain~(\ref{eq:peakless}) is unique. 
Consider another increasing chain
\begin{equation}\label{eq:new_inc_chain}
  e\ \stackrel{\beta'_1}{\longrightarrow}\ \zeta'_1\ 
  \stackrel{\beta'_2}{\longrightarrow}\ \cdots\ 
  \stackrel{\beta'_m}{\relbar\joinrel\longrightarrow}\ \zeta'_m\ =\ \zeta
\end{equation}
and form $\eta'_i$, $\xi'_i$, $\alpha'_i$,and $k'$ as for the original
chain~\ref{eq:inc_chain}). 
If $k=k'$, then the chains~(\ref{eq:inc_chain}),~(\ref{eq:new_inc_chain})
coincide:
The final segments agree, by the uniqueness of~(\ref{eq:peakless}), as do
their initial segments, by Proposition~\ref{prop:unique-increasing}.

If $\delta(\zeta)=0$, then the minimal label in the peakless
chain~(\ref{eq:peakless}) (either $\alpha_{k-1}$ or $\beta_k$) corresponds
to the cover whose reflection has the form $t_a$.
As $\delta(\zeta_{k-1})=1$, this must be $\beta_k$ and so $k=k'$
and the chain~(\ref{eq:inc_chain}) is the unique increasing chain in the
r\'eseau $[e,\zeta]_\prec$.

Suppose now that $\delta(\zeta)=1$ and $k<k'$.
Since $\alpha_i=\alpha'_i$ for $i<k$, $\beta_i=\beta'_i$ for $i\geq k'$, and
$\alpha_1>\cdots>\alpha_{k-1}$ and $\beta_k<\cdots<\beta_m$, we must have
$k'=k+1$. 
But then $\xi_i=\xi'_i$ for $i<k$ and also $\zeta_i=\zeta'_i$ for $i\geq k$,
and so the two chains~(\ref{eq:inc_chain}) and~(\ref{eq:new_inc_chain})
agree except for the label of the cover 
$\zeta_{k-1}\precdot \zeta_k$.
Thus there are at most 2 increasing chains in the r\'eseau
$[e,\zeta]_\prec$ and their underlying permutations coincide.
\QED

\begin{ex}
Suppose $\eta=(1,2,5,3,4)$ is a permutation in ${\mathcal S}_5$.
Consider $\zeta=\Span{1,2,5,3,4}= \eta\ol{\eta}$, a permutation in 
${\mathcal B}_5$.
Figure~\ref{fig:up-down} shows the r\'eseau $[e,\zeta]_\prec$.
\begin{figure}[htb]
$$
\epsfxsize=2.55in \epsfbox{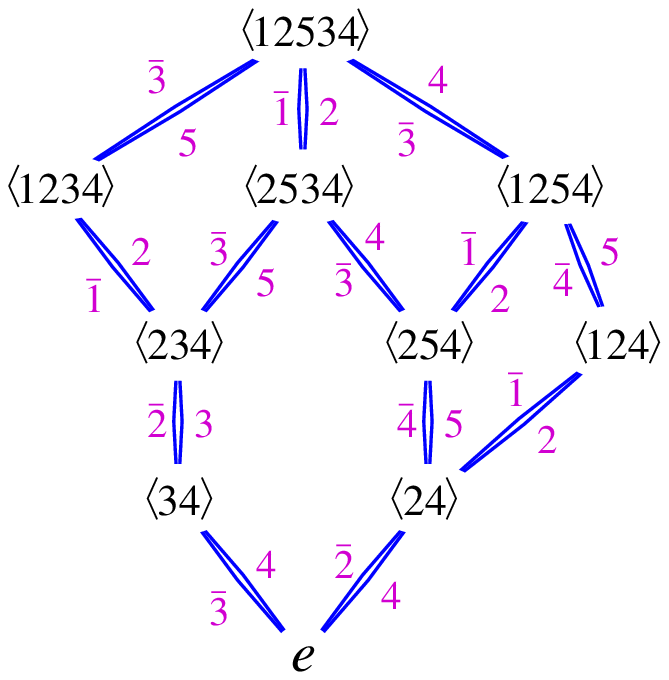}
$$
\caption{The interval $[e,\Span{1,2,5,3,4}]_\prec$. \label{fig:up-down}}
\end{figure}
In this r\'eseau, there are two increasing chains
\begin{eqnarray*}
e\ \stackrel{\ol{3}}{\longrightarrow}\ 
\Span{3,4}\ \stackrel{\ol{2}}{\longrightarrow}\ 
\Span{2,3,4}\ \stackrel{\ol{1}}{\longrightarrow}\ 
\Span{1,2,3,4}\ \stackrel{5}{\longrightarrow}\ 
\Span{1,2,5,3,4}\\
e\ \stackrel{\ol{3}}{\longrightarrow}\ 
\Span{3,4}\ \stackrel{\ol{2}}{\longrightarrow}\ 
\Span{2,3,4}\ \stackrel{2}{\longrightarrow}\ 
\Span{1,2,3,4}\ \stackrel{5}{\longrightarrow}\ 
\Span{1,2,5,3,4}
\end{eqnarray*}
which correspond to the unique peakless chain in $[e,\eta]_\prcs$
with $\beta_1>\beta_2>\beta_3<\beta_4$ 
$$
e\ \stackrel{4}{\longrightarrow}\ 
(3,4)\ \stackrel{3}{\longrightarrow}\ 
(2,3,4)\ \stackrel{2}{\longrightarrow}\ 
(1,2,3,4)\ \stackrel{5}{\longrightarrow}\ 
(1,2,5,3,4)\,.
$$
\end{ex}

Let $I(\zeta)$ count the increasing chains in the r\'eseau
$[e,\zeta]_\prec$.

\begin{lem}
If\/ $\eta,\zeta\in{\mathcal B}_\infty$ are disjoint, then 
$I(\eta\cdot\zeta)=I(\eta)\cdot I(\zeta)$.
\end{lem}

\noindent{\bf Proof. }
As with Lemma~\ref{lem:peak-prod}, this is a consequence of the analogous
bijection concerning increasing words among shuffles of words with disjoint
alphabets.
\QED

For $\zeta\in{\mathcal B}_\infty$, define 
$$
\chi(\zeta)\ =\ \left\{ \begin{array}{lcl}
2^{\#\left\{\mbox{\scriptsize irreducible
           factors $\eta$ of $\zeta$ with $\delta(\eta)=1$}\right\} }&\ &
\mbox{\rm if $\zeta$ is minimal}\\
0&& \mbox{\rm otherwise}\end{array}\right.
$$

Let $D(\zeta)$ enumerate the decreasing chains in the r\'eseau
$[e,\zeta]_\prec$.
By Theorem~\ref{thm:prpoerties}(6),
an increasing chain in $\zeta$ becomes a decreasing chain in $\zeta^{-1}$.
The following result is now immediate.

\begin{cor}\label{cor:chi-descent}
For $\zeta\in{\mathcal B}_\infty$,
$\chi(\zeta)=I(\zeta)=D(\zeta)$.
\end{cor}

\section{The Pieri-type formula}\label{sec:pieri}
Recall that for $\zeta\in{\mathcal B}_n$,
$\delta(\zeta)=1$  if $a>0$ implies $\zeta(a)>0$ and $\delta(\zeta)=0$ 
otherwise. 

\begin{lem}\label{lem:dimension}
Suppose $\zeta\in{\mathcal B}_n$ with $\delta(\zeta)=1$.
Then there is a quadratic form $f$ on 
$V$ with $f|_K\equiv 0$ for every $K\in{\mathcal Y}_\zeta$.
\end{lem}

\noindent{\bf Proof. }
We may suppose that $\supp(\zeta)=[n]$; otherwise pull back the quadratic
form along the map $V\hookrightarrow V\oplus (H_2)^{\oplus t}$ 
($t=n-\#\supp(\zeta)$)
of Section~\ref{sec:identity}.
Let $\Edot,\Epdot$ be opposite isotropic flags.
Since $\Edot, \Epdot$ are opposite flags, we have
$V= E'_{\ol{1}} \oplus E_{\ol{1}}$.
Let $\beta$ be the alternating form on $V$ and define a symmetric bilinear
form $q$ on $V$ by
$$
q(x,y)\ :=\ \beta(x^+,y^-) + \beta(y^+,x^-),
$$
where $v^+, v^-$ are the  projections of $v\in V$ to the summands
$E'_{\ol{1}}$ and $E_{\ol{1}}$.
This form is non-degenerate.
Let $f$ be the associated quadratic form.

Let $u\in {\mathcal B}_n$ with $u\leq_0\zeta u$ and $\zeta u$ a Grassmannian
permutation, as in Remark~\ref{rem:X}.
Then there is a $k\in [n]$ with 
$i\leq k \Rightarrow \zeta u(i)\leq \ol{1}$ and 
$i> k \Rightarrow \zeta u(i)\geq 1$.
Since $\delta(\zeta)=1$, $u(i)\leq \ol{1}$ for $i\leq k$ and 
$u(i)\geq 1$ for $i> k$.
Let $\Fdot\in Y_u\Edot\bigcap Y_{\omega_0\zeta u}\Epdot$.\
Then
$$
\dim F_{\ol{1}}\bigcap E_{\ol{1}}\ \geq\ 
\#\{l\geq 1\mid u(l)\leq \ol{1}\,\} \ =\ k
$$
and 
$$
\dim F_{\ol{1}}\bigcap E'_{\ol{1}}\ \geq\ 
\#\{l\geq 1\mid \zeta u(l)\geq 1\} \ =\ n-k.
$$
Thus 
$F_{\ol{1}}=F_{\ol{1}}\bigcap E_{\ol{1}}\oplus 
F_{\ol{1}}\bigcap E'_{\ol{1}}$.
Since $F_{\ol{1}}$ is isotropic for the alternating form $\beta$, this
decomposition shows it is isotropic for the symmetric form $q$, proving the
lemma. 
\QED

\begin{cor}\label{cor:quad_forms}
Suppose $\zeta=\zeta_1\cdots\zeta_s$ is the factorization of $\zeta$ into
irreducible permutations.
Suppose exactly $r$ of the $\zeta_i$ have $\delta(\zeta_i)=1$.
Then there exist $r$ linearly independent quadratic forms $f_1,\ldots,f_r$ on
$V$ with $f_i|_K\equiv 0$ for every $K\in{\mathcal Y}_\zeta$.
\end{cor}

\noindent{\bf Proof. }
By Theorem~\ref{thm:geom_facts}(1), there exist symplectic spaces
$V_1,\ldots,V_s$ with $V=V_1\oplus\cdots\oplus V_s$
such that $K\in {\mathcal Y}_\zeta$ may be written as
$K=K_1\oplus\cdots\oplus K_s$, where, for $i=1,\ldots,s$, 
 $K_i=K\bigcap V_i$ and
$K_i\in{\mathcal Y}_{\zeta'_i}$, with $\zeta'_i\in{\mathcal B}_{n_i}$
shape equivalent to $\zeta_i$.
(Here, $\dim V_i=2n_i$.)
We may assume that the $\zeta_i$ are ordered so that $\delta(\zeta_i)=1$ for
$i=1,\ldots,r$ and $\delta(\zeta_i)=0$ for $i>r$.
By Lemma~\ref{lem:dimension}, for each $i=1,\ldots,r$, there is a quadratic
form $f_i$ on $V_i$ vanishing on $K_i\in{\mathcal Y}_{\zeta'_i}$.
The pullback of $f_i$ along the orthogonal projection 
$V\twoheadrightarrow V_i$ is a quadratic form on $V$, also denoted $f_i$.
These are the desired forms.
\QED

The special Schubert variety $\Upsilon_{(m)}$ of $\Spn/P_0$ consists of
those Lagrangian subspaces which meet a fixed $(n{+}1{-}m)$-dimensional
isotropic subspace $M$ of $V$.
In what follows, write $\Upsilon_M$ for this Schubert variety, and assume
$M$ is in general position in $V$.

For $\zeta\in{\mathcal B}_n$, define 
$$
\chi(\zeta)\ =\ \left\{ \begin{array}{lcl}
2^{\#\left\{\mbox{\scriptsize irreducible
             factors $\eta$ of $\zeta$ with $\delta(\eta)=1$}\right\} }&\ &
\mbox{\rm if $\zeta$ is minimal}\\
0&& \mbox{\rm otherwise}\end{array}\right.
$$
For $\zeta\in{\mathcal B}_n$, define ${\mathcal W}_\zeta$ to be the zero
locus of the forms $f_1,\ldots,f_r$ of Corollary~\ref{cor:quad_forms}, a
complete intersection.
Then ${\mathcal W}_\zeta$ is the cone over a subvariety of 
${\mathbb P}V$ of degree $2^r$.
If $\zeta$ is minimal, then this degree is $\chi(\zeta)$.

\begin{cor}
Let $\zeta\in{\mathcal B}_n$ with $\Len(\zeta)=m$ and $\supp(\zeta)=[n]$.
Let $M$ be a general isotropic $(n{+}1{-}m)$-dimensional subspace of $V$.
Then $M\bigcap {\mathcal W}_\zeta=\{0\}$ unless 
$\zeta$ is minimal, and in that case, $M\bigcap {\mathcal W}_\zeta$
is $\chi(\zeta)$ reduced lines.
\end{cor}

\noindent{\bf Proof. }
Since $\Spn$ acts transitively on ${\mathbb P}V$, Kleiman's theorem on the
transversality of a general translate~\cite{Kleiman} will imply
the corollary if we show
\begin{equation}\label{eq:dim}
\dim M + \dim {\mathcal W}_\zeta\ \leq \ n+1,
\end{equation}
with equality only if $\zeta$ is minimal.

By Corollary~\ref{cor:quad_forms}, 
$\dim {\mathcal W}_\zeta=2n-r$, where $r$ counts the irreducible factors
$\eta$ of $\zeta$ with $\delta(\eta)=1$.
Since $n=\#\supp(\zeta)$, Corollary~\ref{cor:minimal_perm} implies that 
$\dim M=n+1-\Len(\zeta)\leq 1+\sum_i \delta(\zeta_i)=1+r$,
with equality only if each $\zeta_i$ is a minimal cycle,
establishing~(\ref{eq:dim}).
\QED

\begin{thm}\label{thm:singleton}
Suppose $\zeta\in{\mathcal B}_n$ is minimal
and $\supp(\zeta)=[n]$.
Then a general line $\langle v\rangle$ in ${\mathcal W}_\zeta$ determines a
unique $K\in{\mathcal Y}_\zeta$ with $v\in K$.
\end{thm}

We deduce the Pieri-type formula from Theorem~\ref{thm:singleton}.
First, define $\theta(\zeta)=0$ if $\zeta$ is not
minimal, and for $\zeta$ minimal, set
$$
\theta(\zeta)\ :=\ 2^{\#\{\mbox{\scriptsize 
                  irreducible factors of $\zeta$}\} -1}.
$$
Recall that $b^\zeta_m$ was the structure constant corresponding to  
$c^\zeta_m$ for $\Son/B$.
\medskip

\noindent{\bf Theorem D. }(Pieri-type Formula)
{\it 
Let $\zeta\in{\mathcal B}_\infty$ with $\Len(\zeta)=m$.
Then $c^\zeta_{m}=\chi(\zeta)$ and $b^\zeta_m=\theta(\zeta)$.}\medskip

By Corollaries~\ref{cor:theta-peak} and~\ref{cor:chi-descent}, this implies 
the chain-theoretic version of the Pieri-type formula (Theorem A).
\medskip

\noindent{\bf Proof. }
Let $s(\zeta)$ be the number of sign changes in a permutation 
$\zeta\in{\mathcal B}_\infty$.
Since $s(v_m)=1$ and $s(\zeta)=s(\zeta u)-s(u)$ if $u\leq_0 \zeta u$, 
Equation~(\ref{eq:coeff_compare}) implies that
$b^\zeta_m=2^{s(\zeta)-1}c^\zeta_m$.
Since, for a minimal cycle $\zeta$, $s(\zeta)+\delta(\zeta)=1$,
we have
$$
\theta(\zeta)\ =\ 2^{s(\zeta)-1} \chi(\zeta).
$$
Thus it suffices to show $c^\zeta_{m}=\chi(\zeta)$.

By Theorem~B(2), replacing $\zeta$ by a shape equivalent permutation
if necessary, we may assume that $\supp(\zeta)=[n]$.
Let $M$ be a general isotropic $(n+1-m)$-dimensional subspace of $V$.
By the projection formula, 
$c^\zeta_m=\deg(\Upsilon_M\bigcap{\mathcal Y}_\zeta)$.
Since $K\in \Upsilon_M\bigcap{\mathcal Y}_\zeta$ implies that $K$ meets 
$M\bigcap {\mathcal W}_\zeta$ non-trivially, we see that 
$$
c^\zeta_m\ =\ d\cdot \deg(M\bigcap {\mathcal W}_\zeta),
$$
where $d$ counts the $K\in {\mathcal Y}_\zeta$ which contain a general line
of ${\mathcal W}_\zeta$, and the degree is taken in ${\mathbb P}V$.
By Theorem~\ref{thm:singleton}, $d=1$, which completes the proof.
\QED

\noindent{\bf Reduction of Theorem~\ref{thm:singleton} 
to the case of $\zeta$ a minimal cycle.}

Let $\zeta=\zeta_1\cdots\zeta_s$ be the irreducible factorization 
of $\zeta$.
In the notation of the proof of Corollary~\ref{cor:quad_forms}, 
a general $0\neq v\in {\mathcal W}_\zeta$ has the form
$v=v_1\oplus \cdots \oplus v_s$, where 
$0\neq v_i\in{\mathcal W}_{\zeta'_i}$ for $i=1,\ldots,s$.
Moreover $v\in K\in{\mathcal Y}_\zeta$ if and only if 
$v_i\in K_i\in{\mathcal Y}_{\zeta'_i}$.
Thus is suffices to prove Theorem~\ref{thm:singleton} for 
$\zeta$ a minimal cycle.
We do this in the following sections.

\subsection{Case $\delta(\zeta)=1$}

\begin{thm}
Let $\zeta\in{\mathcal B}_n$ with 
$\supp(\zeta)=[n]$ and $\Len(\zeta)=n-1$ so that $\zeta$ is a minimal cycle
with  $\delta(\zeta)=1$.
Then, for a general $0\neq v\in{\mathcal W}_\zeta$, there is a unique
$K\in{\mathcal Y}_\zeta$ with $v\in K$.
\end{thm}

\noindent{\bf Proof. }
Define $\eta\in{\mathcal S}_\infty$ by $\iota(\eta)=\zeta$.
Set $k:=\#\{a\mid a<\eta(a)\}$.
Recall the notation of Section~\ref{sec:hidden}:
Let $L,L^\perp$ be complementary Lagrangian subspaces of $V$, which are
identified as linear duals.
Define the map $\Phi_k:{\bf G}_k(L)\rightarrow \Spn/B$ by
$H\mapsto(H+H^\perp)$, where $H^\perp\subset L^\perp$ is the annihilator
of $H$. 
Define $\pi_k:{\mathbb F}\ell(L)\rightarrow {\bf G}_k(L)$ by 
$\Edot\mapsto E_k$.
By Corollary~\ref{cor:hidden-identity},  
$\Phi_k:{\mathcal X}_\eta \stackrel{\sim}{\longrightarrow}
{\mathcal Y}_\zeta$
where ${\mathcal X}_\eta :=\pi_k(X_u\Edot\bigcap X_{(\eta u)^\vee}\Epdot)$.

Schubert varieties $\Omega_\varrho$ of the Grassmannian ${\bf G}_k(L)$
are indexed by ordinary partitions ({\it weakly} decreasing sequences)
$\varrho: n-k\geq \varrho_1\geq\cdots\geq
      \varrho_k\geq 0$~\cite{Fulton_tableaux}.
We show
\begin{equation}\label{eq:pullback}
 \Phi_k^{-1}\left(\{K\in\Spn/B\mid v\in K\}\right)
 \ =\ 
 \Omega_{(n-k,1^{k-1})},
\end{equation}
where 
$(n-k,1^{k-1})$ is the hook-shaped partition with first row $n-k$ and first
column $k$.
It follows from the projection formula that 
$$
\deg\left({\mathcal Y}_\zeta\bigcap \{K\mid v\in K\}\right)
\ =\ \deg\left({\mathcal X}_\eta \bigcap \Omega_{(n-k,1^{k-1})}\right)
$$
which is 
$\deg({\mathfrak S}_w \cdot {\mathfrak S}_{\omega_0\eta w}\cdot
\pi_k^* S_{(n-k,1^{k-1})})$,
the product in $H^*{\mathbb F}\ell L$.
By~\cite[Theorem 8]{Sottile96}, this counts 
the chains in the $k$-Bruhat order
$$
w \stackrel{\beta_1}{\longrightarrow} w_1 
\stackrel{\beta_2}{\longrightarrow} \cdots 
\stackrel{\beta_{n-1}}{\relbar\joinrel\longrightarrow} \eta w
$$
with $\beta_1 > \beta_2 >\cdots>\beta_k <\beta_{k+1}<\cdots<\beta_{n-1}$.
The conclusion follows by Lemma~\ref{lem:unique-vee}.

To show (\ref{eq:pullback}), suppose $H\in{\bf G}_k(L)$ and with 
$v\in H+H^\perp$.
Let $v^-$ be the projection of $v$ into $L$ and $v^+$ its projection to
$L^\perp$.
Then $v=v^-\oplus v^+$ and $v^-\in H$ and $v^+\in H^\perp$ so that 
$H\subset (v^+)^\perp$.
Thus 
$\Phi_k^{-1}(\{K\mid v\in K\})=
\{H\mid v^-\in H \mbox{ and }H\subset(v^+)^\perp\}$, which is just the
Schubert variety $\Omega_{(n{-}k,1^{k-1})}$.
\QED

\subsection{Case $\delta(\zeta)=0$}

\begin{thm}\label{thm:last}
Let $\zeta\in{\mathcal B}_n$ with $\delta(\zeta)=0$ and $\Len(\zeta)=n$.
Then $c^\zeta_n=1$.
\end{thm}

This completes the proof of Theorem~\ref{thm:singleton}, as 
$c^\zeta_n=1$ counts the $K\in{\mathcal Y}_\zeta$ which meet a generic line
in $V$.
We use some geometric constructions to reduce the computation of $c^\zeta_n$
to the cohomology of the classical flag manifold.
Let $L, L^\perp$ be complementary Lagrangian subspaces of $V$ with $L^\perp$
identified with the linear dual of $L$ as in Section~\ref{sec:hidden}.
For $1\leq k\leq n$, set
$$
{\mathbb F}_k\ :=\ \{ F_{k-1}\subset F_k\subset L\mid \dim F_i=i\},
$$
a variety of partial flags in $L$.

Let $\wp:{\mathbb F}\ell(L)\rightarrow {\mathbb F}_k$ be the projection.
Then the projections $\pi_{k-1},\pi_k$ of ${\mathbb F}\ell(L)$ to 
${\bf G}_{k-1}(L)$, ${\bf G}_k(L)$ factor through $\wp$.
Let $\pi_{k-1},\pi_k$ also denote the projections of ${\mathbb F}_k$ to
${\bf G}_{k-1}(L)$, ${\bf G}_k(L)$.
Let $\mbox{\it Lag}(V)$ denote $\Spn/P_0$, the Grassmannian of Lagrangian
subspaces of $V$.
Consider the incidence variety $\Gamma$:
$$
\begin{picture}(250,66)(0,-1)
\put(0,3){${\mathbb F}_k$}
\put(32,3){$\mbox{\it Lag}(V)$}
\put(22.5,53){$\Gamma:=\{ (F_{k-1},F_k,K)\mid F_{k-1}\oplus F^\perp_k
\subsetneq K \subsetneq F_k\oplus F_{k-1}^\perp\}$}
\put(4.5,30){$g$}	\put(41,30){$f$}
\put(21,46){\vector(-1,-2){15}} \put(28,46){\vector(1,-2){15}}
\end{picture}
$$
Then $\Gamma$ is a ${\mathbb P}^1$-bundle over ${\mathbb F}_k$ and $f$ is
genericaly 1-1:
The image of $f$ consists of those $K$ with $\dim K\bigcap L\geq k-1$ and
$\dim K\bigcap L^\perp\geq n-k$, which is an intersection of Schbert
varieties. 
Thus a generic $K$ in this intersection determines 
$g(f^{-1}(K))= ( K\bigcap L, (K\bigcap L^\perp)^\perp)$ uniquely.

For $j\in[n]$ and a (not necessarily strict) partition $\varrho$ with 
$n-j\geq \varrho_1\geq \cdots\geq \varrho_j\geq 0$, let 
$\sigma_\varrho\in H^*{\bf G}_j(L)$ be the Schubert class associated to the
partition $\varrho$, as in~\cite{Fulton_tableaux}.
We show:

\begin{lem}\label{lem:push-pull}
$g_*f^* q_n\ =\ \pi_k^*\sigma_{(n-k,1^{k-1})} + 
\pi_{k-1}^*\sigma_{(n-k+1,1^{k-2})}$.
\end{lem}

Let $\zeta\in{\mathcal B}_n$ be minimal with
$\delta(\zeta)=0$, $\Len(\zeta)+n$, and $\supp(\zeta)=[n]$.
We construct a minimal permutation $\eta\in{\mathcal S}_n$ with
$||\eta||=n-1$, a $k\in[n]$, and $w\in {\mathcal S}_n$ with 
$\eta w\vartriangleleft_k w$, and $\eta w\not\vartriangleleft_{k-1} w$.
Let ${\mathcal X}_{\eta^{-1}}$ be $\wp(X_{\eta w}\bigcap X'_{w^\vee})$.
Then 
$[{\mathcal X}_{\eta^{-1}}]=
\wp_*({\mathfrak S}_\eta w\cdot{\mathfrak S}_{w^\vee})$.
We show

\begin{lem}\label{lem:geom-pull-push}
$[{\mathcal Y}_\zeta]\ =\ f_*g^*[{\mathcal X}_{\eta^{-1}}]$.
\end{lem}

Theorem~\ref{thm:last} follows from these Lemmas.
\begin{eqnarray*}
c^\zeta_n\ =\ \deg([{\mathcal Y}_\zeta]\cdot q_n) 
&=& \deg (f_*g^* [{\mathcal X}_{\eta^{-1}}]\cdot q_n)\\
&=& \deg( [{\mathcal X}_{\eta^{-1}}]\cdot g_*f^*  q_n),
\end{eqnarray*}
by the projection formula and Lemma~\ref{lem:geom-pull-push}.
By Lemma~\ref{lem:push-pull},
this is 
$$
\deg \left( {\mathfrak S}_{\eta w}\cdot {\mathfrak S}_{w^\vee}
\cdot(\pi_k^*\sigma_{(n-k,1^{k-1})} + \pi_{k-1}^*\sigma_{(n-k+1,1^{k-2})})
\right).
$$
Since $\eta w\not\vartriangleleft_{k-1} w$, only the first term is
non-zero.
By~\cite[Theorem 8]{Sottile96} and Lemma~\ref{lem:unique-vee}, this degree
is 1.
\QED

\noindent{\bf Proof of Lemma~\ref{lem:push-pull}}
The class $q_n\in H^*\mbox{\it Lag}(L)$ is represented by the Schubert
variety 
$$
\Upsilon_v\ :=\ \{K\in\mbox{\it Lag}(V)\mid v\in K\}\ =\ 
            \{K\mid \beta(v,K)\equiv 0\},
$$
where $0\neq v\in V$ and $\beta$ is the alternating form.
Then $g_*f^* q_n$ is represented by $g(f^{-1} \Upsilon_v)$ which is 
$$
  \{ F_{k-1}\subset F_k\mid \exists K \mbox{ with }
    v\in K\mbox{ and } F_{k-1}\oplus F^\perp_k
    \subsetneq K \subsetneq F_k\oplus F_{k-1}^\perp\}.
$$
Since $V=L\oplus L^\perp$, we may write a general $v$ uniquely as 
$v=w\oplus u$ with $w\in L$ and $u\in L^\perp$
and so $g(f^{-1} \Upsilon_v)$ is a subset of 
\begin{equation}\label{eq:int_Schu_var}
  \{F_{k-1}\subset F_k\mid w\in F_k\} \ \bigcap \
   \{F_{k-1}\subset F_k\mid F_{k-1}\subset u^\perp\}.
\end{equation}
This is an intersection of Schubert varieties (in general position if $v$ is
general) of codimensions $n-k$ and $k-1$, respectively.
These Schubert varieties have classes $\pi_k^*\sigma_{(n-k)}$ and 
$\pi_{k-1}^*\sigma_{(1^{k-1})}$.
Since $\Upsilon_v$ has codimension $n$, $g(f^{-1} \Upsilon_v)$ equals
this intersection if the map 
$g: f^{-1} \Upsilon_v \rightarrow g(f^{-1} \Upsilon_v)$ is finite.
Thus
$$
g_* f^* q_n\ =\ 
d (\pi_k^*\sigma_{(n-k)} \cdot \pi_{k-1}^*\sigma_{(1^{k-1})}),
$$
where $d$ is the degree of the map 
$g: f^{-1} \Upsilon_v \rightarrow g(f^{-1} \Upsilon_v)$
(which is 0 if the map is not finite).

To compute $d$, let $F_{k-1}\subset F_k$ satisfy $w\in F_k$ and 
$F_{k-1}\subset u^\perp$ with $F_k\not\subset u^\perp$
and $w\not\in F_{k-1}$.
Then $F_k\oplus F_{k-1}^\perp \not\subset v^\perp$, and so 
$f(g^{-1}(F_{k-1}, F_k))=F_k\oplus F_{k-1}^\perp\bigcap v^\perp$,
which shows $d=1$.
Lastly, the Pieri-type formula~\cite{Sottile96} in $H^*{\mathbb F}_k$ shows
$$
\pi_k^*\sigma_{(n-k)} \cdot \pi_{k-1}^*\sigma_{(1^{k-1})}
\ =\ \pi_k^*\sigma_{(n-k),1^{k-1})} + 
\pi_{k-1}^*\sigma_{(n-k+1,1^{k-2})}. 
\qquad \QED
$$

\noindent{\bf Proof of Lemma~\ref{lem:geom-pull-push}}
We first make some definitions.
Replacing $\zeta$ by $\zeta^{-1}$ if necessary, we may assume that
$\alpha>0$ is the unique number with $t_{\alpha}\zeta\precdot\zeta$,
by Lemma~\ref{lem:sign-cover}.
Since $\delta(t_{\alpha}\zeta)=1$, we define $\eta\in{\mathcal S}_n$ by
$\iota(\eta)=t_{\alpha}\zeta$. 
Set  
$$
k\ :=\ \#\{i\mid i>\eta(i)\}\ =\ \{i>0\mid \zeta(i)>i\}.
$$
Let $u\in{\mathcal B}_n$ satisfy $u\leq_0\zeta u$ with $\zeta u$ a
Grassmannian permutation ({\it cf.}~Remark~\ref{rem:X}).
Since $\delta(\zeta)=0$, 
$\zeta u(n-k+2)>0>\zeta u(n+1-k)$.
Define $j\geq k$ by $\zeta u(n+1-j)=\ol{\alpha}$.
Then $u(n+1-j)>0$.
Let $w\in{\mathcal S}_n$ be defined by
\begin{equation}\label{def:w}
  w(i)\ =\ \left\{\begin{array}{lcl}
       u(n+1-i)&\quad& i< k       \\
       u(n+1-j)&& i=k             \rule{0pt}{14pt}\\\
       \ol{u(i-k)}&&k<i<n+1+k-j   \rule{0pt}{14pt}\\\
       \ol{u(i+1-k)}&&n+1+k-j\leq i\leq n \rule{0pt}{14pt}\end{array}\right. .
\end{equation}
We claim that $\eta w\vartriangleleft_k w$.
Since $\supp(\eta)=[n]$, we cannot have both $\eta w\vartriangleleft_k w$ and 
$\eta w\vartriangleleft_{k-1} w$.
Since $u<_0\iota(\eta) u$, condition (1) of Proposition~\ref{prop:one} is
satisfied.
For condition (2), as $\zeta u$ is a Grassmannian permutation, 
$$
\eta w(1)>\cdots>\eta w(k-1)\quad\mbox{ and }\quad
\eta w(k+1)>\cdots>\eta w(n),
$$
so we only need show $i<k$ with $\eta w(i)<\eta w(k)=\alpha$ implies
$w(i)<w(k)$.
Let $l=n+1-i>n+1-j$.
Then $\eta w(i)=\iota(\eta)u(l) < \iota(\eta)u(n+1-j)=\alpha$,
which implies $u(l)<u(n+1-j)$ and hence $w(i)<w(k)$, as 
$u<_0\iota(\eta) u$.

\begin{ex}
Let $\zeta=\Perm{1,2,3,5,7,6,\ol{4}}\in{\mathcal B}_7$.
Then $\alpha=4$ and $t_{4}\zeta=\Span{1,2,3,5,7,6,4}\precdot \zeta$
so that $\eta=(1,2,3,5,7,6,4)\in{\mathcal S}_7$.
Here, $k=3$.
If we set $u=\ol{5}\,\ol{3}6\ol{2}\,\ol{1}47$, then 
$\zeta u=\ol{7}\,\ol{5}\,\ol{4}\,\ol{3}\,\ol{2}16$ is a Grassmannian
permutation and $u<_0 \zeta u$.
We see that $j=5\geq 3=k$ and $w=746\,5321$ so that 
$\eta w=614\,7532$ and $\eta w\vartriangleleft_3 w$.
\end{ex}

Lemma~\ref{lem:geom-pull-push} is a consequence of the following construction:

\begin{lem}\label{lem:big-diagram}
Let $\zeta,u,\eta,w$,and $k$ be as above.
Then there exists a commutative diagram
$$
\begin{picture}(250,145)(0,-4)
\put(24,0){${\mathcal X}_{\eta^{-1}}$}
\put(0,70){$X_{\eta w}\bigcap X'_{w^\vee}$}
\put(70,55){$\Gamma|_{{\mathcal X}_{\eta^{-1}}}$}
\put(50,125){$\left(\wp^*\Gamma\right)\!|_{X_{\eta w}\bigcap X'_{w^\vee}}$}
\put(215,55){${\mathcal Y}_\zeta$}
\put(200,125){$Y'_u\bigcap Y_{\omega_0\zeta u}$}
\put(69,47){\vector(-1,-1){34}} 
\put(65,117){\vector(-1,-1){32}} 
\put(108,58){\vector(1,0){100}} 
\put(135,128){\line(1,0){5}}   \put(145,128){\line(1,0){5}} 
\put(155,128){\line(1,0){5}}   \put(165,128){\line(1,0){5}} 
\put(175,128){\line(1,0){5}}   \put(185,128){\vector(1,0){8}} 
\put(28,62){\vector(0,-1){49}} 
\put(80,117){\vector(0,-1){45}} 
\put(220,117){\vector(0,-1){45}} 
\put(18,40){$\wp$}  \put(82,95){$\wp$}   \put(225,95){$\pi$}  
\put(47,38){$g$}    \put(43,105){$g$}   
\put(150,63){$f$}    \put(160,132){$h$} 
\end{picture}   
$$
where the maps $h,f,\wp$, and $\pi$ are isomorphisms on Zariski dense sets.
\end{lem}

Since $g^{-1}({\mathcal X}_{\eta^{-1}})=\Gamma|_{{\mathcal X}_{\eta^{-1}}}$
and ${\mathcal Y}_\zeta$ have the same dimension, and the maps are
generically 1-1, 
$f(g^{-1}({\mathcal X}_{\eta^{-1}}))={\mathcal Y}_\zeta$,
which proves Lemma~\ref{lem:geom-pull-push}.
\smallskip

\noindent{\bf Proof. }
We first define an injective map 
$h:(\wp^*\Gamma)|_{X^\circ_{\eta w}}\rightarrow Y_{\omega_0\zeta u}$,
then show the restriction of $h$ to 
$\left(\wp^*\Gamma\right)\!|_{X_{\eta w}\bigcap X'_{w^\vee}}$
has image contained in $Y'_u\bigcap Y_{\omega_0\zeta u}$.
Since the maps $\pi$ and $\wp$ are generically 1-1, by
the analog of Theorem~\ref{thm:geometric_shapes} (Theorem 5.1.4
of~\cite{BS98a}) for ${\mathbb F}\ell(L)$,
the concluson follows.
\medskip

Let $\Edot\in{\mathbb F}\ell(L)$ be a complete flag and let
$(\Fdot,K)\in\wp^*\Gamma|_{X^\circ_{\eta w}\Edot}$, so that 
$\Fdot\in X^\circ_{\eta w}\Edot$ and $K\in\mbox{\it Lag}(V)$ satisfy
$$
 F_{k-1}\oplus F^\perp_k\ \subsetneq\ K\ \subsetneq\ F_k\oplus F_{k-1}^\perp.
$$
Define $h(\Fdot,K)\in\Spn/B$ by
$$
 h(\Fdot,K)_{\ol{\imath}}\ :=\ \left\{\begin{array}{lcl}
     F_{n+1-i}&\quad& n\geq i\geq n+2-k\\
     F_{k-1}+ F^\perp_{i+k-2}&& n+1-k\geq i\geq n+2-j\\
     K\bigcap (L+E_{n+\zeta u(i)})&& n+1-j\geq i
     \end{array}\right. .
$$
We show this defines a flag in $Y_{\omega_0\zeta u}\varphi_n\Edot$.
Note first that 
$h(\Fdot,K)_{\ol{\imath}}=\varphi_{k-1}(\Fdot)_{\ol{\imath}}$
for $i\geq n+2-j$.
Since $\epsilon_{k-1}(\eta w)(i)=\omega_0\zeta u(i)$ for $\geq n+2-j$,
Lemma~\ref{lem:hidden_subset} shows that 
$$ 	
   \dim(\varphi_n\Edot)_a\bigcap h(\Fdot,K)_{\ol{\imath}}\ \geq\ 
   \#\{l\geq i\mid a\geq \omega_0\eta u(i)\},
$$  	
for $i\geq n+2-j$.

We show this for $i\leq n+1-j$, which shows 
$h(\Fdot,K)\in Y_{\omega_0\zeta u}\varphi_n\Edot$.
Since $\zeta u$ is a Grassmannian permutation, and, for $a>0$
$(\varphi_n\Edot)_a=L+E^\perp_{n-a}$, we need only show that 
\begin{equation}\label{eq:condition}
   \dim K\bigcap (L+E^\perp_{n+\zeta u(i)})\  \geq\ n+1-i
\end{equation}
for $i\leq n+1-j$ (as $\zeta u(i)<0$ in this range).

It suffice to show this for the dense subset of $(\Fdot,K)$ with 
$K\bigcap L=F_{k-1}$.
Then $F_{k-1}^\perp$ is the image of $K$ under the projection
$V\twoheadrightarrow L$.
Thus, for $a>0$
\begin{eqnarray*}
\dim K\bigcap (L+E_{n-a}^\perp)
   &=& k-1 + \dim F_{k-1}^\perp\bigcap E_{n-a}^\perp \\
   &=& a + \#\{k-1\geq l\mid \eta w(l) \geq a+1\}
\end{eqnarray*}
If $k-1\geq l$, then $\eta w(l)=\zeta u(n+1-l)$ and these exhaust the
positive values of $\zeta u$.
If $i\leq n+1-j$ and we set $a=\ol{\zeta u(i)}$, we see that
$\dim K\bigcap (L+E^\perp_{n+\zeta u(i)})$ is
$\ol{\zeta u(i)} + \#\{l\mid \zeta u(l)\geq \ol{\zeta u(i)}+1\}$.
Since $\zeta u\in{\mathcal B}_n$ is a Grassmannian permutation,
$$
  \{\ol{\zeta u(i)},\ldots,n\}\ =\ 
  \{\ol{\zeta u(i)},\ldots,\ol{\zeta u(1)}\}\  \coprod\ 
  \{\zeta u(l)\mid \zeta u(l)>\ol{\zeta u(i)}\}.
$$
Thus 
$n+1-\ol{\zeta u(i)}= i+\#\{l\mid \zeta u(l)\geq \ol{\zeta u(i)}+1\}$
and so 
$\dim K\bigcap (L+E^\perp_{n+\zeta u(i)})=n+1-i$,
which completes the proof that 
$h(\Fdot,K)\in Y_{\omega_0\zeta u}\varphi_n\Edot$ for
$(\Fdot,K)\in g^{-1}X^\circ_{\eta w}\Edot$.
\medskip

We do a useful calculation 
before we finish the proof of Lemma~\ref{lem:big-diagram}.

\begin{lem}
If $\Fdot\in X^\circ_{\eta w}\Edot$ and $1\leq i\leq n+1-j$, then
$F_k^\perp\bigcap E_{n+\zeta u(i)}^\perp = F_{k+i-1}^\perp$.
\end{lem}

\noindent{\bf Proof. }
Similar to the last paragraph,
$\dim F^\perp_{k-1}\bigcap E_{n+\zeta u(i)}^\perp = n+1-i-k
   =\dim F_{k+i-1}^\perp$.
Thus, we need only show 
$F_{k+i-1}^\perp\subset F^\perp_{k-1}\bigcap E_{n+\zeta u(i)}^\perp$.

Note that 
$$
 \dim F_{k+i-1}\bigcap E_{n+\zeta u(i)}\ =\ 
 \#\{ k+i-1\geq l \mid \eta w(l) \geq \ol{\zeta u(i)}+1\}
$$
Since $\eta w(k+1)>\cdots> \eta w(n)$ and $\eta w(k+i)=\ol{\zeta u(i)}$, 
this equals
$$
\#\{ k+i>l\mid \eta w(l)>\eta w(k+i)\}\ =\ 
n-\eta w (k+i) = n+\zeta u(i),
$$
which shows $E_{n+\zeta u(i)}\subset F_{k+i-1}$.
Since $F_k\subset F_{k+i-1}$, this completes the proof.
\QED

We complete the proof of Lemma~\ref{lem:big-diagram},
showing that if we further require  $\Fdot\in X'_{w^\vee}$, then
$h(\Fdot,K)\in Y'_u$. 
Let $\Epdot$ be a flag opposite to $\Edot$ in $L$.
We show
$h(\Fdot,K)\in Y_u\varphi_0\Epdot $ for $(\Fdot,K) \in
\left(\wp^*\Gamma\right)\!|_{X_{\eta w}\bigcap X'_{w^\vee}}$.
Note first that $\Fdot$ satisfies
$$
  \dim E'_a\bigcap F_b\ =\ \#\{b\geq l\mid w(l)\leq a\}.
$$
As before
\begin{equation}\label{eq:last?}
  \dim(\varphi_0\Epdot)_a \bigcap h(\Fdot,K)_{\ol{\imath}}
    \ \geq \ \#\{i\leq l\mid a\geq u(l)\}
\end{equation}
for $i>n+1-j$ since $\epsilon_{k-1}(w^\vee)$ and $u$ agree in this range, as
do $h(\Fdot,K)_{\ol{\imath}}$ and $(\varphi_{k-1}\Fdot)_{\ol{\imath}}$.
For $i\leq n+1-j$, it suffices to establish~(\ref{eq:last?}) for the dense
subset of those $(\Fdot,K)$ with $K\bigcap L^\perp= F_k^\perp$.

Suppose $a>0$. 
Since $K\bigcap L^\perp= F_k^\perp$, $F_k$ is the image of $K$ under the
projection $V\twoheadrightarrow L$.
Similarly, $F_k$ is the image of 
$K\bigcap (L+E_{n+\ol{\zeta u(i)}}^\perp)=h(\Fdot,K)_{\ol{\imath}}$
under this projection.
Since the kernel of this is 
$F_k^\perp\bigcap E_{n+\zeta u(i)}^\perp = F_{k+i-1}^\perp$ and
$(\varphi_0\Epdot)_a=(L^\perp+E'_a)$,
\begin{eqnarray*}
  \dim (\varphi_0\Epdot)_a \bigcap h(\Fdot,K)_{\ol{\imath}}
  &=& \dim F_{k+i-1}^\perp + \dim F_k\bigcap E'_a,\\
  &=& n+1-k-i + \#\{ k\geq l\mid w(l)\leq a\},\\
  &=& n+1-i -\#\{k\geq l\mid w(l)>a\}.
\end{eqnarray*}
Since $\{w(1),\ldots,w(k)\}=\{u(n),\ldots,u(n+2-k),u(n+1-j)\}$
are the positive values of $u$ and $i\leq n+1-j$, this is 
$$
n+1-i-\#\{i\leq l\mid u(l)>a\}\ =\ \#\{ i\leq l\mid u(l)\leq a\},
$$
which shows~(\ref{eq:last?}).
\smallskip

Now suppose $a<0$.
Then $(\varphi\Epdot)_a=(E'_{\ol{a}-1})^\perp \subset L^\perp$
and so
\begin{eqnarray*}
 \left(\varphi_0\Epdot\right)_a \bigcap h(\Fdot,K)_{\ol{\imath}}
  &=& F_k^\perp \bigcap E_{n+\ol{\zeta u(i)}}^\perp \bigcap 
      (E'_{\ol{a}-1})^\perp\\
  &=& F_{k+i-1}^\perp\bigcap (E'_{\ol{a}-1})^\perp.
\end{eqnarray*}
This has dimension
$$
n-k-i+2-\ol{a} +\#\{k+i-1\geq l\mid w(l)\leq \ol{a}-1\}
\ =\ \#\{k+i\leq l\mid w(l)\geq \ol{a}\}.
$$
The values $\ol{w(l)}$ for $l\geq k+i$ are simply the negative values 
$u(l)$ for $l\geq i$.
Thus this is $\#\{i\leq l\mid u(l)\leq a\}$, which completes the proofs of
Lemma~\ref{lem:big-diagram} and the Pieri-type formula.
\QED

\end{document}